\documentclass[aos,preprint]{imsart}

\RequirePackage{amsthm,amsmath,amsfonts,amssymb}
\RequirePackage[authoryear]{natbib} 
\RequirePackage[colorlinks,citecolor=blue,urlcolor=blue]{hyperref}
\RequirePackage{graphicx}
\usepackage{algorithm}
\usepackage{caption}
\usepackage[caption=false]{subfig}
\usepackage[noend]{algpseudocode}
\usepackage{color} 
\usepackage{bm}
\usepackage{latexsym, epsfig, mathrsfs}
\usepackage[figuresright]{rotating}
\usepackage{rotating, tabularx,paralist}
\usepackage[OT1]{fontenc}
\usepackage{multirow,multicol, siunitx}
\usepackage[latin1]{inputenc}

\usepackage{bbm}
\usepackage{array}
\usepackage{tabu}
\usepackage{booktabs}
 
\makeatletter
 
\newtheorem*{assumption*}{Condition}
\newtheorem{lemma}{Lemma}
\newtheorem{proposition}{Proposition}
\newtheorem{theorem}{Theorem}
\newtheorem{example}{Example}
\newtheorem*{example*}{Example}

\renewcommand{\theequation}{
	\arabic{equation}%
}
\makeatother

\newcommand{\ignore}[1]{}{}

\newcommand{\var}{\mathrm{var}}

\newcommand{\cov}{\mathrm{cov}}

\newcommand{\diag}{\mathrm{diag}}

\def\independenT#1#2{\mathrel{\setbox0\hbox{$#1#2$}%
		\copy0\kern-\wd0\mkern4mu\box0}}

\theoremstyle{plain}

\renewcommand{\theenumi}{\roman{enumi}}

\newcommand{\R}{\mathbb R}
\newcommand{\e}{\mathbb E}
\newcommand{\V}{\mathcal V}
\newcommand{\Rc}{\mathcal R}
\newcommand{\Uc}{\mathcal U}
\newcommand*{\dif}{\mathop{}\!\mathrm{d}}
\newcommand{\noi}{\noindent}
\newcommand{\nn} {\nonumber}
\newcommand{\I}{\mathbf{1}}

\def\independenT#1#2{\mathrel{\rlap{$#1#2$}\mkern2mu{#1#2}}}

\begin{document}

\begin{frontmatter}
\title{Asymptotic Distributions of High-Dimensional Distance Correlation Inference}
\runtitle{HDCI}

\begin{aug}
\author[A]{\fnms{Lan} \snm{Gao}\ead[label=t2,mark]{gaolan@marshall.usc.edu}},
\author[A]{\fnms{Yingying} \snm{Fan}\ead[label=e2,mark]{fanyingy@marshall.usc.edu}},
\author[A]{\fnms{Jinchi} \snm{Lv}\ead[label=e3,mark]{jinchilv@marshall.usc.edu}}
\and
\author[B, C]{\fnms{Qi-Man} \snm{Shao}\ead[label=e4]{shaoqm@sustech.edu.cn}}
\address[A]{Data Sciences and Operations Department, Marshall School of Business, University of Southern California,\printead{t2,e2,e3}}
\address[B]{Department of Statistics and Data Science, Southern University of Science and Technology, \printead{e4}}
\address[C]{Department of Statistics, The Chinese University of Hong Kong}

\runauthor{L. Gao, Y. Fan, J. Lv and Q.-M. Shao}

\thankstext{}{Lan Gao is Postdoctoral Scholar, Data Sciences and Operations Department, Marshall School of Business, University of Southern California, Los Angeles, CA 90089. 
	Yingying Fan is Professor and Dean's Associate Professor in Business Administration, Data Sciences and Operations Department, Marshall School of Business, University of Southern California, Los Angeles, CA 90089. 
	Jinchi Lv is Kenneth King Stonier Chair in Business Administration and Professor, Data Sciences and Operations Department, Marshall School of Business, University of Southern California, Los Angeles, CA 90089. 
	Qi-Man Shao is Chair Professor, Department of Statistics and Data Science, Southern University of Science and Technology, Shenzhen, China and Choh-Ming Li Professor of Statistics, Department of Statistics, The Chinese University of Hong Kong, Shatin, N.T., Hong Kong. 
}
\end{aug}

\begin{abstract}
Distance correlation has become an increasingly popular tool for detecting the nonlinear dependence between a pair of potentially high-dimensional random vectors. Most existing works have explored its asymptotic distributions under the null hypothesis of independence between the two random vectors when only the sample size or the dimensionality diverges. Yet its asymptotic null distribution for the more realistic setting when both sample size and dimensionality diverge in the full range remains largely underdeveloped. In this paper, we fill such a gap and develop central limit theorems and associated rates of convergence for a rescaled test statistic based on the bias-corrected distance correlation in high dimensions under some mild regularity conditions and the null hypothesis. Our new theoretical results reveal an interesting phenomenon of blessing of dimensionality for high-dimensional distance correlation inference in the sense that the accuracy of normal approximation can increase with dimensionality. Moreover, we provide a general theory on the power analysis under the alternative hypothesis of dependence, and further justify the capability of the rescaled distance correlation in capturing the pure nonlinear dependency under moderately high dimensionality for a certain type of alternative hypothesis. The theoretical results and finite-sample performance of the rescaled statistic are illustrated with several simulation examples and a blockchain application.

\end{abstract}



\end{frontmatter}

\section{Introduction} \label{sec: intro} 
In many big data applications nowadays, we are often interested in measuring the level of association between a pair of potentially high-dimensional random vectors giving rise to a pair of large random matrices. There exist a wide spectrum of both linear and nonlinear dependency measures. Examples include the Pearson correlation \citep{P1895}, rank correlation coefficients \citep{K1938,S1904},  coefficients based on the cumulative distribution functions or density functions \citep{H1948,BKR1961,R1975}, measures based on the characteristic functions \citep{F1993,SRB2007,SR2009}, the kernel-based dependence measure \citep{GHSBS2005}, and sign covariances \citep{BD2014,WDM2018}. See also \cite{SP2018,BWBS2019} for some recent developments on determining the conditional dependency through the test of conditional independence. In particular, nonlinear dependency measures have been popularly used since independence can be fully characterized by zero measures. Indeed test of independence between two random vectors is of fundamental importance in these applications.


Among all the nonlinear dependency measures, distance correlation introduced in \cite{SRB2007} has gained growing popularity in recent years due to several appealing features. First, zero distance correlation completely characterizes the independence between two random vectors. Second, the pair of random vectors can be of possibly different dimensions and possibly different data types such as a mix of continuous and discrete components. Third, this nonparametric approach enjoys computationally fast implementation. In particular, distance-based nonlinear dependency measures have been applied to many high-dimensional problems. Such examples include dimension reduction \citep{VTE2018}, independent component analysis \citep{MT2017}, interaction detection 
\citep{KLFL2017}, feature screening 
\citep{LZZ2012,SZ2014}, and variable selection \citep{KWW2015,SZ2014}. See also the various extensions for testing the mutual independence \citep{YZS2018},  testing the multivariate mutual dependence \citep{JM2018,CZ2019}, testing the conditional mean and quantile independence \citep{ZYS2018}, the partial distance correlation \citep{SR2014}, the conditional distance correlation \citep{WPHTZ2015}, measuring  the nonlinear dependence in time series \citep{Z2012,DMMW2018}, and measuring the dependency between two stochastic processes \citep{MMS2017,DMMW2018}.

To exploit the distance correlation for nonparametric inference of test of independence between two random vectors $ X \in \R ^{p} $  and $ Y \in \R^{q}  $ with $p, q \geq 1$, it is crucial to determine the significance threshold. Although the bootstrap or permutation methods can be used to obtain the empirical significance threshold, such approaches can be computationally expensive for large-scale data. Thus it is appealing to obtain its asymptotic distributions for easy practical use. There have been some recent developments along this line. For example, for the case of fixed dimensionality with independent $ X $ and $ Y $, \cite{SRB2007} showed that the standardized sample distance covariance by directly plugging in the empirical characteristic functions converges in distribution to a weighted sum of chi-square random variables as the sample size $n$ tends to infinity. A bias-corrected version of the distance correlation was introduced later in \cite{SR2013, SR2014} to address the bias issue in high dimensions. 
\cite{HS2016} proved that for fixed dimensionality and independent $ X $ and $ Y$, the standardized unbiased sample distance covariance converges to a weighted sum of centralized chi-square random variables asymptotically. In contrast, \cite{SR2013} considered another scenario when the dimensionality diverges with sample size fixed and showed that for random vectors each with exchangeable components, the bias-corrected sample distance correlation converges to a suitable 
$t$-distribution. Recently \cite{ZYZS2019} extended the result to more general assumptions and obtained the central limit theorem in the high-dimensional medium-sample-size setting.  

{Despite the aforementioned existing results, the asymptotic theory for sample distance correlation between $X$ and $Y$ under the null hypothesis of independence in \textit{general} case of $n, p$ and $q$ diverging in an arbitrary fashion remains largely unexplored. As the first major contribution of the paper, we  provide a more complete picture of the precise limiting distribution in such setting.}
In particular, under some mild regularity conditions and the independence of $X$ and $Y$, we obtain central limit theorems for a rescaled test statistic based on the bias-corrected sample distance correlation in high dimensions (see Theorems \ref{thm1} and \ref{thm2}). Moreover, we derive the explicit rates of convergence to the limiting distributions (see Theorems \ref{thm3} and \ref{thm4}). To the best of our knowledge, the asymptotic theory built in Theorems \ref{thm1}--\ref{thm4} is new to the literature. 
{Our theory requires no constraint on the relationship between sample size $n$ and dimensionalities $p $ and $ q$. Our results show that the accuracy of normal approximation can increase with dimensionality, revealing an interesting phenomenon of blessing of dimensionality.}

{The second major contribution of our paper is to provide a general theory on the power analysis of the rescaled sample distance correlation. We show in Theorem \ref{thm-power} that as long as the population distance correlation and covariance do not decay too fast as sample size increases, the rescaled sample distance correlation diverges to infinity with asymptotic probability one, resulting in a test with asymtotic power one.  We further consider in Theorem \ref{prop-power} a specific alternative hypothesis where $X$ and $Y$ have pure nonlinear dependency in the sense that their componentwise Pearson correlations are all zero,   and show that the rescaled sample distance correlation achieves asymptotic power one when $p=q=o(\sqrt{n})$. This reveals an interesting message that in moderately high-dimensional setting, the rescaled sample distance correlation is capable of detecting pure nonlinear dependence with high power.   
}


Among the existing literature, the most closely related paper to ours is the one by \cite{ZYZS2019}. Yet, our results are significantly different from theirs. 
For clarity we discuss the differences under the null and alternative hypotheses separately.   Under the null hypothesis of $X$ and $Y$ being independent, our results differ from theirs in four important aspects: 
1) \cite{ZYZS2019} considered the scenario where sample size $ n $ grows at a slower rate compared to dimensionalities $ p $ and $q $, while our results make no assumption on the relationship between $ n $ and $ p, q $;  2) \cite{ZYZS2019} assumed that 
$ \min\{p, q\} \rightarrow \infty$, 
whereas our theory relies on a more relaxed assumption of $ p + q \rightarrow \infty $; 3) there is no rate of convergence provided in the work of \cite{ZYZS2019}, while explicit rates of convergence are developed in our theory; 4) the proof in \cite{ZYZS2019} is based on the componentwise analysis, whereas our technical proof is based on the joint analysis by treating the high-dimensional random vectors as a whole; See Table \ref{table-dif} in Section \ref{corollaries} for a summary of these key differences under the illustrative example of $m$-dependent components.

{
	The difference under the alternative hypothesis of dependence is even more interesting.   
	\cite{ZYZS2019} showed that under the alternative hypothesis of dependence, when both dimensionalities $ p  $ and $q$ grow much faster than sample size $n$, the sample distance covariance asymptotically measures the linear dependence between two random vectors satisfying certain moment conditions, and fails to capture the nonlinear dependence in high dimensions. 
	To address this issue, a marginally aggregated distance correlation statistic was introduced therein to deal with high-dimensional independence testing.   However, as discussed above, we provide a specific alternative hypothesis under which the rescaled sample distance correlation is capable of identifying the pure nonlinear relationship when $p=q=o(\sqrt{n})$. These two results complement each other and indicate that the sample distance correlation can have rich asymptotic behavior in different diverging regimes of $(n, p, q)$. The complete spectrum of the alternative distribution as a function of $(n,p,q)$ is still largely open and can be challenging to study. In simulation Example \ref{ex-power3} in Section \ref{Sec4.3}, we give an example showing that the marginally aggregated distance correlation statistic can suffer from power loss if the true dependence in data is much more than just marginal.} 

It is also worth mentioning that our Propositions \ref{prop1}--\ref{prop3} (see Section \ref{pr-thm4} of Supplementary Material), which serve as the  crucial ingredient of the proofs for Theorems \ref{thm2} and \ref{thm4}, provide some explicit bounds on certain key moments identified in our theory under fairly general conditions, which can be of independent interest.

The rest of the paper is organized as follows. Section \ref{background} introduces 
the distance correlation and reviews the existing limiting distributions. We present a rescaled test statistic, its asymptotic distributions, and a power analysis for high-dimensional distance correlation inference in Section \ref{sec:result}. Sections \ref{simulation} and \ref{Sec5} provide several simulation examples and a blockchain application justifying our theoretical results and
illustrating the finite-sample performance of the rescaled test statistic. We discuss some implications and extensions of our work in Section \ref{discussion}.  
All the proofs and technical details are provided in the Supplementary Material.

\section{Distance correlation and distributional properties} \label{background}


\subsection{Bias-corrected distance correlation} \label{Sec2.1}
Let us consider a pair of random vectors $ X   \in \R^p $ and $ Y   \in \R^q $ with integers $ p, q \geq 1$ that are of possibly different dimensions and possibly mixed data types such as continuous or discrete components. For any vectors $ t \in \R^p $ and $ s \in  \R^q $, denote by $ \langle t, X \rangle  $ and $ \langle s, Y \rangle   $  the corresponding inner products. 
Let $ f_X (t) = \e e^{ i \langle t, X \rangle }, f_Y (s) = \e e^{ i \langle s, Y \rangle } $, and $ f_{X, Y} (t, s) = \e e^{ i \langle t, X \rangle + i \langle s, Y \rangle } $ be the characteristic functions of $ X $, $ Y $, and the joint distribution $ (X, Y) $, respectively, where $i$ associated with the expectations represents the imaginary unit $(-1)^{1/2}$. \citet{SRB2007} defined the squared distance covariance $ \V^2 ( X, Y )  $ as
\begin{equation}
\V^2 (X, Y) = \int_{\R^{ p + q} }  \frac {| f_{X, Y} (t,s)  - f_X (t) f_Y (s) |^2 }  { c_p c_q  \Vert t \Vert ^{p + 1}  \Vert s \Vert ^{q + 1} }  \dif t \dif s,  \label{DisCov}
\end{equation} 
where
$$ 
c_p = \frac { \pi^{ ( p+ 1)/2 } } { \Gamma(  ( p + 1 ) / 2  ) }
$$ 
with $ \Gamma(\cdot) $ the 
gamma function and $ \Vert \cdot \Vert   $ stands for the Euclidean norm of a vector. Observe that $2 c_p$ and $2 c_q$ are simply the volumes of $p$-dimensional and $q$-dimensional unit spheres in the Euclidean spaces, respectively. In view of the above definition, it is easy to see that $ X  $ and $Y $ are independent if and only if $ \V^2 (X , Y)  = 0 $. Thus distance covariance characterizes completely  the independence.

The specific weight in \eqref{DisCov} gives us an explicit form of the squared distance covariance (see \cite{SRB2007})
\begin{align}
\V^2 ( X, Y)    = & \e [ \Vert X_1 - X_2  \Vert  \Vert  Y_1 - Y_2 \Vert   ]    -  2  \e [ \Vert  X_1 - X_2 \Vert  \Vert  Y_1 - Y_3 \Vert ]   \nn \\
& +  \e [ \Vert  X_1 - X_2 \Vert ]  \e [ \Vert Y_1 - Y_2 \Vert ],  \label{dcov1}
\end{align}
where $ (X_1, Y_1)$ , $( X_2, Y_2 )$,  and $(X_3, Y_3) $ are independent copies of $ ( X, Y ) $. Moreover, \citet{Lyons2013} showed that
\begin{equation}
\V^2 (X, Y) = \e [ d( X_1, X_2 ) d (Y_1, Y_2) ]  \label{dcov2}
\end{equation}
with the double-centered distance
\begin{align}
d ( X_1, X_2 ) = \Vert X_1 - X_2 \Vert - \e [ \Vert X_1 - X_2 \Vert | X_1  ] - \e [ \Vert X_1 - X_2 \Vert | X_2  ]  + \e [ \Vert X_1 - X_2 \Vert ] \label{d-func}
\end{align}
and $ d( Y_1, Y_2 ) $ defined similarly. Let $ \V ^2 ( X  )  = \V^2 ( X, X ) $ and $ \V^2 ( Y ) = \V ^2 ( Y, Y ) $ be the squared distance variances of $ X $ and $ Y $, respectively. Then the squared distance correlation $ \Rc ( X, Y ) $ is defined as 
\begin{equation} \label{neweq.001}
\Rc^2 (X, Y) = \left\{  
\begin{aligned}
&\frac { \V ^2 (X, Y)} { \sqrt{ \V^2 ( X )   \V^2 ( Y ) } } & \ \text{ if } \V^2 ( X ) \V^2 ( Y ) > 0 , \\
&  0     & \ \text{ if } \V^2 ( X ) \V^2 ( Y )  = 0.
\end{aligned}  
\right.
\end{equation}

Now assume that we are given a sample of $ n $ independent and identically distributed (i.i.d.) observations $ \{ (X_i, Y_i), 1 \leq i \leq n \}$ from the joint distribution $ ( X, Y ) $. In \citet{SRB2007}, the squared sample distance covariance $ \V_n^2 (X, Y) $ was constructed by directly plugging in the  empirical characteristic functions as 
\begin{align}
\V_n^2 ( X, Y )  =  \int_{\R^{ p + q} }  \frac {| f_{X, Y}^n  (t,s)  - f_X^n (t) f_Y^n  (s) |^2 }  { c_p c_q  \Vert t \Vert ^{p + 1}  \Vert s \Vert ^{q + 1} }  \dif t \dif s, \label{sam-dcov}
\end{align}
where $ f_Y^n ( t )$, $f_Y^n ( s ) $, and $ f_{ X, Y }^n ( t, s )  $ are the corresponding empirical characteristic functions.  
Thus the squared sample distance correlation is given by 
\begin{align} \label{neweq.002}
\Rc_n^2 ( X, Y  ) = \left\{  
\begin{aligned}
&\frac { \V_n ^2 (X, Y)} { \sqrt{ \V_n^2 ( X )   \V_n^2 ( Y ) } } & \ \text{ if }  \V_n^2 ( X ) \V_n^2 ( Y ) > 0 , \\
&  0     & \ \text{ if } \V_n^2 ( X ) \V_n^2 ( Y )  = 0.
\end{aligned}  
\right.
\end{align}
Similar to \eqref{dcov1} and \eqref{dcov2}, the squared sample distance covariance admits the following explicit form 
\begin{align} \label{neweq.003}
\V_n^2 ( X, Y )  =  \frac {1} { n^2 } \sum_{ k, l = 1 }^n  A_{k, l} B_{k, l},
\end{align} 
where $ A_{k, l} $ and $ B_{ k, l } $ are the double-centered distances defined as 
\begin{align*}
A_{k, l}  & = a_{k, l} - \frac 1 n  \sum_{i = 1}^n   a_{ i, l  }  -  \frac 1 n  \sum_{j = 1}^n   a_{ k, j } + \frac { 1 } { n^2 } \sum_{ i, j = 1 }^n  a_{i, j} , \\
B_{k, l}  & = b_{k, l} - \frac 1 n  \sum_{i = 1}^n   b_{ i, l  }  -  \frac 1 n  \sum_{j = 1}^n   b_{ k, j } + \frac { 1 } { n^2 } \sum_{ i, j = 1 }^n  b_{i, j}
\end{align*}
with $ a_{ k, l }  = \Vert X_k - X_l \Vert $ and $ b_{k , l}  = \Vert Y_k - Y_l  \Vert  $. It is easy to see that the above estimator is an empirical version of the right hand side of \eqref{dcov2}. The double-centered population distance $ d ( X_k, X_l ) $ is estimated by the double-centered sample distance $ A_{k, l} $ and then $ \e [ d ( X_1, X_2 ) ] $ is estimated by the mean of all the pairs of double-centered sample distances. 

Although it is natural to define the sample distance covariance in \eqref{sam-dcov}, \citet{SR2013} later demonstrated that such an estimator is biased and can lead to interpretation issues in high dimensions. They revealed that for independent random vectors $ X \in \R^p $ and $ Y \in \R^q $ with i.i.d. components and finite second moments, it holds that 
\begin{align*}
\Rc_n^2 ( X, Y )  \xrightarrow[ p, q \rightarrow \infty ] {}  1
\end{align*}
when sample size $ n $ is fixed, but we naturally have $ \Rc^2 ( X, Y )= 0  $ in this scenario. To address this issue, \citet{SR2013, SR2014} introduced a modified unbiased estimator of the squared distance covariance and the bias-corrected sample distance correlation given by 
\begin{equation} \label{est-DisCov}
\V^*_n (X, Y) = \frac {1} { n  (n - 3) } \sum_{ k \neq l } A_{k, l}^* B_{k, l}^*     
\end{equation}
and 
\begin{align} \label{neweq.004}
\Rc_n^* ( X, Y  ) = \left\{  
\begin{aligned}
&\frac { \V_n ^* (X, Y)} { \sqrt{ \V_n^* ( X )   \V_n^* ( Y ) } } & \ \text{ if } \V_n^* ( X ) \V_n^* ( Y ) > 0 , \\
&  0     & \ \text{ if }  \V_n^* ( X ) \V_n^* ( Y )  = 0,
\end{aligned}  
\right.
\end{align}
respectively, where the $\Uc$-centered distances $ A_{k, l}^* $ and $ B_{k, l}^* $ are defined as
\begin{align*}
A_{k, l}^*  & = a_{k , l}  - \frac {1} { n - 2 }  \sum_{i = 1}^n   a_{i, l } -  \frac { 1 } { n - 2 }  \sum_{j = 1}^n   a_{ k, j } + \frac { 1 } { ( n - 1 ) ( n - 2 ) } \sum_{ i, j = 1 }^n  a_{ i, j  }, \\
B_{k, l}^*  & = b_{k , l}  - \frac {1} { n - 2 }  \sum_{i = 1}^n   b_{i, l } -  \frac { 1 } { n - 2 }  \sum_{j = 1}^n   b_{ k, j } + \frac { 1 } { ( n - 1 ) ( n - 2 ) } \sum_{ i, j = 1 }^n  b_{ i, j  }.
\end{align*}
Our work will focus on the bias-corrected distance-based statistics $ \V_n^*  ( X, Y ) $ and $ \Rc_n^*  ( X, Y ) $ given in (\ref{est-DisCov}) and (\ref{neweq.004}), respectively.  

\subsection{Distributional properties} \label{Sec2.2}

In general, the exact distributions of the  distance covariance and distance correlation are intractable. Thus it is essential to investigate the asymptotic surrogates in order to apply the distance-based statistics for the test of independence. 
With dimensionalities $ p, q $ fixed and sample size $ n \rightarrow \infty $, \citet{HS2016} validated that $ \V^*_n (X, Y) $ is a U-statistic and then under the independence of $ X $ and $ Y $, it admits the following asymptotic distribution
\begin{equation}
n \V^*_n (X, Y)  \xrightarrow  [ n\rightarrow \infty ]  {\mathscr{D}} \sum_{ i = 1 }^{ \infty } \lambda_i  ( Z_i^2 - 1 ), \label{finite-p}
\end{equation}
where $ \{ Z_i, i \geq 1  \} $ are i.i.d. standard normal random variables and $  \{ \lambda_i , i \geq 1 \}  $  are the eigenvalues of some operator. 

On the other hand, \cite{SR2013} showed that when the dimensionalities $ p $ and $ q $ tend to infinity and sample size $ n \geq 4 $ is fixed, if $ X  $ and $ Y $ both consist of i.i.d. components, then under the independence of $ X $ and $ Y $ we have 
\begin{equation}  
T_R : = \sqrt{ n (n - 3)/2 - 1} \frac { \Rc^*_n (X, Y) } { \sqrt{1 - ( \Rc^*_n (X, Y) )^2 } } \xrightarrow  [ p, \, q \rightarrow \infty ]  {\mathscr{D}}t_{n (n - 3)/2 - 1}. \label{fix-p}
\end{equation}
However, it still remains to investigate the limiting distributions of distance correlation when both sample size and dimensionality are diverging simultaneously. It is common to encounter datasets that are of both high dimensions and large sample size such as in biology, ecology, medical science, and networks. 
When $ p \land q \rightarrow \infty$ and $ n \rightarrow \infty $ at a slower rate compared to $ p, q $, under the independence of $ X $ and $ Y $ and some conditions on the moments \cite{ZYZS2019} showed that 
\begin{align}
T_R \xrightarrow     {\mathscr{D}}  N ( 0, 1 ), \label{high-dim}
\end{align}
where  $ p \land q  $ denotes the minimum value of $p$ and $q$. Their result was obtained by approximating the unbiased sample distance covariance with the aggregated marginal distance covariance, which can incur stronger assumptions including $ n \rightarrow \infty $ at a slower rate compared to $ p, q $ and $ p \land q  \rightarrow \infty$.

The main goal of our paper is to fill such a gap and make the asymptotic theory of distance correlation more complete. Specifically, we will prove central limit theorems for $  \Rc_n^* (X, Y)  $ when $ n \rightarrow \infty $ and $  p + q   \rightarrow \infty $. In contrast to the work of \cite{ZYZS2019}, we analyze the unbiased sample distance covariance directly by treating the random vectors as a whole. {Our work will also complement the recent power analysis in \cite{ZYZS2019}, where distance correlation was shown to asymptotically measure only linear dependency in the regime of fast growing dimensionality ($ \min \{ p, q    \} / n^2 \to \infty $) and thus the marginally aggregated distance correlation statistic was introduced. However, as shown in 
	Example \ref{ex-power3} in Section \ref{Sec4.3}, the marginally aggregated statistic can be less powerful than the  joint distance correlation statistic when the dependency between the two random vectors far exceeds the marginal contributions. To understand such a phenomenon, we will develop a general theory on the power analysis for the rescaled distance correlation statistic in Theorem \ref{thm-power} and further justify its capability of detecting nonlinear dependency in Theorem \ref{prop-power} for the regime of moderately high dimensionality.
}

\section{High-dimensional distance correlation inference} \label{sec:result}

\subsection{A rescaled test statistic} \label{Sec3.1}

To simplify the technical presentation, we assume that $ \e [ X ] =  0 $ and $ \e [ Y ] = 0 $ since otherwise we can first subtract the means in our technical analysis. Let $ \e [ X X^T  ] = \Sigma_x $  and $\e [ Y Y^T  ] = \Sigma_y $ be the covariance matrices of random vectors $X$ and $Y$, respectively. To test the null hypothesis that $X$ and $Y$ are independent, in this paper we consider a rescaled test statistic defined as a rescaled distance correlation 
\begin{align}
T_n :=  \sqrt{ \frac { n (n - 1) } { 2   } }  \Rc^*_n (X, Y) =  \sqrt{  \frac { n (n - 1) } { 2   }  }    \frac { \V^*_n  (X, Y)} { \sqrt{ \V^*_n  (X )   \V^*_n  (Y ) } }.   \label{Tn}
\end{align}
It has been shown in \cite{HS2016} that $ \V_n^* (X, Y) $ is a U-statistic. A key observation is that by the Hoeffding decomposition for U-statistics, the dominating part is a martingale array under the  independence of $ X $ and $ Y $. Then we can apply the martingale central limit theorem and calculate the specific moments involved. 

More specifically, \cite{HS2016} showed that 
\begin{align}
\V^*_n (X, Y) =  { n \choose 4 } ^{-1} \sum_{1 \leq i_1 < i_2 < i_3 < i_4 \leq n} h ( ( X_{i_1}, Y_{i_1} ), \cdots, ( X_{i_4}, Y_{i_4} )  ),  \label{U-stat}
\end{align}
where the kernel function is given by 
\begin{align}
& h ( ( X_{1}, Y_{1} ),  ( X_{2}, Y_{2} ), ( X_{3}, Y_{3} ), ( X_{4}, Y_{4} )  )    \nn \\
& =  \frac 14 \sum_{ \substack {  1 \leq i, j \leq 4,  \\ i \neq j  } }   \Vert X_i - X_j \Vert  \Vert Y_i -  Y_j \Vert    - \frac 1 4 \sum_{i = 1}^4  \bigg(  \sum_{ \substack {  1 \leq  j \leq 4, \\  j \neq i } }    \Vert X_i - X_j \Vert   \sum_{ \substack {  1 \leq  j \leq 4, \\  j \neq i  } } \Vert Y_i - Y_j \Vert   \bigg) \nn \\
&   \quad + \frac {1} {24} \sum_{ \substack {  1 \leq i, j \leq 4,  \\ i \neq j } } \Vert X_i - X_j \Vert   \sum_{ \substack {  1 \leq i, j \leq 4,  \\ i \neq j  } }   \Vert Y_i - Y_j \Vert . \label{kernel}
\end{align}
Let us define another functional 
\begin{equation} 
g(X_1, X_2, X_3, X_4) := d(X_1, X_2) d( X_1, X_3 ) d( X_2, X_4 ) d(X_3, X_4 ),  \label{def-g}
\end{equation} 
where $ d ( \cdot, \cdot ) $ is the double-centered distance defined in \eqref{d-func}. The above technical preparation enables us to derive the main theoretical results.

\subsection{Asymptotic distributions} \label{Sec3.2}

\begin{theorem} \label{thm1}
	Assume that $ \e \Vert X \Vert ^{ 2 + 2 \tau }  + \e \Vert Y \Vert ^{ 2 + 2 \tau } < \infty $ for some constant $ 0 < \tau \leq 1 $.  If 
	\begin{equation}
	\frac { \e ( | d (X_1, X_2 ) |^{  2 + 2 \tau } )  \e ( | d  (Y_1, Y_2 ) |^{  2 + 2 \tau } ) } { n^{ \tau }  [ \V^2 ( X )  \V^2 ( Y ) ]^{ 1 + \tau   }  } \rightarrow 0  \label{cond1}
	\end{equation}  
	and
	\begin{align}
	\frac { \e [ g(X_1 , X_2 , X_3 , X_4) ] \e [ g( Y_1, Y_2 , Y_3 , Y_4 ) ] } {  [ \V^2 ( X )  \V^2 ( Y ) ]^2 } \rightarrow 0   \label{cond2}
	\end{align}
	as $ n \rightarrow \infty $ and $ p + q  \rightarrow \infty $,
	then under the independence of $ X $ and $ Y $ we have
	$
	T_n \stackrel{\mathscr{D}}{\rightarrow} N  (0,1).
	$
\end{theorem}

Theorem \ref{thm1} presents a general theory and relies on the martingale central limit theorem. In fact, when $ X $ and $ Y $ are independent, via the Hoeffding decomposition we can find that the dominating part of $ \V_n^* (X, Y) $ forms a martingale array which admits asymptotic normality under conditions \eqref{cond1} and \eqref{cond2}. Moreover, it also follows from \eqref{cond1} that
\begin{align*}
\frac { \V_n^* ( X ) } { \V^2 (X) } \rightarrow 1  \quad \mbox{and}\quad  \frac { \V_n^* (Y)  } { \V^2 (Y)  } \rightarrow 1 \quad \mbox{in~probability.} 
\end{align*}
Thus an application of Slutsky's lemma results in the desired results.

Although Theorem \ref{thm1} is for the general case, the calculation of the moments involved such as $ \e [ g ( X_1, X_2, X_3, X_4 ) ] $, $ \V^2 ( X )$, and $\e ( | d ( X_1, X_2 ) |^{ 2 + 2 \tau } ) $ for the general underlying distribution can be challenging. To this end, we provide in Propositions \ref{prop1}--\ref{prop3} in Section \ref{pr-thm4} some bounds or exact orders of those moments. These results together with Theorem \ref{thm1} enable us to obtain Theorem \ref{thm2} on an explicit and useful central limit theorem with more specific conditions. 
Let us define quantities
\begin{gather*}
B_{X} =  \e [ \Vert X_1 - X_2 \Vert ^2 ] = 2 \e [ \Vert X \Vert^2 ] , \quad  B_{Y} =  \e [ \Vert Y_1 - Y_2 \Vert ^2 ] = 2 \e [ \Vert Y \Vert^2 ] ,  \\
L_{x, \tau} =   \e  \big(  \big| \Vert X \Vert^2 - \e \Vert X \Vert^{ 2} \big|^{ 2 + 2 \tau  } \big) + \e ( | X_1^T X_2 |^{ 2 + 2 \tau }  )  ,  \\
L_{y, \tau} =   \e  \big(  \big| \Vert Y \Vert^2 - \e \Vert Y \Vert^{ 2 } \big|^{ 2 + 2 \tau  } \big) + \e ( | Y_1^T Y_2 |^{ 2 + 2 \tau }  ) ,
\end{gather*}
and
\begin{align*}
E_{ x } & =   \frac { \e [ ( X_1^T \Sigma_x X_2 )^2 ]  + B_X^{ - 2 \tau  } L_{x, \tau}^{ (2 + \tau) / ( 1 + \tau  ) }  } {  (  \e [ ( X_1^T X_2 )^2 ]  )^2 }  , \\
E_{  y } & =  \frac { \e [ ( Y_1^T \Sigma_y Y_2 )^2 ]  + B_Y^{ - 2 \tau  } L_{y, \tau}^{ ( 2 + \tau ) / ( 1 + \tau  )  }  } {  (  \e [ ( Y_1^T Y_2 )^2 ]  )^2 }.   
\end{align*}

\begin{theorem}  \label{thm2}
	Assume that $ \e  [ \Vert X \Vert ^{ 4 + 4 \tau }   ] + \e [  \Vert Y \Vert ^{ 4 + 4 \tau }   ] < \infty  $ for some constant $ 0 <  \tau \leq 1 / 2  $ and as $ n \rightarrow \infty $ and $ p + q   \rightarrow \infty $, 
	\begin{align}
	\frac {  n^{ - \tau  }  L_{x, \tau}  L_{ y, \tau } }  {   \big(  \e [ ( X_1^T  X_2 )^2 ]  \e [ ( Y_1^T Y_2 )^2 ] \big)^{ 1 + \tau }  }     \rightarrow 0.  \label{co1}     
	\end{align}
	In addition, assume that $ E_{x} \rightarrow 0 $ if $ p \rightarrow  \infty $, and $ E_{y} \rightarrow 0 $ if $ q \rightarrow \infty $.
	Then under the independence of $ X $ and $ Y$, we have 
	$
	T_n \stackrel{\mathscr{D}}{\rightarrow} N(0, 1) .
	$
\end{theorem}

Theorem \ref{thm2} provides a user-friendly central limit theorem with mild regularity conditions that are easy to verify and can be satisfied by a large class of distributions. To get some insights into the orders of the moments $ B_X $, $ L_{x, \tau} $, $ \e [ ( X_1^T X_2 )^2 ] $, and $ \e [ ( X_1^T \Sigma_x X_2)^2 ] $, one can refer to Section \ref{corollaries} for detailed explanations by examining some specific examples. In Theorem \ref{thm2}, we show the results only under the scenario of $ 0 < \tau \leq 1 /2  $. In fact, similar results also hold for the case of $ 1/2 < \tau \leq 1 $; see Section \ref{tau} of Supplementary Material for more details.

\subsection{Rates of convergence} \label{rate}
Thanks to the martingale structure of the dominating term of $ \V_n^* ( X, Y ) $ under the independence of $ X $ and $ Y$, we can obtain explicitly the rates of convergence for the normal approximation.

\begin{theorem}    \label{thm3}
	Assume that $ \e \Vert X \Vert ^{ 2 + 2 \tau }  + \e \Vert Y \Vert ^{ 2 + 2 \tau } < \infty $ for some constant $ 0 < \tau \leq 1 $. Then under the independence of $X$ and $Y$, we have 
	\begin{align}
	\sup \limits_{ x \in \mathbb{R} }  | \mathbb{P} ( T_n  \leq  x  )  - \Phi ( x )  |  
	& \leq  C  \bigg\{   \Big ( \frac  {      \e [ g (X_1,X_2, X_3, X_4 ) ]  \e [ g( Y_1, Y_2, Y_3, Y_4 ) ]   }  {  [  \V^2 ( X ) \V^2 ( Y )  ] ^2 }  \Big)^{ \frac { 1 + \tau } { 2 } } \nn \\
	& \hspace{1.5cm}   + \frac {  \e [ | d (X_1, X_2 ) |^{2 + 2 \tau } ]  \e [  | d (Y_1, Y_2) |^{2 + 2 \tau } ] } { n^{ \tau  }  [ \V^2 ( X ) \V^2 ( Y ) ]^{ 1 + \tau  }    } \bigg\}^{  \frac { 1 } { 3 + 2 \tau  } },  \label{rate-general} 
	\end{align}
	where $C$ is some positive constant and $ \Phi (x) $ is standard normal distribution function.
\end{theorem}

In view of the evaluation of the moments in Propositions \ref{prop1}--\ref{prop3}, we can obtain the following theorem as a consequence of Theorem \ref{thm3}.

\begin{theorem} \label{thm4}
	Assume that $ \e [ \Vert X \Vert^ { 4 + 4 \tau }  ]  + \e [ \Vert Y \Vert^{ 4 + 4 \tau } ] < \infty $ for some constant $ 0 < \tau \leq 1 /2 $,
	\begin{align}
	B_X^{ - 2 \tau }  L_{ x, \tau } /  \e [ ( X_1^T X_2 )^2 ]     \leq 1 / 18, \ \text{ and }  \    B_Y^{ - 2 \tau }  L_{ y, \tau }  /  \e [ ( Y_1^T Y_2 )^2 ]   \leq 1 / 18.     \label{co2}
	\end{align}
	Then under the independence of $X$ and $Y$, we have 
	\begin{align}
	\sup\limits_{ x \in \R } | \mathbb{P} ( T_n \leq x )  - \Phi ( x ) |  
	& \leq C  \bigg\{ (  E_x  E_y ) ^{  \frac { 1 + \tau } { 2 } }   +      \frac{ n^{ - \tau  }  L_{ x, \tau }  L_{ y, \tau } } { \big(  \e [ ( X_1^T X_2 )^2 ] \e [ (  Y_1^T Y_2 )^2 ]  \big)^{  1 + \tau }  }  \bigg\}^{ \frac { 1 } { 3 + 2 \tau } },  \label{rate1}
	\end{align}
	where $C$ is some positive constant.
	
\end{theorem}

The counterpart theory for the case of $ 1/2 < \tau \leq 1  $ is presented in Section \ref{tau} of Supplementary Material. In general, larger value of $ \tau $ will lead to better convergence rates and weaker conditions, which will be elucidated by the example of $m$-dependent components in Proposition \ref{cor1} (see Section \ref{corollaries}).

Let us now consider the case when only one of $ p $ and $q $ is diverging, say, $ p $ is fixed and $ q \rightarrow \infty $. Then by the moment assumption $ \e [ \Vert X \Vert ^{4 + 4 \tau}  ] < \infty $, all the moments related to $X$ on the right hand side of \eqref{rate-general} are of bounded values. Thus in light of the proof of Theorem \ref{thm4}, we can see that if $ \e [ \Vert X \Vert ^{ 4 + 4 \tau  }  ]  + \e [ \Vert Y \Vert ^{ 4 + 4 \tau  } ]  < \infty $ for some constant $ 0 < \tau \leq 1/2 $, then there exists some positive constant $ C_X $ depending on the underlying distribution of $ X $ such that under the independence of $ X $ and $ Y $, we have
\begin{align}
\sup \limits_{ x \in \R }  | \mathbb{P} ( T_n < x  )  - \Phi (x)  |  \leq C_X \bigg\{ \Big( E_y \land \frac {1} { 18 } \Big)^{ \frac { 1 + \tau } { 2 } } +   \frac{ n^{ - \tau  }   L_{ y, \tau } } { \big(   \e [ (  Y_1^T Y_2 )^2 ]  \big)^{  1 + \tau }  }    \bigg \} ^{ \frac { 1 } { 3 + 2 \tau } }.  \label{rate-single}
\end{align}

It is worth mentioning that the bounds obtained in \eqref{rate-general} and \eqref{rate1} are nonasymptotic results that quantify the accuracy of the normal approximation and reveal how the rate of convergence depends on the sample size and dimensionalities. Since we exploit the rate of convergence in the central limit theorem for general martingales \citep{H1988} under the assumption of $ 0 < \tau \leq 1  $, the result may not necessarily be optimal. It is possible that better convergence rate can be obtained for the case of $ \tau > 1 $, which is beyond the scope of the current paper.

An anonymous referee asked a great question on whether similar results as in Theorems \ref{thm1} and \ref{thm3} apply to the studentized statistic $T_R$ defined in \eqref{fix-p}. 
The answer is affirmative.  Combining our Theorem \ref{thm1} with Lemma \ref{le-consist} {
	and \eqref{h-b1}}, it can be shown that {
	$T_R$} enjoys the same asymptotic normality as $T_n$ presented in Theorem \ref{thm1}. Moreover, the rates of convergence in Theorem \ref{thm3} also apply to $T_R$. See Section \ref{SecF} of Supplementary Material for the proof of these results for $T_R$. These results suggest that the studentized statistic $T_R$ can be a good choice in both small and large samples. Yet  the exact phase transition theory for the asymptotic null distribution of $T_R$ in the full diverging spectrum of $(n,p,q)$ remains to be developed.

\subsection{Some specific examples}  \label{corollaries}
To better illustrate the results obtained in the previous theorems, let us consider several concrete examples now. To simplify the technical presentation, we assume in this section that both $ p  $ and $ q  $ tend to infinity as $n$ increases. Our technical analysis also applies to the case when only one of $ p  $ and $ q $ diverges.

\begin{proposition}   \label{cor0}
	Assume that $ \e  ( \Vert X \Vert ^{ 4 + 4 \tau }   ) + \e (  \Vert Y \Vert ^{ 4 + 4 \tau }   ) < \infty  $ for some constant $ 0 <  \tau \leq 1 / 2  $ and there exist some positive constants $ c_1,  c_2 $ such that 
	\begin{gather}
	L_{ x, \tau } \leq c_1 p^{ 1 + \tau  } ,  \quad \e [ ( X_1^T \Sigma_x X_2 )^2 ] \leq c_1 p  , \label{cor0-c1}  \\ 
	\e [   ( X_1^T X_2 )^2  ] \geq c_2 p , \quad  \e [ \Vert X \Vert^2  ] \geq c_2 p ,  \label{cor0-c2}
	\end{gather}
	and  
	\begin{gather}
	L_{y, \tau} \leq c_1 q^{ 1 + \tau  } ,  \quad \e [ ( Y_1^T \Sigma_y Y_2 )^2 ] \leq c_1 q  ,  \label{cor0-c3} \\ 
	\e [   ( Y_1^T Y_2 )^2  ] \geq c_2 q, \quad  \e [ \Vert Y \Vert^2  ] \geq c_2 q. \label{cor0-c4}
	\end{gather}
	Then under the independence of $ X $ and $ Y $, there exists some positive constant A depending upon $ c_1 $ and $ c_2 $ such that for sufficiently large $ p  $ and $ q $, we have 
	\begin{align*}
	\sup\limits_{ x \in \R } | \mathbb{P} ( T_n \leq x )  - \Phi (x ) |   \leq  A  \big[ (   p q )^ {  - \tau ( 1 + \tau  ) / 2  } +   n^{ - \tau }  \big]^{ 1 / ( 3 + 2 \tau ) }.
	\end{align*}     
	Hence as $ n \rightarrow \infty $ and $ p , q  \rightarrow \infty $, it holds that 
	$
	T_n \stackrel{\mathscr{D}}{\rightarrow} N(0, 1). 
	$
\end{proposition}

The first example considered in Proposition \ref{cor0} is motivated by the case of independent components. Indeed, by Rosenthal's inequality for the sum of independent random variables, \eqref{cor0-c1} and \eqref{cor0-c2} are automatically satisfied when $ X $ consists of independent nondegenerate components with zero mean and uniformly bounded $ ( 4 + 4 \tau ) $th moment.

We next consider the second example of $m$-dependent components. For an integer $ m \geq 1  $, a sequence $ \{U_i\}_{0}^{\infty} $ is $m$-depenendent if $ \{ U_i \}_{0}^n  $ and $ \{U_i\}_{n + m + 1}^{\infty} $ are independent for every $ n \geq 0 $. We now focus on a special but commonly used scenario in which $ X $ consists of $ m_1 $-dependent components and $ Y $ consists of $m_2$-dependent components for some integers $ m_1\geq 1  $ and $ m_2 \geq 1$.
Assume that $ ( X_1, Y_1 )  $ and $ ( X_2, Y_2 ) $ are independent copies of $ ( X, Y ) $ and denote by 
\begin{align*}
X_1 = ( X_{1,1}, X_{1, 2}, \cdots, X_{1, p} )^T, & \quad X_2 = ( X_{2, 1},  X_{2, 2}, \cdots, X_{2, p} )^T, \\
Y_1 = ( Y_{1,1}, Y_{1, 2}, \cdots, Y_{1, q} )^T, & \quad Y_2 = ( Y_{2, 1},  Y_{2, 2}, \cdots, Y_{2, q} )^T.   
\end{align*}
We can develop the following proposition by resorting to Theorem \ref{thm4} for the case of $ 0 < \tau \leq 1 / 2  $  and Theorem \ref{thm8} in Section \ref{SecD.1} of Supplementary Material for the case of $ 1/2 < \tau \leq 1 $.  

\begin{proposition}   \label{cor1}
	Assume that $  \e ( | X_{1, i}| ^{4 + 4 \tau} )  < \infty $ and $  \e (  | Y_{1, j}| ^{ 4 + 4 \tau } )  < \infty $ for any $ 1 \leq  i \leq  p , 1 \leq j \leq q  $ with some constant $ 0 < \tau \leq 1  $, and there exist some positive constants $ \kappa_1,  \kappa_2 , \kappa_3, \kappa_4 $ such that 
	\begin{align}
	& \max \Big\{ p ^{ - 1 } \textstyle \sum \nolimits_{i = 1}^{p }  \e [ | X_{1, i} |^{ 4 + 4 \tau } ]   ,  ~ q ^{ - 1 }  \textstyle \sum \nolimits_{j = 1}^{q }   \e [ | Y_{1, j} | ^{ 4 + 4 \tau } ]  \Big\}   \leq  \kappa_1,  \label{cor1-c1}  	\\  	         
	&  \min \big\{  p ^{ - 1} \e [ ( X_1^T X_2 )^2 ]   ,   ~ q ^{ - 1 }  \e [ ( Y_1^T Y_2 )^2 ]  \big\}  \geq \kappa_2,  \label{cor1-c2} \\
	& \min \big\{ p ^{ - 1 } B_X,  ~ q ^{ - 1 } B_Y \big\}  \geq \kappa_3,     \label{cor1-c3}  \\
	& \max  \limits_{ 1 \leq i \leq p  } \e [ X_{1, i}^2 ] \leq \kappa_4, \quad \max  \limits_{ 1 \leq j \leq q } \e [ Y_{1, j}^2 ] \leq \kappa_4. \label{cor1-c4}  	              
	\end{align}
	In addition, assume that $ X $ consists of $ m_1 $-dependent components, $ Y $ consists of $ m_2 $-dependent components, and 
	\begin{align}  
	m_1 = o (   p ^{ \tau / ( 2 + \tau ) }  ),  \quad  m_2 = o (  q ^{ \tau / ( 2 + \tau ) }  ),  \quad m_1 m_2 = o (  n^{ \tau / ( 1 + \tau ) }  ) .  \label{range-m}
	\end{align}
	Then under the independence of $ X $ and $ Y $, there exists some positive constant $ A $ depending upon $ \kappa_1, \cdots,  \kappa_4 $ such that
	\begin{align}
	& \sup \limits_{ x \in \mathbb{R} }  | \mathbb{P} ( T_n \leq x )  -  \Phi ( x )    |   \nn \\ 
	& \leq   A  \Big[  \Big(  [ ( m_1 + 1  ) ( m_2 + 1  ) ]^{ 2 + \tau }  ( p  q   )^{ - \tau }  \Big)^{  \frac { 1 + \tau  } { 2 }  }  +  [ ( m_1 + 1  ) ( m_2 + 1  ) ]^{ 1 + \tau }   n^{ - \tau }  \Big]^{ \frac { 1 } { 3 + 2 \tau  } }.  \label{rate-dep}
	\end{align}
	Hence under condition \eqref{range-m}, we have 
	$ T_n \stackrel { \mathscr{D}} {\rightarrow} N (0, 1)$ as $ n \rightarrow \infty  $ and $ p ,  q \rightarrow \infty $.
\end{proposition} 

\ignore{ 
	As illustrated in this $m$-dependent case, results get better for larger $ \tau $ in terms of both the range of dependence bandwidth $ m_1, m_2 $ and the convergence rate. It is natural that higher moments deserve better asymptotic performance. Indeed, $ \tau / ( 2 + \tau  ) $ and $ \tau / ( 1 + \tau  ) $ are increasing functions with respect to $ \tau  $, so the range of $ m_1 $ and $ m_2 $ increases with $ \tau $; as for the convergence rate, by straightforward calculation it is easy to see under condition \eqref{range-m} that for any $ 0 <  \tau < \tau ' \leq 1  $,
	\begin{align*}
	\frac  {  [ ( m_1  + 1 ) ( m_2 + 1 ) ]^{   \frac { ( 2 + \tau' ) ( 1 + \tau' ) }  { 2  ( 3 + 2 \tau'  )  }   } ( p_n q_n )^{  -  \frac {   ( \tau' ) ( 1 + \tau' ) } { 2 ( 3 + 2 \tau'  ) } }    } {  [ ( m_1  + 1 ) ( m_2 + 1 ) ]^{  \frac { ( 2 + \tau ) ( 1 + \tau  ) }  { 2 ( 3 + 2 \tau  ) }   } ( p_n q_n )^{  -  \frac { \tau ( 1 + \tau  ) } { 2 ( 3 + 2 \tau  )  } }  }    \rightarrow  0,
	\end{align*}
	and 
	\begin{align*}
	\frac { [ ( m_1 + 1  ) ( m_2 + 1  )  ]^{ \frac { 1 + \tau' } {  3 + 2 \tau'  } }  n^{ - \frac { \tau'  } { 3 + 2 \tau'  } } } {   [ ( m_1 + 1  ) ( m_2 + 1  )  ]^{ \frac { 1 + \tau   } {  3 + 2 \tau    } }  n^{ - \frac { \tau    } { 3 + 2 \tau   } }  }    \rightarrow  0.
	\end{align*}
}

{\cite{ZYZS2019} also established the asymptotic normality of the rescaled distance correlation. For clear comparison, we summarize in Table \ref{table-dif} the key differences between our results and theirs under the assumptions of Proposition \ref{cor1} and the existence of the eighth moments ($\tau = 1$).}

\begin{table}[htbp]
	\renewcommand\arraystretch{1.35}
	\centering
	\caption{{Comparison under the assumptions of Proposition \ref{cor1} } } \label{table-dif}
	\begin{tabular}{c c c c}
		\toprule
		&  \multicolumn{3}{c}{Conditions for asymptotic mormality }  \\  \cline{2-4}	
		& $ p \to \infty , q \to \infty$ & & $ p \to \infty$, $q  $ fixed \\
		&&  & (similarly for $p$ fixed, $ q \to \infty $) \\ 
		\cline{2-2} \cline{4-4}
		\multicolumn{1}{c}{\multirow{4}{*}{\cite{ZYZS2019}}} & \multicolumn{1}{c} {$ m_1^3 / p \rightarrow 0 $, $ m_2 ^3 / q \rightarrow 0 ,$ } & & \multirow{4}{*}{No result}\\
		\multicolumn{1}{c}{}& \multicolumn{1}{c}{ $ m_1 / n^{1/4} \rightarrow 0 $, $  m_2 / n^{1/4} \rightarrow 0 ,$ } & &  \\
		\multicolumn{1}{c}{}& \multicolumn{1}{c}{$ n \sqrt{m_1} m_2 / \sqrt{q}  \rightarrow 0 ,$ }& & \\ 
		\multicolumn{1}{c}{}&  \multicolumn{1}{c}{ $ n m_1 \sqrt{m_2} / \sqrt{p} \rightarrow 0 .$ } & & \\  
		\cline{2-2} \cline{4-4} 
		\multicolumn{1}{c}{\multirow{2}{*}{Our work} } & \multicolumn{1}{c} {$  m_1^3 / p \rightarrow 0 $, $ m_2 ^3 / q \rightarrow 0 ,$ } & & $ m_1^3 /p \rightarrow 0 ,$ \\
		\multicolumn{1}{c}{} & \multicolumn{1}{c}{ $ m_1   m_2  / \sqrt{n} \rightarrow 0  .$ } & & $m_1 / \sqrt{n}  \rightarrow 0.$ \\
		\bottomrule
	\end{tabular}
\end{table}

We further consider the third example of multivariate normal random variables. For such a case, we can obtain a concise result in the following proposition. 

\begin{proposition} \label{cor2}
	Assume that $ X \sim N(0, \Sigma_x) $, $ Y \sim N(0, \Sigma_y) $, and the eigenvalues of $ \Sigma_x $ and $ \Sigma_y $ satisfy that $ a_1 \leq  \lambda_1^X \leq \lambda_2^X \leq \cdots \leq \lambda_{ p }^X \leq a_2  $ and $ a_1 \leq  \lambda_1^Y \leq \lambda_2^Y \leq \cdots \leq \lambda_{ q }^Y \leq a_2 $ for some positive constants $ a_1 $ and $ a_2 $. Then under the independence of $ X $ and $Y$, there exists some positive constant $ C   $ depending upon $ a_1, a_2 $ such that 
	\begin{align*}
	\sup \limits_{ x \in \mathbb{R} } | \mathbb{P } ( T_n \leq x  )  - \Phi (x) | \leq C  \big[ ( p  q )^{ -  1/ 5  }  + n^{ - 1/ 5 }  \big].
	\end{align*}
	Hence we have $ T_n \stackrel{ \mathscr{D} }{ \rightarrow } N(0, 1) $ as $ n \rightarrow \infty $ and $ p ,  q \rightarrow \infty $. 
\end{proposition}

We would like to point out that the rate of convergence obtained in Proposition \ref{cor2} can be suboptimal since the error rate $ n^{ - 1/ 5 } $ is 
slower than the classical convergence rate with order $ n^{ - 1/2 } $ of the CLT for the sum of independent random variables. Our results are derived by exploiting the convergence rate of CLT for general martingales \citep{H1988}. It may be possible to improve the rate of convergence if one takes into account the specific intrinsic structure of distance covariance, which is beyond the scope of the current paper.

\subsection{Power analysis} \label{power}
We now turn to 
the power analysis for the rescaled distance correlation. We start with presenting a general theory on power in Theorem \ref{thm-power} below. Let us define two quantities
\begin{align}
L_x  = \e \big(\big | \Vert X \Vert^2 - \e \Vert X \Vert^2 \big|^{  4 }  \big)  + \e \big( | X_1^T X_2 |^4 \big)  , ~
L_y  = \e \big(\big | \Vert Y \Vert^2 - \e \Vert Y \Vert^2 \big|^{  4 }  \big)  + \e \big( | Y_1^T Y_2 |^4 \big) . \label{neweq.006}
\end{align}

\begin{theorem}   \label{thm-power}
	Assume that $ \e  ( \Vert X \Vert ^{ 8 }   ) + \e (  \Vert Y \Vert ^{ 8 }   ) < \infty   $ and \eqref{cond1} holds with $\tau = 1$. If $   n  \Rc^2 (X, Y)   \to \infty $ and $ \sqrt{n} \V^2 (X, Y) / \big(  B_X^{-1/2} B_Y^{-1/2} L_{x}^{1/4} L_{y}^{1/4} \big)  \to \infty$, then {
		for any arbitrarily large constant $C>0$, $\mathbbm{P}(T_n>C)\rightarrow 1$ as $n\rightarrow \infty$.}
	Thus, for any significance level $\alpha$, $ \mathbbm {P} (T_n > \Phi^{-1} (1 - \alpha) ) \to 1 $ as $n \to \infty$, where $ \Phi^{-1} (1 - \alpha) $ represents the $(1 - \alpha)$th quantile of the standard normal distribution.
\end{theorem}

{
	Theorem \ref{thm-power} provides a general result on the power of the rescaled distance correlation statistic. It reveals that as long as the signal strength, measured by $\Rc^2 (X, Y)$ and $\V^2 (X, Y)$, is not too weak, the power of testing independence with the rescaled sample distance correlation can be asymptotically one. In most cases, the population distance variances $ \V^2 (X) $ and $\V^2 (Y)$ are of constant order by Proposition \ref{prop2}. Therefore, if $ B_X^{ - 1/2 } B_Y^{-1/2} L_x^{1/4} L_{y}^{1/4} $ is also of constant order, then the conditions in Theorem \ref{thm-power} will reduce to $ \sqrt n \Rc^2 (X, Y) \to \infty $, which indicates that the signal strength should not decay faster than $ n^{ - 1/2} $.
	To gain some insights, assume that both $X\in \R^p$ and $Y\in \R^q$ consist of independent components with uniformly upper bounded eighth moments and uniformly lower bounded second moments. Then it holds that $ B_X = O (p) $, $B_Y = O (q) $, $ L_{x} = O (p^2) $, $ L_{y} = O(q^2) $, $ \V^2 (X) = O (1) $, and $ \V^2 (Y) = O (1) $. Thus the conditions in Theorem \ref{thm-power} above reduce to $ \e  ( \Vert X \Vert ^{ 8 }   ) + \e (  \Vert Y \Vert ^{ 8 }   ) < \infty   $ and $\sqrt n {\Rc^2} (X, Y) \to \infty $.
	In general, $ \Rc^2 (X, Y) $ and $\V^2 (X, Y) $ depend on the dimensionalities and hence the conditions of Theorem \ref{thm-power} impose certain relationship between $n$ and $p$. 
}



Recently \cite{ZYZS2019} showed that in the asymptotic sense, the distance covariance detects only componentwise linear dependence in the high-dimensional setting when both dimensionalities $p$ and $q$ grow much faster than sample size $n$ (see Theorems 2.1.1 and 3.1.1 therein). In particular, when $X$ and $Y$ both consist of i.i.d. components with certain bounded moments, distance covariance was shown to asymptotically measure linear dependence if  
$ \min \{ p, q    \} / n^2 \to \infty $. 
However, in view of \eqref{DisCov} and \eqref{neweq.001}, the population distance covariance and distance correlation indeed characterize completely the independence between two random vectors in arbitrary dimensions. Therefore, it is natural to ask whether  the sample distance correlation can detect nonlinear dependence in some other diverging regime of $(n, p, q)$. 
The answer turns out to be affirmative in the regime of moderately high dimensionality: {We formally present this result in the following theorem on the asymptotic power and compare with the results in \cite{ZYZS2019} in Table \ref{dif-power}.}


\begin{theorem}  \label{prop-power}
	Assume that we have i.i.d. observations $ \{ (X_i, Y_i) , 1 \leq i \leq n\} $ with $ X_i \in \R^p $ and $ Y_i \in \R^p $, $ X_1 = (X_{1, 1}, \ldots, X_{1, p})  $ with $X $ having a symmetric distribution, and $ \{X_{1, i},  1\leq i \leq p \} $ are m-dependent for some fixed positive integer $m$.  Let $ Y_1  = (Y_{1, 1}, \ldots, Y_{1, p}) $ be given by $ Y_{1, j} = g_j (X_{1, j})  $ for each $1 \leq j \leq p$, where $ \{g_j, 1 \leq j \leq p  \}$ are symmetric functions satisfying $ g_j(x) = g_j( - x) $ for $x \in \R$ and $ 1 \leq j \leq p  $. Assume further that $ \e ( X_{1, j}^{ 12 }) \  +  \e (Y_{1, j}^{12}) \leq c_1^{12}$, $  \var(X_{1, j}) \geq c_2^2 $, and $ \var (Y_{1, j}) \geq c_2^2 $ for some positive constants $c_1, c_2$. Then there exists some positive constant $A$ depending on $c_1, c_2$, and $m$ such that 
	\begin{align*}
	\V^2 (X, Y) & \geq A p^{- 1} + O (p^{ - 3/2} ) \\
	\mbox{and}\quad 
	\Rc^2 (X, Y) & \geq A p^{-1} + O(p^{ - 3/2}).
	\end{align*}   
	Consequently, if $ p = o (\sqrt n) $, then for any arbitrary large constant $C> 0$, $\mathbbm{P} (T_n > C) \to 1$ as $n \to \infty$, and thus the test of independence between $X$ and $Y$ based on the rescaled sample distance correlation $T_n$ has asymptotic power one. 
\end{theorem} 

Under the symmetry assumptions in Theorem \ref{prop-power}, we can show that there is no linear dependence between $X$ and $Y$ by noting that $ \cov (X_{1, i}, Y_{1, j})  = 0 $ for each $1 \leq i, j \leq p$. It is worth mentioning that we have assumed the $m$-dependence for some fixed integer $ m \geq 1 $ to simplify the technical analysis. In fact, $m$ can be allowed to grow slowly with sample size $n$ and our technical arguments are still applicable. 

\begin{table}[htbp]
	\centering
	\caption{Comparison of power analysis in detecting pure nonlinear dependency} \label{dif-power}
	\begin{tabular}{l c}
		\toprule
		\multirow{3}{*}{\cite{ZYZS2019}} & \multicolumn{1}{l}{Asymptotically no power when $p $ and $q$ grow much faster than $n$ } \\
		& \multicolumn{1}{l}{(especially it requires $\min\{p, q\} \gg n^2$ when $X, Y$ consist of}  \\
		&  \multicolumn{1}{l}{~ i.i.d. components)} \\ \midrule
		\multirow{2}{*}{Our work} & \multicolumn{1}{l}{Asymptotically can achieve power one when $p =q = o (\sqrt{n}) $} \\
		& \multicolumn{1}{l}{(under the conditions of Theorem \ref{prop-power}) } \\
		\bottomrule
	\end{tabular}
\end{table}

\section{Simulation studies} \label{simulation}

In this section, we conduct several simulation studies to verify our theoretical results on sample distance correlation and illustrate the finite-sample performance of our rescaled test statistic for the test of independence.

{
	\subsection{Normal approximation accuracy} \label{hist}
	We generate two independent multivariate normal random vectors $ X \in \R^p $ and $ Y \in \R^{p} $ in the following simulated example and calculate the rescaled distance correlation $ T_n $ defined in \eqref{Tn}. 
	
	\begin{example} \label{ex0}
		Let $\Sigma = ( \sigma_{i, j} )  \in \R^{p \times p }$ with $ \sigma_{i, j} = 0.7^{| i - j |} $, and $ X \sim N(0, \Sigma) $ and $ Y \sim N(0,  \Sigma ) $ be independent. {We consider the settings of $ n = 100$ and $p= 10, 50, 200, 500 $.}
	\end{example}
	
	{We conduct $ 5000 $ Monte Carlo simulations and generate the histograms of the rescaled test statistic $  T_n  $  to investigate its empirical distribution. Histograms with a comparison of the kernel density estimate (KDE) and the standard normal density function are shown in Figure \ref{fig:n100}.}
	From the histograms, we can see that the distribution of $    T_n  $ mimics very closely the standard normal distribution under different settings of dimensionalities. Moreover, for more refined comparison, the maximum pointwise distances between the KDE and the standard normal density function under different settings are presented in Table \ref{table-error}. It is evident that the accuracy of the normal approximation increases with dimensionality, which is in line with our theoretical results.

	\begin{table} 
		\footnotesize
		\centering
		\caption{Distances between the KDE and standard normal density function in Example \ref{ex0}.}  
		\begin{tabular}{c c  c | c c c}
			\toprule
			$n   $      &  $ p $      &  Distance   & $n$ & $p$ & Distance  \\ \midrule
			100 & 10 & 0.0955    & 100 & 200 & 0.0288 \\
			100 & 50 &  0.0357  & 100 & 500 & 0.0181 \\
			\bottomrule
		\end{tabular}
		\label{table-error}
	\end{table}
	
\begin{figure}[t!]
	\centering
	\includegraphics[scale=0.45]{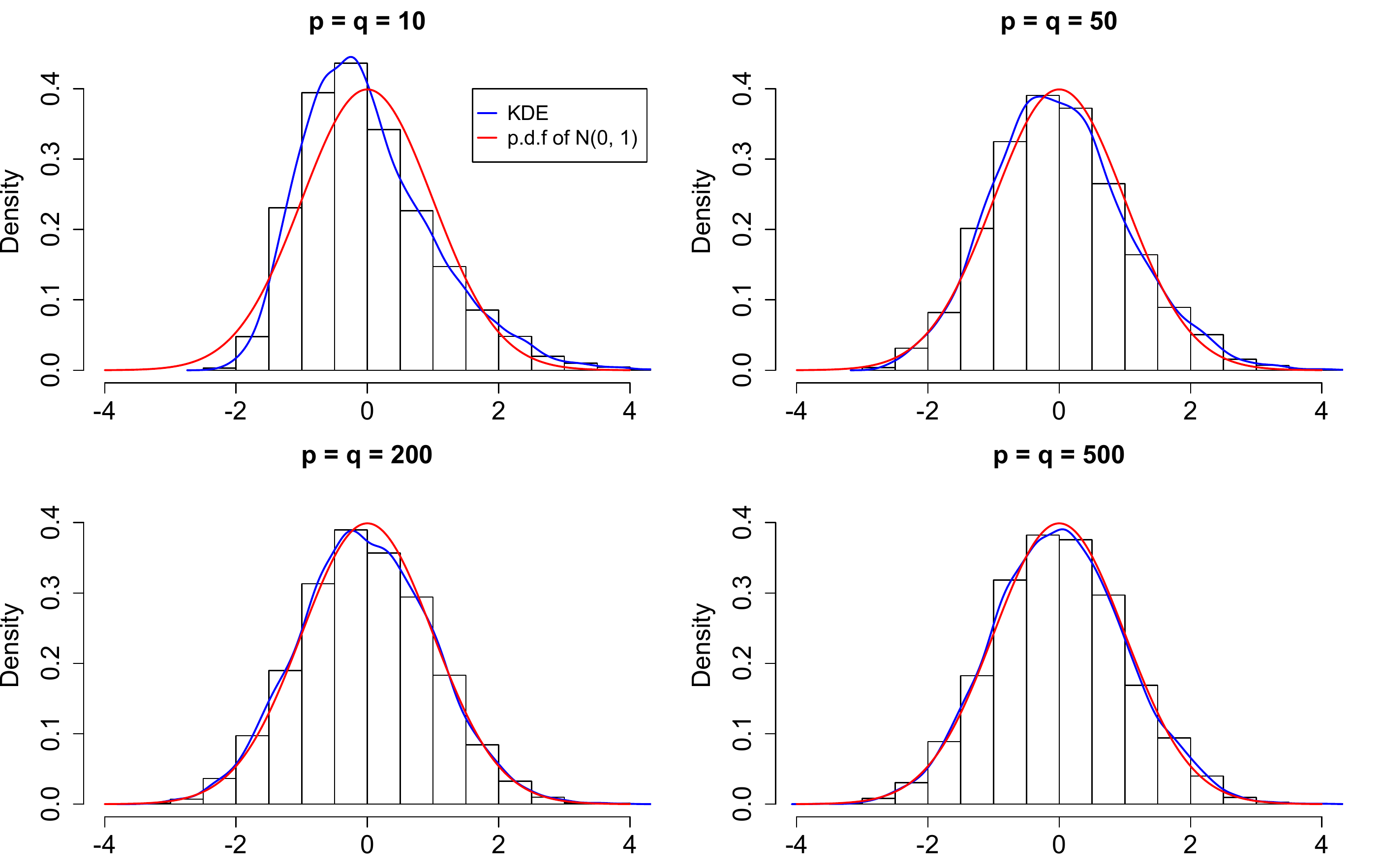}
	\caption{\it Histograms of the rescaled test statistic $  T_n $ in Example \ref{ex0}. 
		The blue curve represents the kernel density estimate and the red curve represents the standard normal density.}
	\label{fig:n100}
\end{figure}



\subsection{Test of independence} \label{test}  
To test the independence of random vectors $ X $ and $ Y $ in high dimensions, based on the asymptotic normality developed for the rescaled distance correlation statistic $T_n$, under significance level $\alpha$ we can reject the null hypothesis when
\begin{equation} \label{neweq.005}
T_n  =  \sqrt{ \frac { n  ( n - 1 ) } { 2 } } \Rc_n^* (X, Y) > \Phi^{ - 1   } ( 1 - \alpha  ),
\end{equation}
since the distance correlation is positive under the alternative hypothesis. {To assess the performance of our normal approximation test, we also include the gamma-based approximation test \citep{HH2017} 
	and normal approximation for studentized sample distance correlation $ T_R $ defined in \eqref{fix-p} \citep{ZYZS2019} in the numerical comparisons.} 

} 
The gamma-based approximation test assumes that the linear combination $ \sum_{ i = 1}^{ \infty  } \lambda_i  Z_i^2  $ involved in the limiting distribution of the standardized sample distance covariance $ n \V_n^* ( X, Y ) $ under fixed dimensionality (see \eqref{finite-p}) can be approximated \textit{heuristically} by a gamma distribution $ \Gamma (\beta_1, \beta_2) $ with matched first two moments. In particular, the shape and rate parameters are determined as 
\begin{align*}
\beta_1 =   \frac { \big( \sum_{ i = 1 }^{ \infty } \lambda_i   \big)^2 } { 2 \sum_{ i = 1 }^ { \infty } \lambda_i^2  } = \frac { \big( \e   \Vert X - X' \Vert   \e   \Vert Y - Y' \Vert \big)^2  } { 2  \V^2 (X) \V^2 ( Y ) }
\end{align*}
and 
\begin{align*}
\beta_2 = \frac {  \sum_{ i = 1 }^{ \infty } \lambda_i  } { 2 \sum_{ i = 1 }^ { \infty } \lambda_i^2  } = \frac { \e   \Vert X - X' \Vert   \e   \Vert Y - Y' \Vert } { 2  \V^2 (X) \V^2 ( Y ) }.
\end{align*}
Thus given observations $ ( X_1, Y_1 ), \cdots, ( X_n, Y_n ) $, $ \beta_1  $ and $ \beta_2 $ can be estimated by their empirical versions
\begin{align*}
\hat{\beta_1}  = \frac {\mu^2  } { 2  \V_n^* ( X ) \V_n^* ( Y ) } \quad \mbox{and} \quad
\hat{\beta_2} = \frac { \mu } { 2  \V_n^* ( X ) \V_n^* ( Y ) },
\end{align*}
where $ \mu =   \frac {1} { n^2 ( n - 1)^2 } \sum_{ i \neq j } \Vert X_i - X_j \Vert \sum_{ i \neq j } \Vert Y_i - Y_j \Vert  $.
Then the null hypothesis is rejected at the significanve level $\alpha$ if $ n \V^*_n (X,Y) > \Gamma_{1 - \alpha} (\hat{\beta}_1, \hat{\beta}_2) - \mu $, where $ \Gamma_{1 - \alpha} (\hat{\beta}_1, \hat{\beta}_2)  $ is the $ (1 - \alpha) $th quantile of the distribution $ \Gamma  (\hat{\beta}_1, \hat{\beta}_2)  $. 
The gamma-based approximation test still lacks rigorous theoretical justification. 

{

\ignore{When the sample size is fixed and the dimensionalities tend to infinity, \cite{SR2013} proved the student's t approximation for the studentized sample distance correlation $T_R $ in \eqref{fix-p} under the assumption that the components of $X$ and $Y$ are i.i.d. with finite variance, which was relaxed later by \cite{ZYZS2019}. Thus under significance level $\alpha$, the null hypothesis is rejected when $ T_R > t_{n (n - 3)/ 2 - 1} (1 - \alpha) $, where $   t_{n (n - 3)/ 2 - 1} (1 - \alpha)  $ represents the $(1 - \alpha)$th quantile of student's t distribution with $ n (n - 3)/ 2 - 1 $ degrees of freedom. Yet the student's t approximation for $T_R$ is proved only when sample size is fixed.}
When the sample size and dimensionalities tend to infinity simultaneously, in view of our main result in Theorem \ref{thm2} and the consistency of $ \Rc_n^* (X, Y) $ (recall Lemma \ref{le-consist} and \eqref{h-b1} in Section \ref{SecC.1} of Supplementary Material), one can see that under the null hypothesis, 
$
T_R \xrightarrow{\mathscr{D}} N(0, 1).
$
Therefore, we can reject the null hypothesis at significance level $\alpha$ if $ T_R > \Phi^{-1} (1 - \alpha) $.

We consider two simulated examples to compare the aforementioned three approaches for testing the independence between two random vectors in high dimensions. The significance level is set as $ \alpha =  0.05 $ and $ 2000 $ Monte Carlo replicates are carried out to compute the empirical rejection rates.

\begin{example}  \label{ex2}
	Let $\Sigma = (\sigma_{i, j})  \in \R^{p \times p } $ with $ \sigma_{i, j} = 0.5^{| i - j |} $. Let $ X $ and $Y$ be independent and  $ X \sim N (0, \Sigma ) $, $ Y \sim N(0, \Sigma) $.
\end{example}

\begin{example} \label{ex3}
	Let $ \Sigma = (\sigma_{i, j})  \in \R^{p \times p} $ with $ \sigma_{i, j} = 0.5^{| i - j |} $. Let $ X = (X^{(1)}, \ldots, X^{(p)}) \sim N (0, \Sigma ) $ and $ Y = (Y^{(1)}, \ldots, Y^{(p)}) $ with $ Y^{(i)} = 0.2 \big( X^{(i)} + (X^{(i)} )^2 \big) + \varepsilon_i  $ and $ \varepsilon_i \stackrel{i.i.d.}{\sim}  t_4 $. 
\end{example}

Type-I error rates in Example \ref{ex2} under different settings of $n$ and $p$ are presented in Figure \ref{fig:ex2}. From Figure \ref{fig:ex2}, it is easy to see that the rejection rates of the normal approximation test for $T_n$ tend to be closer and closer to the preselected significance level as the dimensionalities and the sample size grow. The same trend applies to the other two approches too. The empirical powers of the three tests in Example \ref{ex3} are shown in Figure \ref{fig:ex3}. We can observe from the simulation results in Figures \ref{fig:ex2} and \ref{fig:ex3} that these three tests perform asymptotically almost the same, which is sensible. Empirically, the gamma approximation for $ n \V^*_n (X, Y) $ and normal approximation may be asymptotically equivalent to some extent and more details on their connections are discussed in Section \ref{normal-gamma} of Supplementary Material. However, the theoretical foundation of the gamma approximation for $ n \V_n^* (X, Y)$ remains undeveloped.
As for the asymptotic equivalence between $T_n$ and the studentized sample distance correlation $T_R $, Lemma \ref{le-consist} and \eqref{h-b1} imply that under the null hypothesis and some general conditions, $ \Rc^*_n (X, Y) \to 0 $ in probability and hence $ T_R $ can be asymptotically equivalent to $T_n $ when $ n \to \infty $.

\begin{figure}[t!]
	\centering
	\includegraphics[scale=0.45]{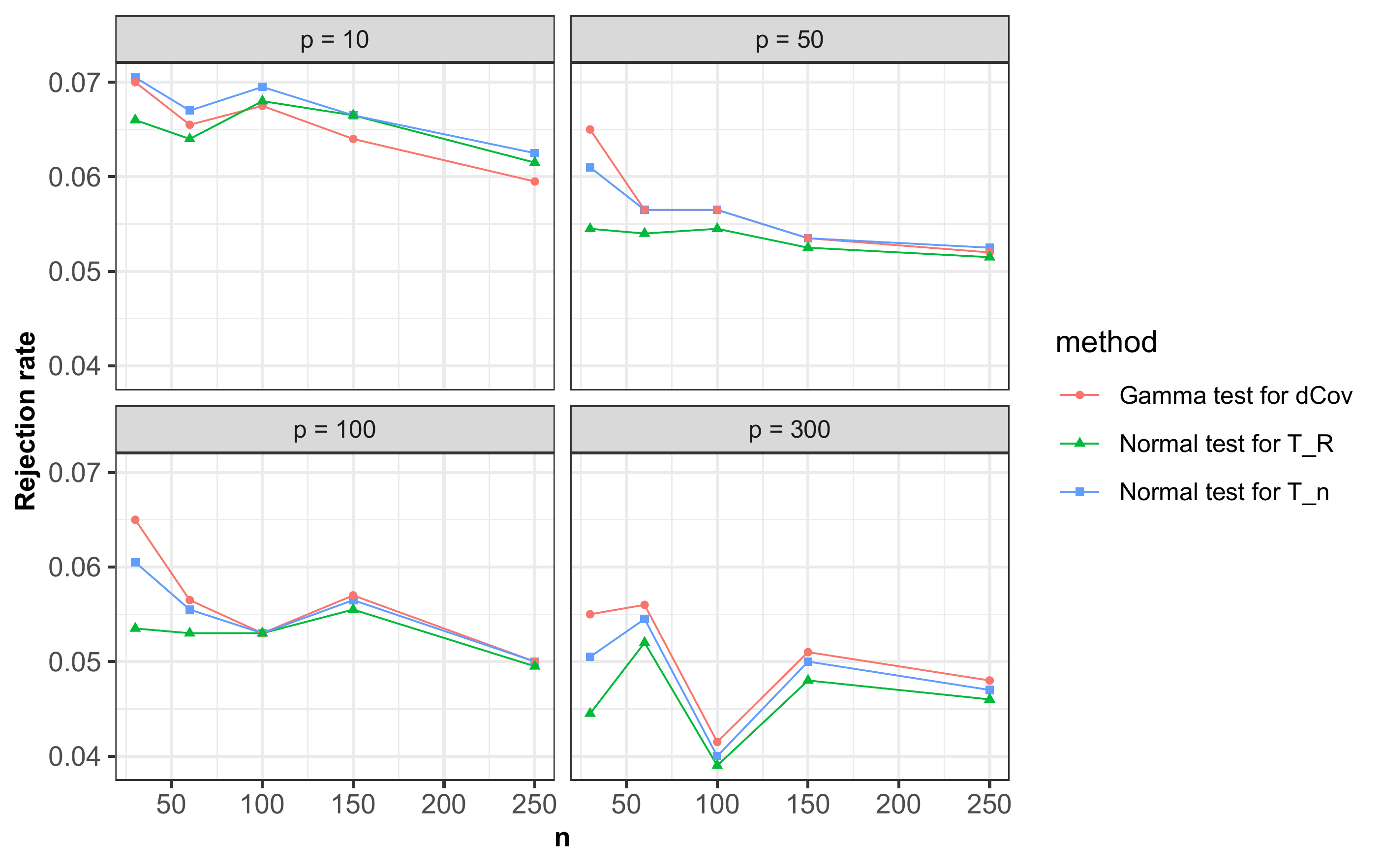}
	\caption{\it  Rejection rates of the three approaches under different settings of $n$ and $p$ in Example \ref{ex2}.}
	\label{fig:ex2}
\end{figure}

\begin{figure}[t!]
	\centering
	\includegraphics[scale=0.45]{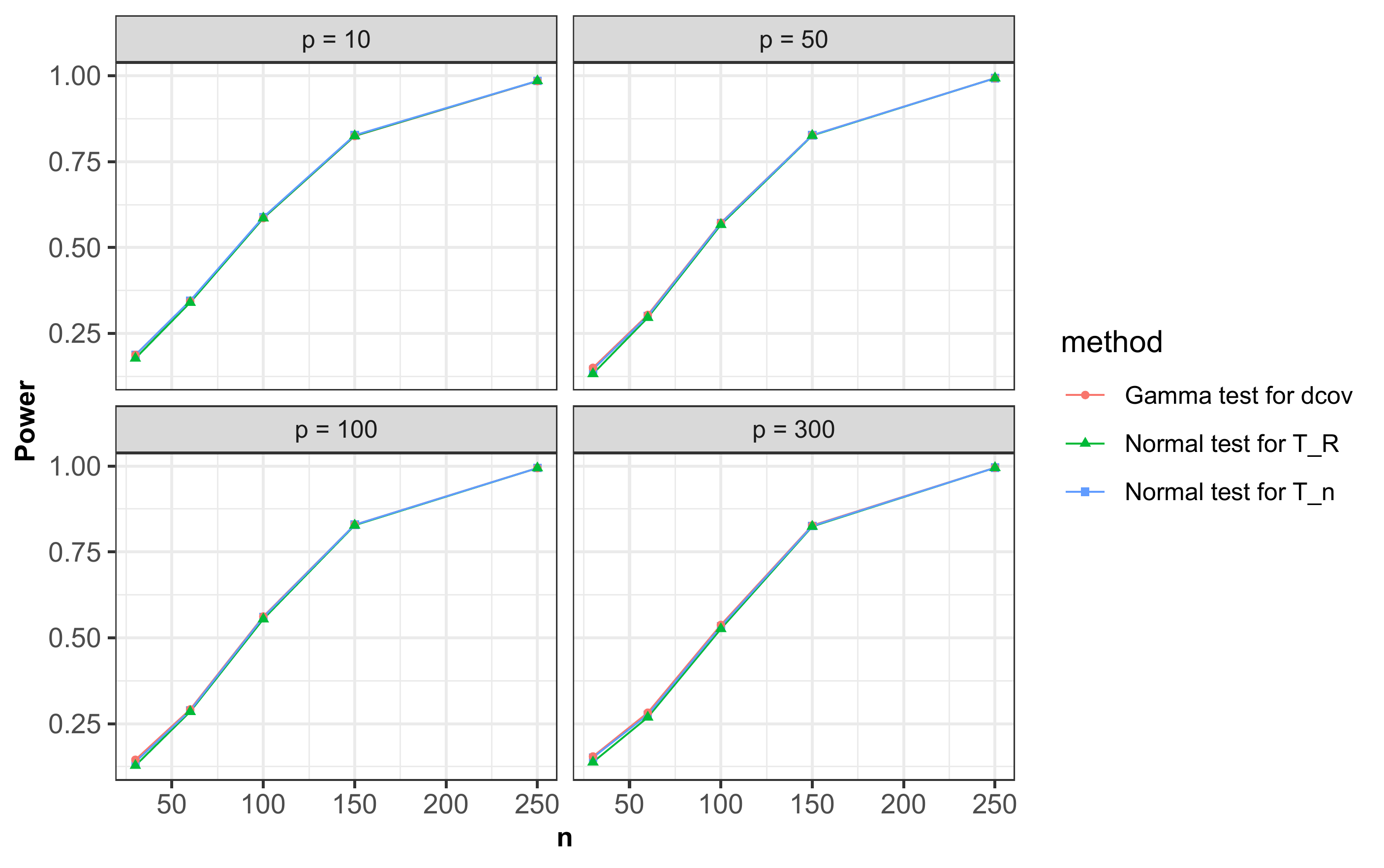}
	\caption{\it Power of the three approaches under different settings of $n$ and $p$ in Example \ref{ex3}.}
	\label{fig:ex3}
\end{figure}

\begin{table}
	\footnotesize
	\centering
	\caption{Power of our rescaled test statistic with $p = 2 [\sqrt{n}]$ in Examples \ref{ex-power1} and \ref{ex-power2} (with standard errors in parentheses). }  
	\begin{tabular}{c c  c | c c c}
		\toprule
		\multicolumn{3}{c}{Example \ref{ex-power1}} & \multicolumn{3}{c}{Example \ref{ex-power2}}   \\ \midrule
		$n   $      &  $ p $      &  Power & $n $ & $p$ & Power   \\ \midrule 
		10 & 6  &  0.2765 (0.0100)   &	10 & 6  &   0.3060 (0.0103)  \\
		40 & 12 &  0.5165 (0.0112)  & 40 & 12 &   0.7005 (0.0102)  \\
		70 & 16 &   0.6970 (0.0103) & 70 & 16 &   0.9380 (0.0054)   \\
		100 & 20   & 0.8220 (0.0086)  & 100 & 20   & 0.9885 (0.0024) \\
		130 & 22 & 0.9270 (0.0058)  & 130 & 22   & 0.9995 (0.0005)   \\
		160 & 26 & 0.9550 (0.0046) & 160 & 26 & 0.9990 (0.0007) \\
		\bottomrule
	\end{tabular}
	\label{table-p1}
\end{table}

\subsection{Detecting nonlinear dependence} \label{Sec4.3}
We further provide several examples to justify the power of the rescaled distance correlation statistic in detecting nonlinear dependence in the regime of moderately high dimensionality. In the following simulation examples, the significance level of test is set as 0.05 and 2000 Monte Carlo replicates are conducted to compute the rejection rates. 

\begin{example} \label{ex-power1}
	Let $ X = (X^{(1)}, \ldots, X^{(p)})^T \sim N (0, I_p) $ and $ Y = (Y^{(1)}, \ldots, Y^{(p)})^T  $ satisfying $ Y^{(i)} = (X^{(i)})^2 $. 
\end{example}

\begin{example} \label{ex-power2}
	Set $\Sigma = (\sigma_{i, j}) \in \R^{p \times p} $ with $ \sigma_{i, j}  = 0.5^{ | i - j | }$. Let $ X = (X^{(1)}, \ldots, X^{(p)}) \sim N(0,  \Sigma ) $ and $ Y = (Y^{(1)}, \ldots, Y^{(p)})^T $ with $ Y^{(i)} = (X^{(i)} )^2 $. 
\end{example}

For the above two examples, it holds that $ \cov(X^{(i)}, Y^{(j)}) =0  $ for each $   1 \leq i, j \leq p $. Simulation results on the power under Examples \ref{ex-power1} and \ref{ex-power2} for different settings of $n$ and $p$ are summarized in Table \ref{table-p1}. 
{Guided by Theorem \ref{prop-power}, we set $ p = 2 [\sqrt{n}] $ with $[\cdot]$ denoting the integer part of a given number.}
From Table \ref{table-p1}, we can see that even though there is only nonlinear dependency between $X$ and $Y$, the power of rescaled distance correlation can still approach one when the dimensionality $p$ is moderately high. One interesting phenomenon is that the power in Example \ref{ex-power2} is higher than that in Example \ref{ex-power1}, which suggests that the dependence between components may strengthen the dependency between $ X $ and $Y$.

Moreover, we investigate the setting when one dimensionality is fixed and the other one tends to infinity.

\begin{example} \label{ex-power3}
	Set $\Sigma = (\sigma_{i, j} ) \in \R^{p \times p}$ with $\sigma_{i, j} = 0.7^{ | i - j | } $. Let $ X = (X^{(1)}, \ldots, X^{(p)})  \sim N (0, \Sigma)$ and $ Y = ( \sum_{i = 1}^p X^{(i)} )^2/ p $. 
\end{example}

\begin{figure}[t!]
	\centering
	\includegraphics[scale=0.45]{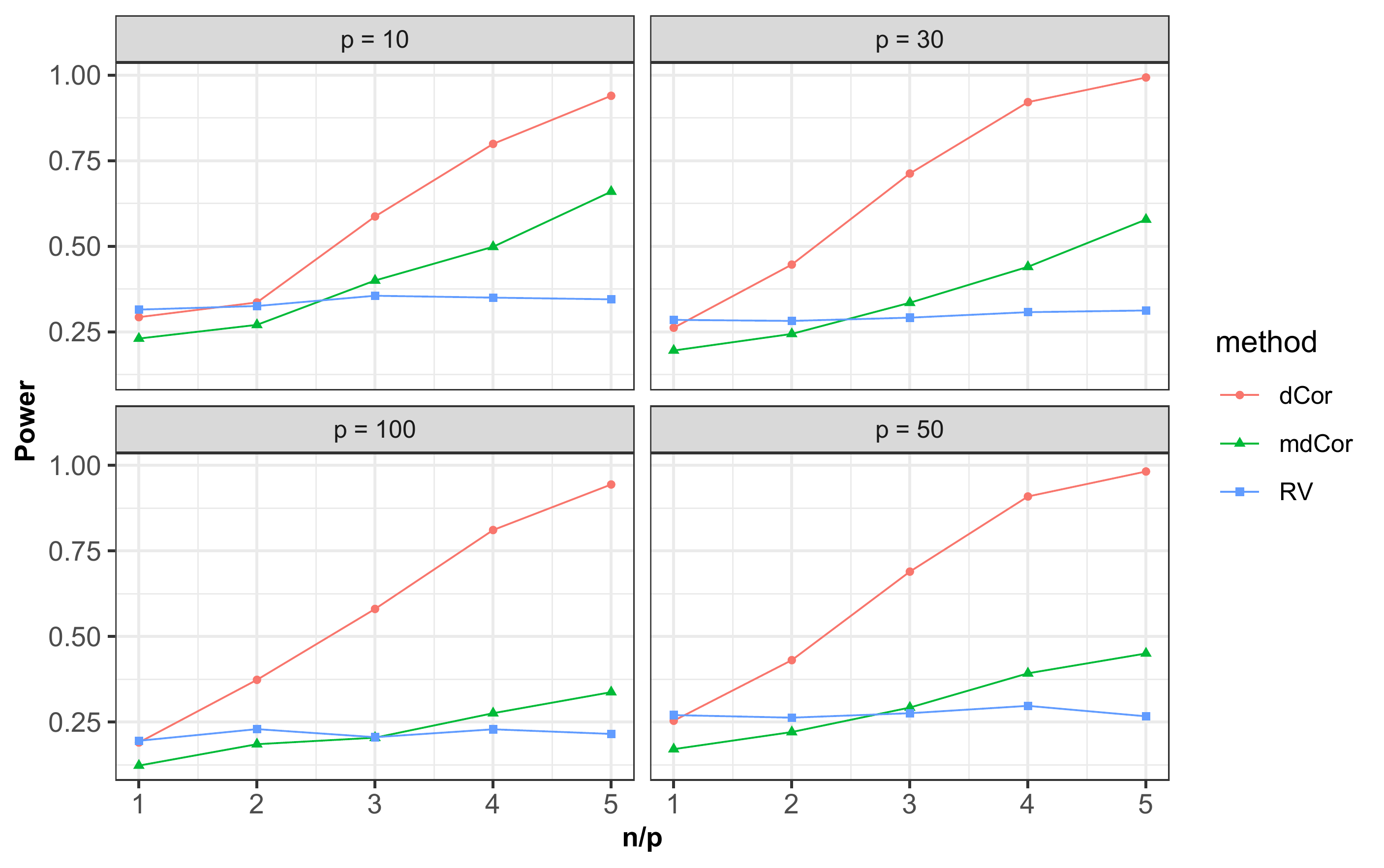}
	\caption{\it Comparison of power under different settings of $n$ and $p$ in Example \ref{ex-power3}. }
	\label{fig:ex6}
\end{figure}

For Example \ref{ex-power3}, it holds that $ \cov(X^{(i)}, Y) = 0  $ for each $ 1 \leq i \leq p  $ and thus the dependency is purely nonlinear. We compare the power of our rescaled distance correlation statistic with the marginally aggregated distance correlation (mdCor) statistic \citep{ZYZS2019} and the linear measure of RV coefficient \citep{E1973,RE1976}. The comparison under different settings of $p$ and $n$ are presented in Figure \ref{fig:ex6}. We can observe from Figure \ref{fig:ex6} that under this scenario, the rescaled distance correlation statistic significantly outperforms the marginally aggregated distance correlation statistic. This is because the marginally aggregated statistic can detect only the marginal dependency between $X$ and $Y$, while $Y$ depends on the entire $X$ jointly 
in this example. Since the RV coefficient measures the linear dependence, its power stays flat and low when the sample size increases.  

These simulation examples demonstrate the capability of distance correlation in detecting nonlinear dependence in the regime of moderately high dimensionality, which is in line with our theoretical results on the power analysis in Theorem \ref{prop-power}. Moreover, when $X$ and $Y$ depend on each other far from marginally, the marginally aggregated distance correlation statistic can indeed be less powerful than the rescaled distance correlation statistic. 
}

\section{Real data application} \label{Sec5}

We further demonstrate the practical utility of our normal approximation test for bias-corrected distance correlation on a blockchain application, which has gained increasing public attention in recent years. Specifically, we would like to understand the nonlinear dependency between the cryptocurrency market and the stock market through the test of independence. Indeed investors are interested in testing whether there is any nonlinear association between these two markets since they want to diversify their portfolios and reduce the risks. In particular, we collected the historical daily returns over recent three years from 08/01/2016 to 07/31/2019 for both stocks in the Standard \& Poors 500 (S\&P 500) list (from https://finance.yahoo.com) and the top 100 cryptocurrencies (from https://coinmarketcap.com). As a result, we obtained a data matrix of dimensions $ 755 \times 505 $ for stock daily returns and a data matrix of dimensions $ 1095 \times 100 $ for cryptocurrency daily returns, where the rows correspond to the trading dates and the columns represent the stocks or cryptocurrencies. Since stocks are traded only on Mondays through Fridays excluding holidays, we adapted the cryptocurrency data to this restriction and picked a submatrix of cryptocurrency data matrix to match the dates. Moreover, because some stocks and cryptocurrencies were launched after 08/01/2016, there are some missing values in the corresponding columns. We removed those columns containing missing values. Finally, we obtained a data matrix $X_{T \times N_1} $ for stock daily returns and a data matrix $Y_{T \times N_2}$ for cryptocurrency daily returns, where $ T = 755 $, $ N_1 = 496 $, and $ N_2 = 22 $. Although the number of cryptocurrencies drops to 22 after removing the missing values, the remaining ones are still very representative in terms of market capitalization, which include the major cryptocurrencies such as Bitcoin, Ethereum, Litecoin, Ripple, Monero, and Dash.

\begin{figure}
	\centering
	\includegraphics[scale=0.48]{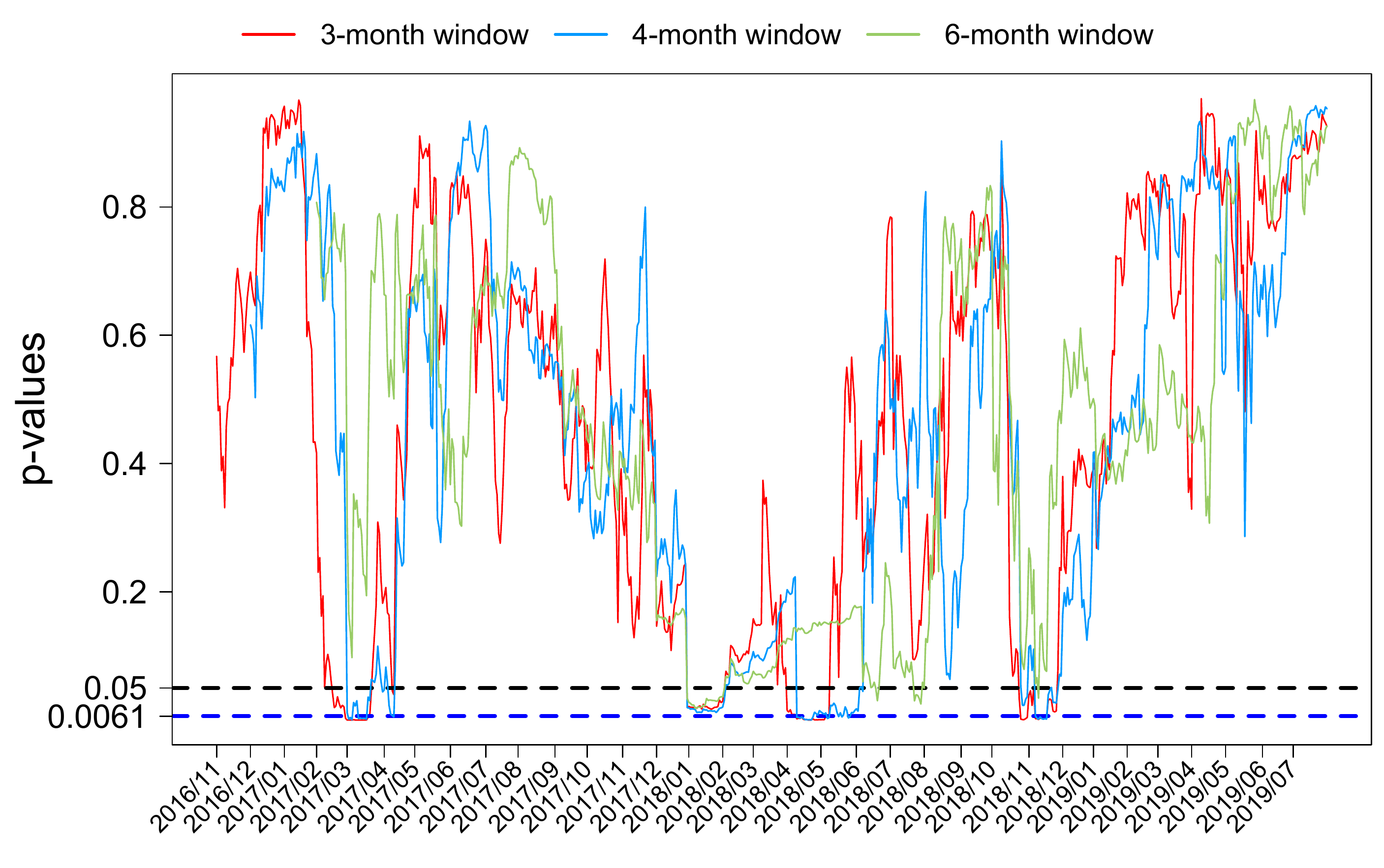}
	\caption{  \it Time series plots of p-values from 11/01/2016 to 07/31/2019 using three-month, four-month, and six-month rolling windows, respectively. 
	}
	\label{fig:real_data}
\end{figure}

To test the independence of the cryptocurrency market and the stock market, we choose three-month rolling windows (66 days). 
Specifically, for each trading date $t$ from 11/01/2016 to 07/31/2019, we set $ X_{F_t \times N_1} $ as a submatrix of $ X_{T \times N_1} $ that contains the most recent three months before date $t$, where $ F_t  $ is the set of 66 rows right before date $t$ (including date $t$). The data submatrix $ Y_{F_t \times N_2 } $ is defined similarly. Then we apply the rescaled test statistic $ T_n  $ defined in \eqref{Tn} to $ X_{F_t \times N_1} $ and $ Y_{F_t \times N_2} $. Thus the sample size $n= 66$ and the dimensions of the two random vectors are $ N_1 = 496 $ and $ N_2 = 22 $, respectively. For each trading date, we obtain a p-value calculated by $ 1 - \Phi ( T_n^{(t)} ) $, where $ T_n^{ (t) } $ is the value of the test statistic based on $ X_{F_t \times N_1} $ and $ Y_{ F_t \times N_2 } $ and $ \Phi (\cdot) $ is the standard normal distribution function.  As a result, we end up with a p-value vector consisting of $ T_n^{ (t) } $ for trading dates $ t $ from 11/01/2016 to 07/31/2019. In addition, we use the ``fdr.control" function in R package ``fdrtool," which applies the algorithms in \cite{BH1995} and \cite{S2002} to calculate the p-value cut-off for controlling the false discovery rate (FDR) at the $ 10 \% $ level. Based on the p-value vector, we obtain the p-value cut-off of 0.0061. {The time series plot of the p-values is shown in Figure \ref{fig:real_data} (the red curve).}

\begin{figure}
	\centering
	\parbox{10cm}{
		\includegraphics[width=10cm]{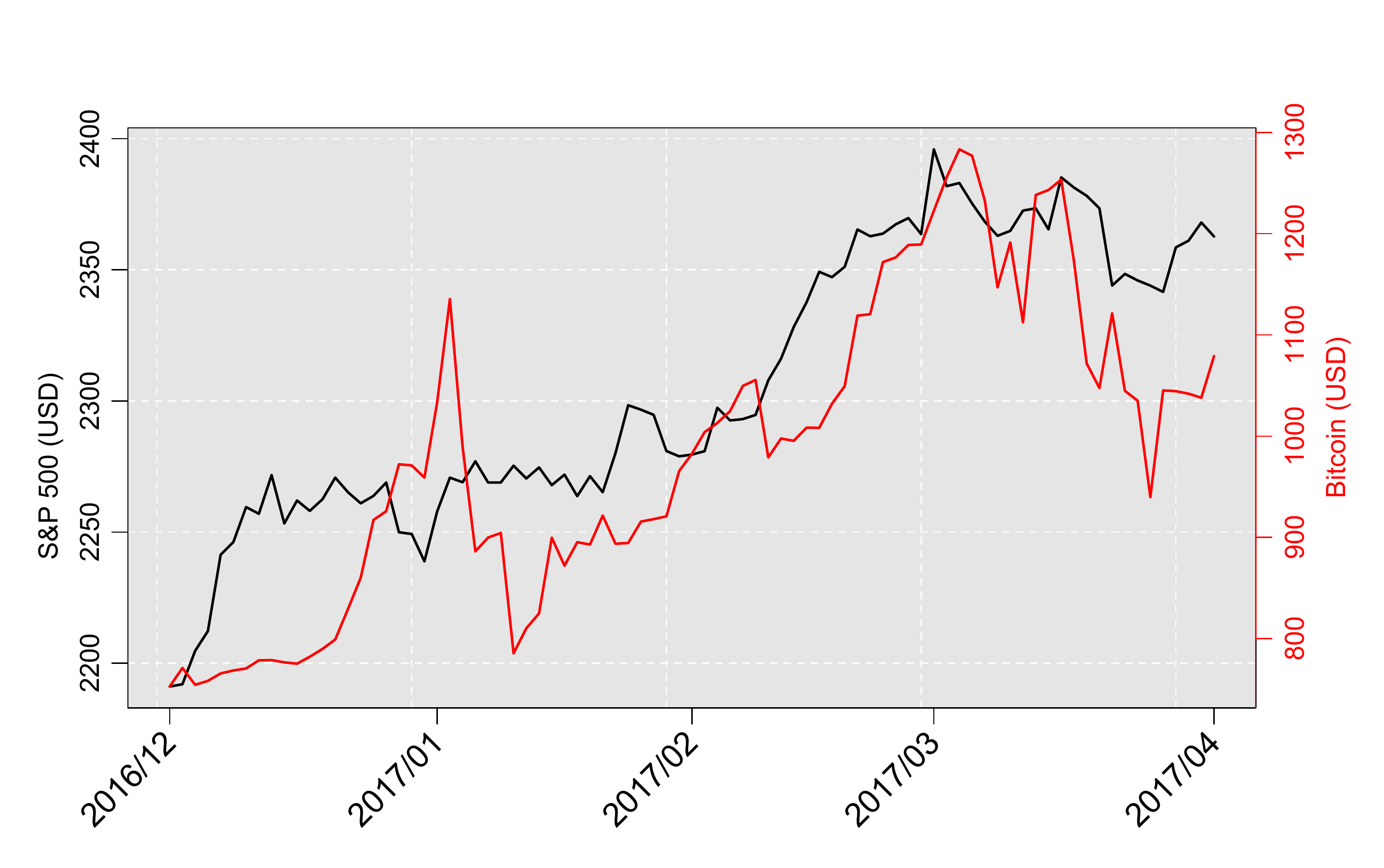}}	 
	\begin{minipage}{10cm}
		\includegraphics[width=10cm]{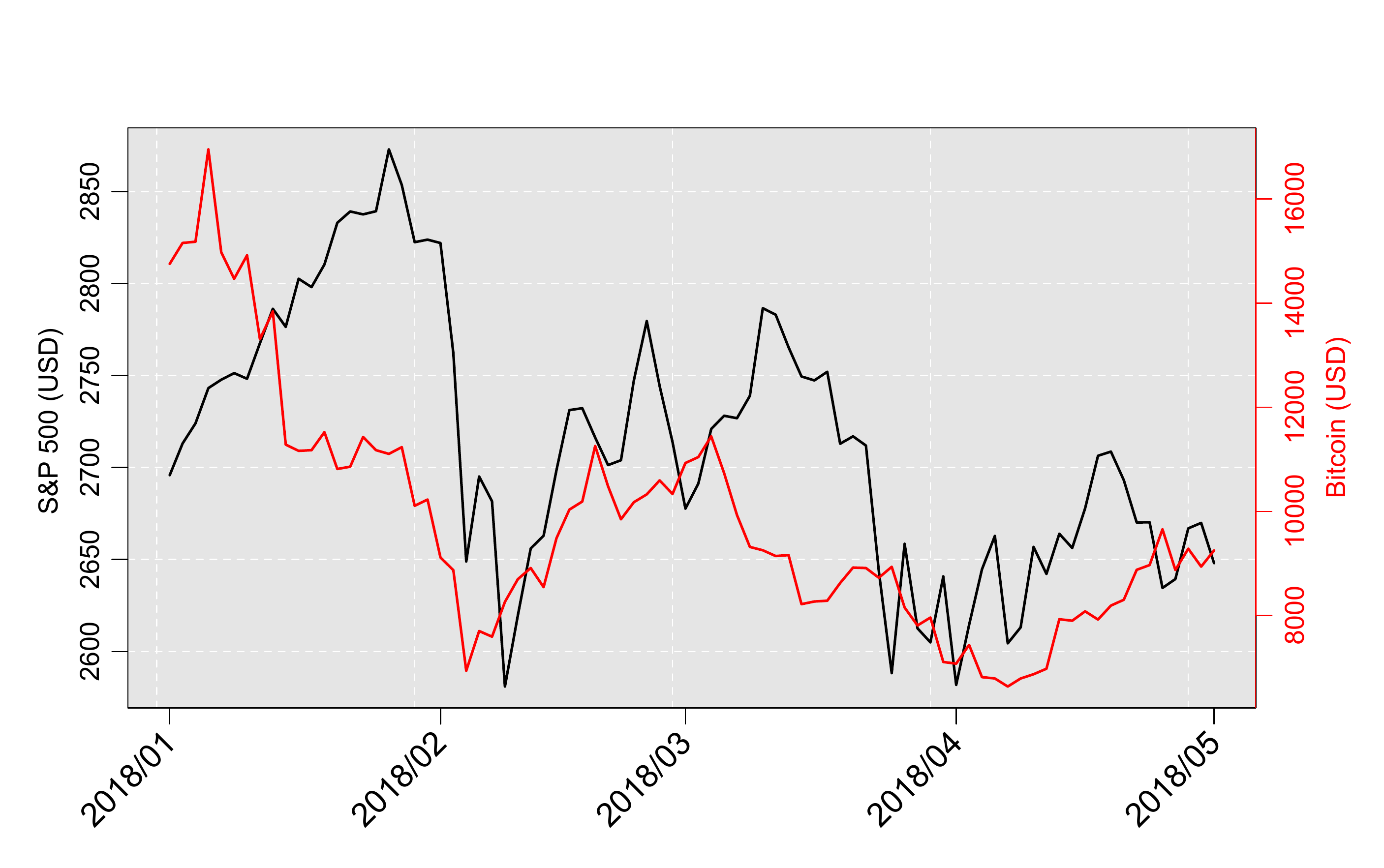}	
	\end{minipage}
	\caption{Closing prices of Standard $\&$ Poors 500 Index and Bitcoin during the time periods 12/01/2016--03/31/2017 and 01/01/2018--04/30/2018, respectively.  The black curve is for Standard $ \&$ Poors 500 Index and the red one is for Bitcoin.}
	\label{fig:realdata}
\end{figure}

The red curve in Figure \ref{fig:real_data} indicates that most of the time the cryptocurrency market and the stock market tend to move independently. There are apparently two periods during which the p-values are below the cut-off point 0.0061, roughly March 2017 and April 2018. Since we use the three-month rolling window right before each date to calculate the p-values, the significantly low p-values in the aforementioned two periods might suggest some nonlinear association between the two markets during the time intervals 12/01/2016--03/31/2017 and 01/01/2018--04/30/2018, respectively. To verify our findings, noticing that Bitcoin is the most representative cryptocurrency and the S\&P 500 Index measures the overall performance of the 500 stocks on its list, we present in the two plots in Figure \ref{fig:realdata} the trend of closing prices of Bitcoin and that of S\&P 500 Index during the periods 12/01/2016--03/31/2017 and 01/01/2018--04/30/2018, respectively. The first plot in Figure \ref{fig:realdata} shows that the trends of the two prices shared striking similarity starting from the middle of January 2017 and both peaked around early March 2017. From the second plot in Figure \ref{fig:realdata}, we see that both the prices of S\&P 500 Index and Bitcoin dropped sharply to the bottom around early Febrary 2018 and then rose to two rekindled peaks followed by continuingly falling to another bottom. Therefore, Figure \ref{fig:realdata} indicates some strong dependency between the two markets in the aforementioned two time intervals and hence demonstrate the effective discoveries of dependence by our normal approximation test for biased-corrected distance correlation.

In addition, to show the robustness of our procedure and choose a reasonable length of rolling window, we also apply four-month and six-month rolling windows before each date $t$ to test the independence between the cryptocurrency market and the stock market. The time series plots of the resulting p-values are presented as the blue curve and the green curve in Figure \ref{fig:real_data}, respectively. From Figure \ref{fig:real_data}, we see that the p-values from using the three different rolling windows (three-month, fourth-month, and six-month) move in a similar fashion. For the four-month rolling window, the p-value cut-off for FDR control at the 10\% level is 0.0053. We observe that the time periods with significantly small p-values by applying four-month rolling window are almost consistent with those by applying three-month rolling window. However, when the six-month rolling window is applied, the p-value cut-off for FDR control at the 10\% level is 0 and hence there is no significant evidence for dependence identified at any time point. This suggests that the long-run dependency between the cryptocurrency market and the stock market might be limited, but there could be some strong association between them in certain special periods. These results show that to test the short-term dependence, the three-month rolling window seems to be a good choice.

\begin{figure}
	\centering
	\includegraphics[scale=0.55]{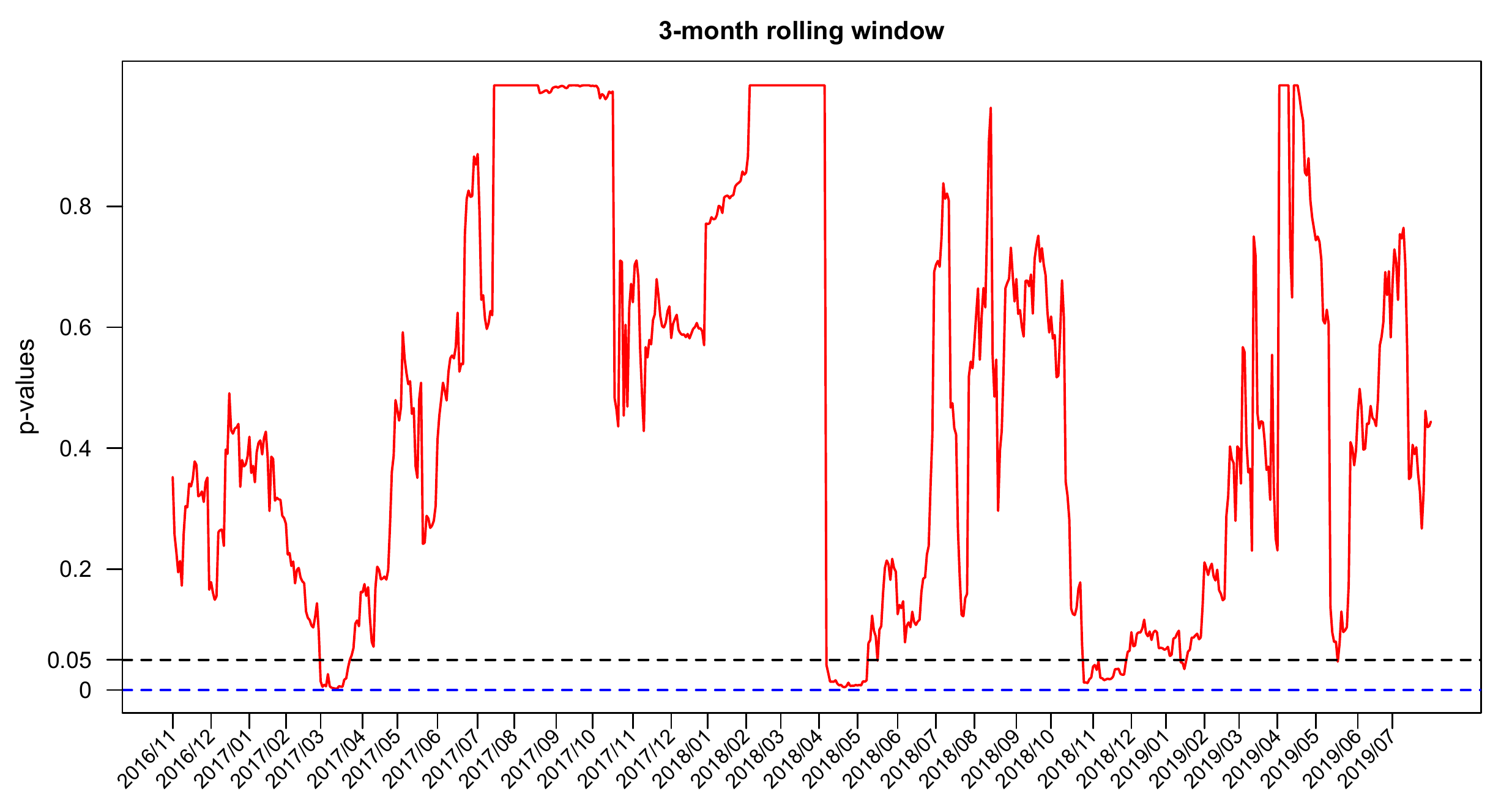}
	\caption{  \it Time series plot of p-values based on RV coefficient from 11/01/2016 to 07/31/2019 using three-month rolling window. 
	}
	\label{fig:RV}
\end{figure}

As a comparison, we conduct the analysis with the rescaled sample distance correlation statistic $T_n$ replaced by the RV coefficient, which measures only the linear dependence between two random vectors. The three-month rolling window is utilized as before. We apply the function `coeffRV' in the R package `FactoMineR' to calculate the p-values of the independence test based on the RV coefficient. The time series plot of the resulting p-values is depicted in Figure \ref{fig:RV}. From Figure \ref{fig:RV}, we see that there are three periods in which the p-values are below the significance level 0.05, while there are four such periods in Figure \ref{fig:real_data} for p-values based on the rescaled sample distance correlation $T_n$ from using three-month rolling window. Moreover, the four periods detected by $T_n$ roughly cover the three periods detected by the RV coefficient. On the other hand, for the p-values based on the RV coefficient, the p-value cut-off for the Benjamini--Hochberg FDR control at the 10\% level is 0, which implies that no significant periods can be discovered with FDR controlled at the 10\% level. However, as mentioned previously, if we use $T_n$ the corresponding p-value cut-off with the three-month rolling window is 0.0061 and two periods, roughly March 2017 and April 2018, are still significant. The effectiveness of these two periods are demonstrated in Figure \ref{fig:realdata}. Therefore, compared to the linear measure of RV coeffcient, the nonlinear dependency measure of rescaled distance correlation is indeed more powerful in this real data application.

\section{Discussions} \label{discussion}

The major contributions of this paper are twofold. First, we have obtained central limit theorems for a rescaled  distance correlation statistic for a pair of high-dimensional random vectors and the associated rates of convergence under the independence when both sample size and dimensionality are diverging. 
Second, we have also developed a general power theory for the sample distance correlation and demonstrated its ability of detecting nonlinear dependence in the regime of moderately high dimensionality. 
These new results shed light on the precise limiting distributions of distance correlation in high dimensions and provide a more complete picture of the asymptotic theory for distance correlation. To prove our main results, Propositions \ref{prop1}--\ref{prop3} in Section \ref{pr-thm4} of Supplementary Material have been developed to help us better understand the moments therein in the high-dimensional setting, which are of independent interest. 

In particular, Theorem \ref{prop-power} unveils that the sample distance correlation is capable of measuring the nonlinear dependence when the dimensionalities of $X$ and $Y$ are diverging. It would be interesting to further investigate the scenario when only one of the dimensionalities tends to infinity and the other one is fixed. 
Moreover, it would also be interesting to extend our asymptotic theory to the conditional or partial distance correlation and investigate more scalable high-dimensional nonparametric inference with theoretical guarantees, for both i.i.d. and time series data settings. These problems are beyond the scope of the current paper and will be interesting topics for future research.

\section*{Acknowledgements}
The authors would like to thank the anonymous referees, an Associate
Editor and the Editor for their constructive comments that improved the paper significantly.

Fan, Gao and Lv's research was supported by NIH Grant 1R01GM131407-01, NSF Grant DMS-1953356, a grant from the Simons Foundation, and Adobe Data Science Research Award. Shao's research was partially suppported by NSFC12031005.

\begin{supplement}
\textbf{Supplement to ``Asymptotic Distributions of High-Dimensional Distance Correlation Inference''}.
The supplement \cite{GFLS2020} contains all the proofs and technical details.  
\end{supplement}

\bibliographystyle{imsart-nameyear} 
\bibliography{references}       

 
 \newpage
 \appendix
 \setcounter{page}{1}
 \setcounter{section}{0}
 \renewcommand{\theequation}{A.\arabic{equation}}
 \renewcommand{\thesubsection}{A.\arabic{subsection}}
 \setcounter{equation}{0}
 
 \begin{center}{\bf \large Supplementary Material to ``Asymptotic Distributions of High-Dimensional Distance Correlation Inference"}
 	
 	\bigskip
 	
 	Lan Gao, Yingying Fan, Jinchi Lv and Qi-Man Shao
 \end{center}

\noindent This Supplementary Material contains all the proofs and technical details. Section \ref{SecA} presents the proofs of the main results in Theorems \ref{thm1}--\ref{prop-power} and Propositions \ref{prop1}--\ref{prop3} in Section \ref{pr-thm4}. We provide the proofs of Propositions \ref{cor0}--\ref{cor2}, some key lemmas with their proofs, and additional technical details in Sections \ref{proof}--\ref{SecF}. In particular, Section \ref{tau} presents the parallel versions of Theorems \ref{thm2} and \ref{thm4} for the case of $ 1/2 < \tau \leq 1 $ and their proofs, while Section \ref{normal-gamma} discusses the connections between the normal approximation for $ T_n  $ and the gamma approximation for $ n \V^* (X, Y) $. Moreover, we provide the proof of the asymptotic normality and associated rates of convergence for $T_R$ in Section \ref{SecF}. 
Throughout the paper, $ C  $ stands for some positive constant whose value may change from line to line.

\section{Proofs of main results} \label{SecA}

\subsection{Proof of Theorem \ref{thm1}}      \label{pr-thm1}
Note that \cite{HS2016} showed that $ \V_n^* (X, Y) $ is a U-statistic. The main idea of our proof is to apply the Hoeffding decomposition for U-statistics and the martingale central limit theorem.
Lemmas \ref{le-consist}--\ref{le-MarCLT} in Sections \ref{SecC.1}--\ref{SecC.4} of Supplementary Material, respectively, draw an outline of the proof.
In particular, Lemma \ref{le-consist} provides the ratio consistency of $ \V_n^* ( X ) $ and $ \V_n^* ( Y ) $. Thus by \eqref{X-con} and \eqref{Y-con}, the denominator of $ T_n $ can be replaced with the corresponding population counterpart in Lemma \ref{le-consist}. In consequence, by Slutsky's lemma it suffices to analyze the limiting distribution of the following random variable
\begin{equation}
\breve{T}_n =   \sqrt{ \frac { n ( n - 1 )} { 2   } }   \frac { \V^*_n  (X, Y)} { \sqrt{ \V^2  ( X )   \V^2  ( Y ) } }. \label{N-Tn}
\end{equation}
Moreover, we have the conclusion in Lemma \ref{le-decom} by the Hoeffding decomposition. In fact, Lemma \ref{le-decom} implies that under the independence of $X $ and $Y$, $ \breve{T}_n $ can be decomposed into two parts $ W_n^{ (1) } ( X, Y ) $ and $ W_n^{ (2) } ( X, Y ) $, where the former is the leading term and the latter is asymptotically negligible. Hence to obtain the limiting distribution of $ \breve{T}_n $, it suffices to focus on $  W_n^{ (1) } ( X, Y ) $ defined in \eqref{W_n_1}. 

Recall the definition of the double-centered distance $ d ( \cdot, \cdot ) $ in \eqref{d-func}. Define $ \zeta_{n, 1} = 0 $ and for $ k \geq 2 $, 
\begin{equation}
\zeta_{n,k} = \sqrt{ \frac {2 } { n (n - 1)  } } \sum_{i= 1}^{ k - 1}  \frac { d (X_i, X_k) d (Y_i, Y_k)} { \sqrt{ \V^2( X ) \V^2 ( Y )  } }. \label{mart-dif}
\end{equation}
It is easy to see that $ W_n^{(1)} (X, Y) = \sum_{k = 1}^n \zeta_{n, k} $. Then by Lemmas \ref{le-martin} and \ref{le-MarCLT}, \eqref{cond1} and \eqref{cond2} directly lead to 
$$
\sum_{k =1}^n \e [ \zeta_{n,k}^2 \vert \mathscr{F}_{k-1} ]   \rightarrow 1 \  \mbox{ in probability}
$$
with $\mathscr{F}_{k}$ a $\sigma$-algebra defined in Lemma \ref{le-martin}, and for any $\varepsilon > 0$, 
$$
\sum_{ k = 1}^n \e [ \zeta_{n,k}^2  \mathbf{1} \{  | \zeta_{n,k} | > \varepsilon  \} ]  \rightarrow 0. 
$$
Therefore, by the Lindeberg-type central limit theorem for martingales (see, for example, \cite{Brown1971}), we can obtain $ W_n^{(1)} (X, Y) \xrightarrow { \mathscr{D} } N(0, 1) $. This completes the proof of Theorem \ref{thm1}.

\subsection{Proof of Theorem \ref{thm2}}  \label{pf-thm2}
The main idea of the proof is based on the conclusion of Theorem \ref{thm4}. In view of the definitions of $ E_x  $ and $ L_{ x, \tau } $, by the Cauchy--Schwarz inequality we can obtain that 
\begin{align*}
E_x \geq \frac { B_X^{ - 2 \tau } L_{ x, \tau }^{ ( 2 + \tau  ) / ( 1 + \tau ) } } {  \{ \e [ ( X_1^T X_2 )^2 ] \}^2  } \geq \frac { B_X^{  - 2 \tau } L_{ x, \tau }  } { \e [ ( X_1^T X_2 )^2 ] }.
\end{align*}
In the same manner,  we can deduce
\begin{align*}
E_y \geq \frac { B_Y^{  - 2 \tau } L_{ y, \tau }  } { \e [ ( Y_1^T Y_2 )^2 ] }.
\end{align*}
Note that $ p + q \rightarrow \infty $ implies that at least one of $ p  $ and $ q  $ tends to infinity. First let us assume that both $ p \rightarrow \infty $ and $ q \rightarrow  \infty $. Then by assumption, we have $ E_x \rightarrow 0 $ and $ E_y \rightarrow 0 $. Thus for sufficiently large $ p $ and $ q $, it holds that 
\begin{align*}
B_x^{ - 2 \tau  }  L_{ x, \tau } / \e [ ( X_1 ^T X_2  )^2 ] \leq \frac {1} { 18 } \ \text{ and } \ B_Y^{ - 2 \tau  } L_{ y, \tau } /  \e [ ( Y_1^T Y_2 )^2 ]  \leq \frac {1} {  18 }.
\end{align*}
It follows from Theorem \ref{thm4} that if \eqref{co1} holds, $ E_x \rightarrow 0 $, and $ E_y \rightarrow 0 $, then we have 
\begin{align*}
\sup \limits_{ x\in \R } | \mathbb{P} ( T_n \leq x  ) - \Phi ( x ) |  \rightarrow 0 
\end{align*}
with $ \Phi (x) $ the standard normal distribution function, which yields $ T_n \stackrel{\mathscr{D}}{\rightarrow } N( 0, 1 ) $.

We now consider the scenario when only one of $ p $ and $ q $ tends to infinity. Without loss of generality, assume that $ p $ is bounded and $ q \rightarrow \infty $. Then by assumption, we have $ E_y \rightarrow 0 $. In addition, note that 
$L_{x, \tau }  \geq  \big(  \e [ ( X_1^T X_2 )^2 ] \big)^{ 1 + \tau }$.
Thus it follows from \eqref{co1} that 
\begin{align*}
\frac { n ^{ - \tau  } L_{ y, \tau } } { \big\{ \e [ ( Y_1^T Y_2 )^2 ] \big\}^{ 1 + \tau } } \rightarrow 0.
\end{align*}
Consequently, an application of bound \eqref{rate-single} 
results in 
\begin{align*}
\sup \limits_{ x\in \R } | \mathbb{P} ( T_n \leq x  ) - \Phi ( x ) |  \rightarrow 0,
\end{align*}
which concludes the proof of Theorem \ref{thm2}.

\subsection{Proof of Theorem \ref{thm3}} \label{SecA.3}
The key ingredient of the proof is to replace the denominator with the population counterpart and apply the convergence rate in the martingale central limit theorem. In light of the definition in \eqref{N-Tn}, we can write
\begin{align*}
|  \mathbb{P} ( T_n \leq x ) - \Phi (x) |  = \Big| \mathbb{P} \Big( \breve{T}_n \cdot \sqrt{ \frac { \V^2 ( X  ) \V^2 ( Y ) } { \V_n^* ( X ) \V_n^* ( Y )  }   } \leq x \Big) - \Phi ( x ) \Big|.
\end{align*} 
Note that Lemma \ref{le-consist} entails that $ \V_n^* (X)  / \V^2 (X) $ and $ \V_n^* (Y)  / \V^2 (Y) $ converge to one in probability. Thus we can relate the distance between $ \mathbb P ( T_n \leq x ) $ and $ \Phi (x) $ to that between $ \mathbb {P} ( \breve{T}_n \leq x ) $ and $ \Phi (x) $. Specifically, for small quantities $ \gamma_1 > 0 $ and $ \gamma_2 > 0 $ it holds that 
\begin{align}
|  \mathbb{P} ( T_n \leq x ) - \Phi (x) | & \leq P_1 + P_2 + \mathbb{P} \Big( \Big| \frac { \V_n^* ( X ) } { \V^2 ( X ) }  - 1 \Big| > \gamma_1  \Big)  + \mathbb{P} \Big( \Big| \frac { \V_n^* ( Y ) } { \V^2 ( Y ) }  - 1 \Big| > \gamma_2  \Big),  \label{bound_all}
\end{align} 
where
\begin{align*}
P_1 & =  \big| \mathbb{P} \big( \breve{T}_n \leq x ( 1 + \gamma_1 ) ( 1 + \gamma_2 ) \big) - \Phi (x) \big|, \nn \\
P_2 & =  \big| \mathbb{P} \big( \breve{T}_n \leq x ( 1 - \gamma_1 ) ( 1 -\gamma_2 ) \big) - \Phi (x) \big|.
\end{align*}
Let us choose
$$ 
\gamma_1 =  \Big\{  \frac { \e [ | d ( X_1, X_2 ) |^{ 2 + 2 \tau } ] } { n^{ \tau }  [ \V^2 ( X ) ]^{ 1 + \tau  }  }  \Big\}^{ 1/ ( 2 + \tau ) } , \quad  \gamma_2 = \Big\{ \frac {  \e [ | d ( Y_1, Y_2 ) |^{ 2 + 2 \tau } ] } { n^{ \tau }    [ \V^2 ( Y ) ]^{ 1 + \tau }  }  \Big\}^{ 1 / ( 2 + \tau ) }.
$$
Without loss of generality, assume that $ \gamma_1 \leq 1/2  $ and $ \gamma_2 \leq 1/2 $. Otherwise since 
\begin{equation}
\e [ | d ( X_1, X_2 ) |^{ 2 + 2 \tau } ]  \geq \big\{ \e [ d^2 ( X_1, X_2 ) ]  \big\}^{ 1 + \tau }  = [  \V^2 ( X )  ]^{ 1 + \tau }  \label{Jensen}
\end{equation}
and similar result holds for $Y$,
we have 
\begin{align*}
\frac {  \e [ | d (X_1, X_2 ) |^{2 + 2 \tau } ]  \e [  | d (Y_1, Y_2) |^{2 + 2 \tau } ] } { n^{ \tau  }  [ \V^2 ( X ) \V^2 ( Y ) ]^{ 1 + \tau  }    } & \geq  \max \Big\{    \frac { \e [ | d ( X_1, X_2 ) |^{ 2 + 2 \tau } ] } { n^{ \tau }  [ \V^2 ( X ) ]^{ 1 + \tau  }  } ,  \frac {  \e [ | d ( Y_1, Y_2 ) |^{ 2 + 2 \tau } ] } { n^{ \tau }    [ \V^2 ( Y ) ]^{ 1 + \tau }  } \Big\} \nn \\
& \geq 2^{ - ( 2 + \tau ) }
\end{align*}
and thus the desired result \eqref{rate-general} is trivial. 

Now we bound the four terms on the right hand side of \eqref{bound_all}. By \eqref{bound-1}, it holds that 
\begin{align}
\mathbb{P} \Big( \Big| \frac { \V_n^* ( X ) } { \V^2 ( X ) }  - 1 \Big| > \gamma_1  \Big) & \leq  \frac { C \e [ | d ( X_1, X_2 ) |^{ 2 + 2 \tau } ] } { n^{ \tau } \gamma_1^{ 1 + \tau }  [ \V^2 ( X ) ]^{ 1 + \tau } } =  C  \Big\{  \frac { \e [ | d ( X_1, X_2 ) |^{ 2 + 2 \tau } ] } { n^{ \tau }  [ \V^2 ( X ) ]^{ 1 + \tau  }  }  \Big\}^{ 1/ ( 2 + \tau ) }   \label{P3}
\end{align} 
and similarly, 
\begin{align}
\mathbb{P} \Big( \Big| \frac { \V_n^* ( Y ) } { \V^2 ( Y  ) }  - 1 \Big| >  \gamma_2  \Big) \leq  C \Big\{ \frac {  \e [ | d ( Y_1, Y_2 ) |^{ 2 + 2 \tau } ] } { n^{ \tau }    [ \V^2 ( Y ) ]^{ 1 + \tau }  }  \Big\}^{ 1 / ( 2 + \tau ) }.  \label{P4}
\end{align}

Then we deal with term $ P_1 $. By symmetry, term $ P_2 $ shares the same bound as term $ P_1 $. By Lemma \ref{le-decom}, $ \breve{T}_n  $ can be decomposed into two parts, one being the dominating martingale array and the other being an asymptotically negligible error term. In details, for $ 0 < \gamma_3  = n^{ - 1/3 } / 4  \leq 1 /4  $ we have 
\begin{align*}
P_1 \leq P_{11} + P_{12} + \mathbb{P} ( | W_{n}^{ (2) } ( X, Y ) | > \gamma_3 ), 
\end{align*}
where
\begin{align*}
P_{11} & = \big| \mathbb{P} \big( W_{n}^{ ( 1 ) } ( X, Y ) \leq x ( 1 + \gamma_1 ) ( 1 + \gamma_2 ) - \gamma_3  \big) - \Phi ( x ) \big|, \nn \\
P_{12} & = \big| \mathbb{P} \big( W_{n}^{ ( 1 ) } ( X, Y ) \leq x ( 1 + \gamma_1 ) ( 1 + \gamma_2 ) + \gamma_3  \big) - \Phi ( x ) \big|. 
\end{align*} 
It follows from Lemma \ref{le-decom} that
\begin{align}
\mathbb{P} ( | W_{n}^{ (2) } ( X, Y ) | > \gamma_3 ) \leq \frac { 1 } { n \gamma_3^2 } \leq 16  n^{ - 1 / 3 }. \label{P13}
\end{align}
Since terms $ P_{11} $ and $ P_{12} $ share the same bound, it suffices to show the analysis for term $ P_{11} $. It holds that 
\begin{align}
P_{11} \leq \sup_{x \in \R } \big| \mathbb{P} \big( W_n^{ (1) } ( X, Y ) \leq x  \big)  - \Phi ( x ) \big| + \sup_{ x \in \R }  \big| \Phi [ x ( 1 + \gamma_1  ) ( 1 + \gamma_2 ) - \gamma_3  ] - \Phi (x) \big|.  \label{bound_P11}
\end{align} 

Observe that by definitions, we have $ \gamma_1 \leq  1/2  $, $ \gamma_2 \leq  1/ 2 $, and $ \gamma_3 \leq 1/4  $. When $ | x | \leq 2 $, it is easy to see that 
\begin{align*}
| \Phi [ x ( 1 + \gamma_1  ) ( 1 + \gamma_2 ) - \gamma_3  ] - \Phi (x) | \leq C ( \gamma_1 + \gamma_2 + \gamma_3 ).
\end{align*} 
When $ | x | > 2 $, we have $ | x ( 1 + \gamma_1 ) ( 1 + \gamma_2 )  | / 2 >  \gamma_3 $ and thus 
\begin{align*}
| \Phi [ x ( 1 + \gamma_1  ) ( 1 + \gamma_2 ) - \gamma_3  ] - \Phi (x) | & \leq C ( x \gamma_1 + x \gamma_2 + \gamma_3 ) e^{ -     x^2 / 128  } \nn \\
& \leq C ( \gamma_1 + \gamma_2 + \gamma_3 ).
\end{align*} 
Consequently, it follows that 
\begin{align}
\sup_{ x \in \R } | \Phi [ x ( 1 + \gamma_1  ) ( 1 + \gamma_2 ) + \gamma_3  ] - \Phi (x) | \leq C  ( \gamma_1 + \gamma_2 + \gamma_3 ). \label{P112}
\end{align}

As for the bound of $  \big| \mathbb{P} \big( W_n^{ (1) } ( X, Y ) \leq x  \big)  - \Phi ( x ) \big| $, note that $ W_n^{ (1) } (X, Y) = \sum_{ k = 1 }^n \zeta_{n,k} $ and Lemma \ref{le-martin} states that $  \{ (\zeta_{n, k}, \mathscr{F}_k) , k \geq 1  \} $ is a martingale difference array under the independence of $X$ and $Y$. Hence by Theorem 1 in \citet{H1988} on the convergence rate of the martingale central limit theorem and Lemma \ref{le-MarCLT}, we can obtain
\begin{align}
& \sup_{x \in \R } | \mathbb{P} \big( W_n^{ (1) } ( X, Y ) \leq x  \big)  - \Phi ( x ) \big| \nn \\
& \leq C  \Big\{ \textstyle \sum_{ k  = 1 }^n  \e [ | \zeta_{n, k} |^{ 2 + 2 \tau }  ] + \e \big(  \big| \sum_{k = 1}^n \e [ \zeta_{n, k}^2 \vert \mathscr{F}_{ k - 1} ] - 1 \big|^{ 1 + \tau }  \big)  \Big\}^{ 1 / (  3 + 2 \tau ) } \nn \\
& \leq C  \bigg\{   \Big ( \frac  {      \e [ g (X_1,X_2, X_3, X_4 ) ]  \e [ g( Y_1, Y_2, Y_3, Y_4 ) ]   }  {  [  \V^2 (X) \V^2 (Y)  ] ^2 }  \Big)^{   ( 1 + \tau  ) / 2   } \nn \\
& \hspace{3cm}   + \frac {  \e [ | d (X_1, X_2 ) |^{2 + 2 \tau } ]  \e [  | d (Y_1, Y_2) |^{2 + 2 \tau } ] } { n^{ \tau  }  [ \V^2 (X) \V^2 (Y) ]^{ 1 + \tau  }    } \bigg\}^{ 1/ ( 3 + 2 \tau ) }.   \label{P111}
\end{align}
By \eqref{Jensen} and $ 0 < \tau \leq 1 $, it holds that 
\begin{align*}
\frac {  \e [ | d (X_1, X_2 ) |^{2 + 2 \tau } ]  \e [  | d (Y_1, Y_2) |^{2 + 2 \tau } ] } { n^{ \tau  }  [ \V^2 (X) \V^2 (Y) ]^{ 1 + \tau  }    }  \geq \frac {  \e [ | d (X_1, X_2 ) |^{2 + 2 \tau } ]   } { n^{ \tau  }  [ \V^2 (X)   ]^{ 1 + \tau  }    }
\end{align*}
and 
\begin{align*}
\Big\{   \frac {  \e [ | d (X_1, X_2 ) |^{2 + 2 \tau } ]   } { n^{ \tau  }  [ \V^2 (X)   ]^{ 1 + \tau  }    } \Big\}^{ 1/ (2 + \tau ) } \geq  n^{ - \tau / ( 2 + \tau  ) } \geq n^{ - 1 /3  }.
\end{align*}  
Finally, the desired result \eqref{rate-general} can be derived by plugging in \eqref{P3}--\eqref{P111} and noting that all the error terms can be absorbed into \eqref{P111}. This completes the proof of Theorem \ref{thm3}.

\subsection{Proof of Theorem \ref{thm4}} \label{pr-thm4}
The proof is mainly based on the conclusion of Theorem \ref{thm3}. It is quite challenging to calculate the exact form of the moments that appear in conditions \eqref{cond1} and \eqref{cond2}. Nevertheless, the bounds of these moments can be worked out in concise form under some general conditions. These bounds are summarized in the following three propositions, respectively.

\begin{proposition} \label{prop1}
	If $ \e [ \Vert X \Vert ^{ 4 + 4 \tau  }  ]  < \infty $ for some constant $ \tau > 0 $, then there exists some absolute positive constant $ C_{\tau}$ such that 
	\begin{align}
	\e ( | d ( X_1, X_2 ) |^{ 2 + 2 \tau } ) \leq  C_{\tau}    B_X ^{ - (  1 + \tau ) }  L_{x, \tau}  .  \label{d-mo}
	\end{align} 
\end{proposition}

\begin{proposition}   \label{prop2}
	If $ \e [ \Vert X \Vert ^{ 4 + 4 \tau }   ]  < \infty $ for some constant $ 0 < \tau \leq 1/2 $, then  it holds that 
	\begin{align}
	\big|  \V^2 ( X )  - B_X^{ - 1} \e [ ( X_1 ^T X_2 )^2 ]   \big| \leq 9  B_X^{ - (1 + 2 \tau ) }  L _{x, \tau} .     \label{V}
	\end{align}
\end{proposition}

\begin{proposition} \label{prop3}
	If $ \e [ \Vert X \Vert ^{ 4 + 4 \tau }   ]  < \infty $ for some constant $ 0 < \tau \leq 1/2 $, then there exists some absolute positive constant $ C  $ such that 
	\begin{align}
	\big| \e [ g ( X_1, X_2, X_3, X_4 ) ] \big|  \leq   B_X^{ - 2 } \e [ ( X_1^T \Sigma_x X_2 )^2  ]  + C B_X^{ - ( 2 + 2 \tau ) }  L_{ x, \tau }^{  (2 + \tau) / ( 1 + \tau  )    }.  \label{g1}
	\end{align}
\end{proposition}     

The proofs of Propositions \ref{prop1}--\ref{prop3} are presented in Sections \ref{pf-prop1}--\ref{pf-prop3}, respectively. We now proceed with the proof of Theorem \ref{thm4}. Note that condition \eqref{co2} entails that 
\begin{align*}
9  B_X^{ - ( 1 + 2 \tau ) } L_{ x, \tau }   \leq  \frac 1 2 B_X^{ - 1  }  \e [ ( X_1^T X_2 )^2  ] \ \text{ and } \ 
9  B_Y^{ - ( 1 + 2 \tau ) } L_{ y, \tau }   \leq  \frac 1 2 B_Y^{ - 1  }  \e [ ( Y_1^T Y_2 )^2  ].
\end{align*}
Therefore, it follows from Proposition \ref{prop2} that 
\begin{align}
\V^2 ( X ) \geq  \frac 1 2 B_X^{ -1  } \e [  ( X_1 ^T X_2 )^2 ] \ \text{ and } \    \V^2 ( Y )  \geq \frac 1 2 B_Y^{ - 1 } \e [ ( Y_1^T Y_2  )^2 ],  \label{lb-V2}
\end{align}
which together with Propositions \ref{prop1} and \ref{prop3} yield the desired results \eqref{rate1} by Theorem \ref{thm3}. This concludes the proof of Theorem \ref{thm4}.

\subsection{Proof of Theorem \ref{thm-power}} \label{pf-thm-power}
{
	Recall that 
	$$ T_n =  \sqrt{ \frac { n (n - 1) } {2} } \frac { \V_n^* (X, Y) } { \sqrt{\V_n^* (X) \V_n^* (Y) }} $$
	and it has been proved in Lemma \ref{le-consist} in Section  \ref{SecC.1} that under condition \eqref{cond1}, we have 
	$   \V_n^* (X) / \V^2 (X)  \to 1 $ and $  \V_n^* (Y) / \V^2 (Y) \to 1$ in probability. Thus it suffices to show that for any arbitrarily large constant $C > 0 $,
	$$
	\breve{T}_n :=  \sqrt{\frac{ n (n - 1) }{2}}  \frac { \V_n^* (X, Y) } { \sqrt{\V^2 (X) \V^2 (Y) }} > C  \ \mbox{ with asymptotic probability 1}.
	$$
	Observe that 
	\begin{align*}
	\Big| \breve{T}_n  - \sqrt{ \frac { n ( n - 1) } { 2 } } \Rc^2 (X, Y)  \Big|   =    \sqrt{ \frac { n ( n - 1) } { 2 } } \frac {  | \V_n^*(X, Y) - \V^2 (X, Y) | } { \sqrt{ \V^2 (X) \V^2 (Y) } } .
	\end{align*} 
	
	It follows from \eqref{h-b1}, \eqref{ex-h-d} and Proposition \ref{prop1} that there exists some absolute positive constant $C$ such that 
	\begin{align*}
	\e [ ( \V_n^* (X, Y) - \V^2 (X, Y) )^2  ]&  \leq C n^{-1}  \e [ h^2 ( (X_1, Y_1 ), (X_2, Y_2), (X_3, Y_3), (X_4, Y_4) ) ] \\
	& \leq C n^{-1} ( \e [ d^4 (X_1, X_2) ] \e [ d^4 (Y_1, Y_2) ] )^{1/2} \\
	& \leq C n^{-1} B_X^{-1} B_Y^{-1} L_{x, 1}^{1/2} L_{y, 1}^{1/2}.
	\end{align*}
	Therefore, if $ \sqrt{n} \V^2 (X, Y) / \big(  B_X^{-1/2} B_Y^{-1/2} L_{x, 1/2}^{1/4} L_{y, 1}^{1/4} \big)  \to \infty$, it holds that 
	\begin{align*}
	\frac {  | \V_n^*(X, Y) - \V^2 (X, Y) | } { \sqrt{ \V^2 (X) \V^2 (Y) } }  \Big/ \Rc^2 (X, Y) \to 0 \quad \mbox{in probability}.
	\end{align*}
	This together with $ n \Rc^2 (X, Y) \to \infty $ yields for any arbitrarily large constant $C > 0$, $ \mathbbm{P} ( \breve{T}_n > C) \rightarrow 1 $ and hence  as $ \mathbbm{P} ( T_n > C ) \rightarrow 1$, which completes the proof of Theorem \ref{thm-power}.
}	

\subsection{Proof of Theorem \ref{prop-power}} \label{pf-prop-power}

The main ingredient of the proof is bounding $ \V^2 (X, Y) $ using the decomposition developed in Lemma \ref{le-dcov} in Section \ref{SecC.10}. We will calculate the orders of terms $  I_i, 1 \leq i \leq 5$, introduced in Lemma \ref{le-dcov}. Let us begin with the first term 
\begin{align*}
I_1 =  \frac 1 4 B_X^{1/2} B_Y^{1/2} \big( \e [ W_{12} V_{12} ] -  2 \e [ W_{12} V_{13} ]  \big),
\end{align*}
where $ W_{12} = B_X^{-1} ( \Vert X_1 - X_2 \Vert^2 - B_X ) $ and $ V_{12}  = B_Y^{-1} ( \Vert Y_1 - Y_2 \Vert^2 - B_Y )  $. Denote by $ \widetilde{Y}_1 = Y_1 - \e Y = ( Y_{1,1} - \e Y_{1, 1}, \cdots, Y_{1, p} - \e Y_{1, p} ) ^T $ and $ \widetilde{Y}_2 =  Y_2 - \e Y = (Y_{2, 1}- \e Y_{2, p}, \cdots, Y_{2, p} - \e Y_{2, p})^T $ the centered random variables, and define 
\begin{align*}
\alpha_1 (X) =  \Vert X \Vert^2  - \e \Vert X \Vert^2 , \quad  \alpha_2 (X_1, X_2 ) = {X}_1^T {X}_2, \\
\beta_1 (Y)  = \Vert Y \Vert^2 - \e \Vert Y \Vert^2 , \quad \beta_2 (Y_1, Y_2)  = \widetilde{Y}_1^T \widetilde{Y}_2.
\end{align*}
Since $ \e [ \alpha_1 (X) ]  = \e[ \beta_1 (Y) ] = 0$ and $ \e [ \alpha_2 (X_1, X_2) ]  = \e [ \beta_2 (Y_1, Y_2) ]= 0$, it holds that 
\begin{align*}
\e [ W_{12} V_{12} ] & = \e \Big[ \big(\alpha_1 (X_1) + \alpha_1 (X_2) - 2 \alpha_2 (X_1, X_2) \big)  \big( \beta_1 (Y_1) + \beta_1 (Y_2)  - 2 \beta_2 (Y_1, Y_2)\big)\Big]  \\
& = 2 \e [ \alpha_1 (X_1) \beta_1 (Y_1)] + 4 \e [ \alpha_2 (X_1, X_2) \beta_2 (Y_1, Y_2)].
\end{align*}
Similarly, we have $ \e [ W_{12} V_{13}]  = 2 \e [ \alpha_1 (X_1) \beta_1 (Y_1) ]$. 
Thus it follows that 
\begin{align*}
I_1 =   4  \e [ \alpha_2 (X_1, X_2) \beta_2 (Y_1, Y_2)] =  4 \sum_{i, j = 1}^p  ( \cov(X_{1, i}, Y_{1, j}) )^2 .
\end{align*}
Observe that under the symmetry assumptions, there is no linear dependency between $X$ and $Y$; that is, $ \cov(X_{1, i}, Y_{1, j}) = 0 $ for each $ 1 \leq i, j \leq p $. This together with the representation of $I_1$ above entails that $ I_1  = 0 $.

We now consider the second term $I_2$. Using similar arguments but much more tedious calculations, we can obtain
\begin{align*}
I_2  &= \frac 1 4 B_X^{ - 1/2} B_Y^{ - 3/2}   \Big( 2 \e [ \alpha_2 ( X_1, X_2  ) \beta_2^2 ( Y_1, Y_2 ) ] +  \e [ \alpha_2 (X_1, X_2) \beta_1 (Y_1) \beta_1 (Y_2) ] \\
& \hspace{3cm}    - 4  \e  [ \alpha_2 (X_1, X_2) \beta_1 (Y_1) \beta_2 (Y_1, Y_2) ]     \Big) \\
& \quad  +  \frac 1 4 B_X^{ - 3/2 } B_Y^{ - 1/2} \Big(  2 \e [ \beta_2 (Y_1, Y_2) \alpha_2^2 (X_1, X_2)  ]  + \e [ \beta_2 (Y_1, Y_2)  \alpha_1 (X_1) \alpha_1 (X_2) ] \\
& \hspace{3cm} - 4 \e [ \beta_2 (Y_1, Y_2 ) \alpha_1 (X_1) \alpha_2 (X_1, X_2) ] \Big).
\end{align*}
By assumption, we have $ c_2 p \leq B_X \leq c_1 p $ and $ c_2 p \leq B_Y \leq c_1 p $. Since $X$ has a symmetric distribution and $ Y_{1, j}  = g_{j} (X_{1, j}) $ with $ g_j (x), 1 \leq j \leq p $, symmetric functions, it holds that 
\begin{align*}
\e [ \alpha_2 ( X_1, X_2 )  \beta_2^2 (Y_1, Y_2) ]  = \e [ (X_1^T X_2) ( \widetilde{Y}_1^T \widetilde{Y}_2 )^2 ] = \e [ ( - X_1^T X_2 ) (\widetilde{Y}_1^T \widetilde{Y}_2 )^2   ] = 0 .
\end{align*}
Similarly, with the symmetry assumptions we can show that $ \e [ \alpha_2 (X_1, X_2) \beta_1 (Y_1) \beta_1 (Y_2) ] = 0 $, $ \e  [ \alpha_2 (X_1, X_2) \beta_1 (Y_1) \beta_2 (Y_1, Y_2) ]   =0  $, and $
\e [ \beta_2 (Y_1, Y_2 ) \alpha_1 (X_1) \alpha_2 (X_1, X_2) ]  = 0 
$. Moreover, it holds that 
\begin{align*}
\e [ \beta_2 (Y_1, Y_2) \alpha_2^2 (X_1, X_2)  ]  &= \sum_{i, j, k = 1}^p \big( \e [  {X}_{1, i}  {X}_{1, j} \widetilde{Y}_{1, k} ] \big)^2 \geq  0, \\       
\e [  \beta_2 (Y_1, Y_2)  \alpha_1 (X_1) \alpha_1 (X_2)   ] &  = \sum_{i= 1}^p \Big( \sum_{j=1}^p \e \big[ \widetilde{Y}_{1, i} ( X_{1, j}^2 - \e X_{1, j}^2  ) \big] \Big)^2  \geq 0.
\end{align*}
Thus it follows that $I_2 \geq 0 $.

Let us proceed with terms $ I_3 $ and $ I_4 $. By some tedious 
calculations, we can deduce that 
\begin{align*}
I_3 &=   \frac {1} { 8 } B_X^{-1/2} B_Y^{  - 5/2  }  \Big(  4 \e [ \alpha_2 (X_1, X_2) \beta_2^3 (Y_1, Y_2) ]  + 6 \e [ \alpha_2 (X_1, X_2)  \beta_1(Y_1) \beta_1 (Y_2) \beta_2 (Y_1, Y_2) ]   \\
& \hspace{3cm} + 6 \e [ \alpha_2 (X_1, X_2) \beta_1^2 (Y_1) \beta_2 (Y_1, Y_2)]   - 3 \e [ \alpha_2 (X_1, X_2)  \beta_1 (Y_1) \beta_1^2 (Y_2) ] \\
& \hspace{3cm}  - 12  \e [ \alpha_2 (X_1, X_2)  \beta_1 (Y_1) \beta_2^2 (Y_1, Y_2) ] \Big) \\
& \quad + \frac 1 8 B_X^{ - 5/ 2} B_Y^{ - 1/2}  \Big(   4 \e [ \beta_{2}(Y_1, Y_2) \alpha_2^3 (X_1, X_2) ] +  6 \e [ \beta_{2}(Y_1, Y_2) \alpha_1 (X_1) \alpha_1 (X_2) \alpha_2 (X_1, X_2) ]\\
& \hspace{3cm} + 6 \e [ \beta_{2}(Y_1, Y_2) \alpha_1^2 (X_1) \alpha_2 (X_1, X_2)  ] - 3 \e [ \beta_{2}(Y_1, Y_2) \alpha_1(X_1) \alpha_1^2 (X_2) ]\\
& \hspace{3cm}  - 12 \e [ \beta_{2}(Y_1, Y_2)  \alpha_1 (X_1) \alpha_2^2 ( X_1 , X_2 ) ]
\Big) 
\end{align*}
and 
\begin{align*}
I_4 & = \frac{1}{16} B_X^{ - 3/2} B_Y^{ - 3/2} \Big( 4 \e [ \alpha_2^2 (X_1, X_2)  \beta_2^2 (Y_1, Y_2)]  +  \big(  \e [ \alpha_1 (X) \beta_1(Y) ] \big)^2  \\
&  \hspace{1cm} +  8 \e [ \alpha_2^2 (X_1, X_2)  \beta_2^2 ( Y_1, Y_3 ) ]  +  4 \e [ \alpha_2^2 ( X_1, X_2 ) ] \e [ \beta_2^2 (Y_1, Y_2) ] \\
& \hspace{1cm} +  2 \e [ \alpha_2^2 (X_1, X_2) \beta_1(Y_1) \beta_1 (Y_2) ]  - 8 \e [ \alpha_2^2 (X_1, X_2) \beta_1 (Y_1) \beta_2 (Y_1, Y_2) ] \\
& \hspace{1cm} + 2 \e [ \alpha_1 (X_1) \alpha_1 (X_2) \beta_2^2 (Y_1, Y_2) ] - 4 \e [ \alpha_1 (X_1) \alpha_1(X_2) \beta_1 (Y_1) \beta_2 (Y_1, Y_2) ] \\
& \hspace{1cm}  - 8 \e [ \alpha_1 (X_1) \alpha_2 (X_1, X_2)  \beta_2^2 (Y_1, Y_2) ] - 4 \e [ \alpha_1 (X_1) \alpha_2 (X_1, X_2) \beta_1 (Y_1) \beta_1 (Y_2)] \\
& \hspace{1cm} + 8 \e [ \alpha_1 (X_1) \alpha_2 (X_1, X_2) \beta_1 (Y_1) \beta_2 (Y_1, Y_2) ]  + 8 \e [ \alpha_1 (X_1) \alpha_2 (X_1, X_2) \beta_1 (Y_2) \beta_2 (Y_1, Y_2)   ]   \\
& \hspace{1cm} - 8 \e [ \alpha_2^2 (X_1, X_2) \beta_1 (Y_3) \beta_2 (Y_1, Y_3)] - 8 \e [ \alpha_1 (X_2) \alpha_2(X_1, X_2) \beta_2^2 (Y_1, Y_3) ]  \\
& \hspace{1cm} + 8 \e [ \alpha_1 (X_2) \alpha_2 (X_1, X_2) \beta_1 (Y_3) \beta_2 ( Y_1, Y_3 ) ] \Big) .   
\end{align*}
A useful observation is that under the assumptions that $X_1$ has a symmetric distribution and $ g_j (x)$ with $1 \leq j \leq p$ are symmetric functions, many terms in $ I_3 $ and $ I_4 $ above in fact become zero. In particular, we can show that  
\begin{align*}
I_3  & = \frac 1 8 B_X^{ - 5/2  } B_Y^{ - 1/2 } \Big (  - 3 \e [ \beta_{2}(Y_1, Y_2) \alpha_1(X_1) \alpha_1^2 (X_2) ] \\
&\quad-  12 \e [ \beta_{2}(Y_1, Y_2)  \alpha_1 (X_1) \alpha_2^2 ( X_1 , X_2 ) ] \Big) .
\end{align*}

Denote by $ \mathcal{D} (i) = \{ (j, k, l) : \max (| j - i |, | k - i  | , | l - i | ) \leq 3 m + 1  \} $. Since $ \{X_{1, i} , 1 \leq i \leq p \} $ are $m$-dependent, it holds that 
\begin{align*}
& \e [ \beta_{2}(Y_1, Y_2) \alpha_1(X_1) \alpha_1^2 (X_2) ]  \\
& =  \sum_{i = 1}^p \sum_{ (j, k, l) \in \mathcal{D}(i) }   \e \big[ \widetilde{Y}_{1, i} ( X_{1, j}^2 - \e X_{1, j}^2 ) \big] \e \big[ \widetilde{Y}_{1, i} ( X_{1, k}^2 - \e X_{1, k}^2 ) (X_{1, l}^2 - \e X_{1, l}^2 )  \big]   \\
&= O (   c_1^8 m^3 p)
\end{align*}
and 
\begin{align*}
& \e [ \beta_{2}(Y_1, Y_2)  \alpha_1 (X_1) \alpha_2^2 ( X_1 , X_2 ) ] \\
&   =  \sum_{i = 1}^p \sum_{(j, k, l) \in \mathcal{D} (i) }  \e \big[ \widetilde{Y}_{1, i} X_{1, k} X_{1, l} \big]  \e \big[ \widetilde{Y}_{1, i} ( X_{1, j}^2 - \e X_{1, j}^2 ) X_{1, k} X_{1, l} \big]   \\
&= O (   c_1 ^8 m^3 p ).
\end{align*}
Consequently, it follows that 
\begin{align*}
| I_3| \lesssim (c_1/c_2)^8  m^{3} p^{ - 2 },
\end{align*}
where $\lesssim$ represents the asymptotic order. By the same token, the symmetry assumptions lead to 
\begin{align*}
I_4 & = \frac{1}{16} B_X^{ - 3/2} B_Y^{ - 3/2} \Big( 4 \e [ \alpha_2^2 (X_1, X_2)  \beta_2^2 (Y_1, Y_2)]  +  \big(  \e [ \alpha_1 (X) \beta_1(Y) ] \big)^2  \\
&  \hspace{1cm} +  8 \e [ \alpha_2^2 (X_1, X_2)  \beta_2^2 ( Y_1, Y_3 ) ]  +  4 \e [ \alpha_2^2 ( X_1, X_2 ) ] \e [ \beta_2^2 (Y_1, Y_2) ] \\
& \hspace{1cm} +  2 \e [ \alpha_2^2 (X_1, X_2) \beta_1(Y_1) \beta_1 (Y_2) ]  - 8 \e [ \alpha_2^2 (X_1, X_2) \beta_1 (Y_1) \beta_2 (Y_1, Y_2) ] \\
& \hspace{1cm} + 2 \e [ \alpha_1 (X_1) \alpha_1 (X_2) \beta_2^2 (Y_1, Y_2) ] - 4 \e [ \alpha_1 (X_1) \alpha_1(X_2) \beta_1 (Y_1) \beta_2 (Y_1, Y_2) ] \\
& \hspace{1cm} - 8 \e [ \alpha_2^2 (X_1, X_2) \beta_1 (Y_3) \beta_2 (Y_1, Y_3)]   \Big) .   
\end{align*}

It is easy to see that $ \big(  \e [ \alpha_1 (X) \beta_1(Y) ] \big)^2    \geq 0 $, $ \e [ \alpha_2^2 (X_1, X_2)  \beta_2^2 ( Y_1, Y_3 ) ] \geq 0  $, 
\begin{align*}
\e [ \alpha_2^2 (X_1, X_2)  \beta_2^2 (Y_1, Y_2)] & = \sum_{i, j, k, l = 1}^p \big(  \e [ X_{1, i} X_{1, j} \widetilde{Y}_{1, k} \widetilde{Y}_{1, l} \big)^2  \geq  \sum_{ |i - k|> m}  \big( \e [ X_{1, i}^2 ] \e [ \widetilde{Y}_{1, k}^2 ] \big)^2 \\
&\geq c_2^8 p(p -2 m ),       
\end{align*} 
and 
\begin{align*}
\e [ \alpha_2^2 ( X_1, X_2 ) ] \e [ \beta_2^2 (Y_1, Y_2) ]  & = \sum_{i , j, k, l= 1}^p \big( \e [X_{1, i} X_{1, j}] \big)^2 \big( \e [ \widetilde{Y}_{1, k} \widetilde{Y}_{1, l} ] \big)^2 \\
& \geq \sum_{i, k} \big( \e  [ X_{1, i}^2 ] \big)^2 \big( \e [ \widetilde{Y}_{1, k}^2 ] \big)^2 \geq c_2^8 p^2  .
\end{align*}
Moreover, since $ \{X_{1, i} , 1\leq i \leq p \} $ are $m$-dependent random variables, we can deduce 
\begin{align*}
\e [ \alpha_2^2 (X_1, X_2) \beta_1(Y_1) \beta_1 (Y_2) ]    & = \sum_{i = 1}^p \sum_{(j, k) \in \mathcal{\widetilde{D}}(i)} \e \big[ X_{1, i} X_{1, j} ( Y_{1, k}^2 - \e Y_{1, k}^2  )  \big] \\
&\quad \times  \e \big[ X_{1, i} X_{1, j} ( Y_{1, k}^2 - \e Y_{1, k}^2  )  \big] \\
& = O ( c_1^8 m^3 p ),  
\end{align*}
where $\mathcal{\widetilde{D}}(i)$ is defined similarly as for $ \mathcal{D} (i)$. In the same fashion, we can show that 
\begin{align*} 
\e [ \alpha_2^2 (X_1, X_2) \beta_1 (Y_1) \beta_2 (Y_1, Y_2) ] & = O (c_1^8 m ^3 p ), \\
\e [ \alpha_1 (X_1) \alpha_1 (X_2) \beta_2^2 (Y_1, Y_2) ]  & = O (c_1^8 m^3 p ), \\
\e [ \alpha_1 (X_1) \alpha_1(X_2) \beta_1 (Y_1) \beta_2 (Y_1, Y_2) ] & = O (c_1^8 m^3 p ), \\
\e [ \alpha_2^2 (X_1, X_2) \beta_1 (Y_3) \beta_2 (Y_1, Y_3)]  & =  O (c_1^8 m ^3 p   ).
\end{align*} 
As a result, there exists some positive constant $A$ depending on $c_1, c_2$, and $ m $ such that
\begin{align*}
I_4 \geq A    p^{ - 1} + O (  p^{ - 2} ) .
\end{align*}

Finally, we deal with term $I_5 $. In view of Lemma \ref{le-dcov}, the first term for the order of $I_5$ is 
\begin{align*}
& B_X^{ 1/2 } B_Y^{1/2 } ( \e | W_{12} |^{5} )^{ 2/5 } ( \e | V_{12} |^5 )^{3/5} \\
& = B_X^{ - 3/2} B_Y^{ - 5/2} \big( \e [ | \Vert X_1 - X_2  \Vert^2 - B_X |^5 ] \big)^{2/5} \big( \e [ | \Vert Y_1 - Y_2 \Vert^2 - B_Y |^5 ] \big)^{3/5}.
\end{align*}
Since $\{ X_{1, i}, 1 \leq i \leq p\}$ are $m$-dependent, without loss of generality we assume that $ s = p/(m +1) $ is an integer. For each $ 1 \leq u \leq m + 1 $, define
\begin{align*}
\mathcal{E}_{u} = \{  (m + 1 ) (j -1) + u: 1 \leq j \leq s \} .
\end{align*}
Clearly, $ \{ X_{1, i}: i \in \mathcal{E}_u \} $ are independent random variables for each $ 1 \leq u \leq m + 1$. Then it follows from the basic inequality $ | \sum_{i = 1}^n a_i |^r  \leq n^{r - 1} \sum_{i = 1}^n | a_i |^r $ for $ r \geq 1 $ and Rosenthal's inequality for independent random variables that 
\begin{align*}
\e & [ | \Vert X_1 - X_2 \Vert^2 - B_X | ^5   ] =  \e \Big[ \Big| \sum_{i = 1}^p \big[ ( X_{1, i} - X_{2, i} )^2 - \e (X_{1, i} - X_{2, i})^2 \big] \Big |^5  \Big] \\
& = \e \Big[ \Big| \sum_{u = 1}^{m + 1 } \sum_{ i \in \mathcal{E}_u }   \big[ ( X_{1, i} - X_{2, i} )^2 - \e (X_{1, i} - X_{2, i})^2 \big]  \Big|^5   \Big] \\
& \leq (m + 1)^4 \sum_{u = 1}^{m + 1}  \e \Big[ \Big|  \sum_{ i \in \mathcal{E}_u }    \big[ ( X_{1, i} - X_{2, i} )^2 - \e (X_{1, i} - X_{2, i})^2 \big]  \Big|^5  \Big] \\
& \leq C (m + 1)^4 \sum_{u = 1}^{m + 1} \bigg\{  \bigg[   \sum_{ i \in \mathcal{E}_u }  \e   \Big( ( X_{1, i} - X_{2, i} )^2 - \e (X_{1, i} - X_{2, i})^2 \Big)^2     \bigg]^{ 5/2 } \\
& \quad + \sum_{ i \in \mathcal{E}_u }  \e   \Big| ( X_{1, i} - X_{2, i} )^2 - \e (X_{1, i} - X_{2, i})^2 \Big| ^5  \bigg\}.
\end{align*}

Note that by assumptions, there exists some absolute positive constant $A$ such that $ \e \big(  \big| ( X_{1, i} - X_{2, i} )^2 - \e (X_{1, i} - X_{2, i})^2 \big| ^5 \big) \leq A  c_1^{10} $ and $  \e \big(  \big[ ( X_{1, i} - X_{2, i} )^2 - \e (X_{1, i} - X_{2, i})^2 \big]^2 \big)  \leq A c_1^4  $, and we have $ B_X \geq 2 c_2^2 p $. Then it follows that 
\begin{align*}
\e [ | \Vert X_1 - X_2 \Vert^2 - B_X | ^5   ] \lesssim  c_1^{10} m^4 \cdot m \cdot (p / m)^{5/2}  = c_1^{10} m^{5/2} p^{5/2} .
\end{align*}
Similarly, we can obtain
\begin{align*}
\e [  | \Vert Y_1 - Y_2 \Vert^2 - B_Y |^5   ] \lesssim c_1^{10} m^{5/2} p^{5/2}.
\end{align*}
Hence it holds that 
\begin{align*}
B_X^{1/2} B_Y^{1/2}  ( \e | W_{12} |^{5} )^{ 2/5 } ( \e | V_{12} |^5 )^{3/5} \lesssim m^{5/2} p^{- 3/2}.
\end{align*}
In the same manner, we can deduce that 
\begin{align*}
B_X^{ 1/2 } B_Y^{1/2} ( \e | W_{12} |^5  )^{3/5}  ( \e | V_{12} |^5  )^{ 2/5 } & \lesssim m^{5/2} p^{ - 3/2 },\\
B_X^{1/2} B_Y^{1/2} ( \e | W_{12} |^{5} )^{1/5} ( \e | V_{12} |^5 )^{4/5} & \lesssim m^{5/2} p^{ - 3/2  } , \\
B_X^{1/2} B_Y^{1/2} ( \e | W_{12} |^{5} )^{4/5} ( \e | V_{12} |^5 )^{1/5} & \lesssim m^{5/2} p^{ - 3/2  } , \\
B_X^{1/2} B_Y^{1/2} ( \e | W_{12} |^6 )^{1/2} ( \e | V_{12} |^6 )^{1/2} & \lesssim m^3 p^{ - 2}. 
\end{align*}
Thus substituting the above five inequalities into the order of $I_5$ in Lemma \ref{le-dcov} yields that there exists some positive constant $A$ depending on $c_1, c_2$, and $m$ such that 
\begin{align*}
I_5 \leq A p^{ - 3/2}.
\end{align*}

{
	As a consequence, combining all the bounds above leads to 
	\begin{align}
	\V^2 (X, Y) \geq A p^{ - 1} + O (p^{ - 3 / 2} ). \label{v-bound}
	\end{align}
	Hence this entails that when $ p = o (\sqrt n) $, it holds that $ \sqrt n \V^2 (X, Y) \to \infty$. Furthermore, it follows from Proposition \ref{prop2} that $ \V^2 (X) = B_X^{ - 1} \e [ ( X_1^T X_2 )^2 ] + O ( B_X^{-2} L_{x, 1/2} ) $. By the assumptions $ \e (X_{1, i}^{12} ) + \e ( Y_{1, i}^{12} ) \leq c_1^{12} $, $ \var (X_{1, i} ) \geq c_2^2 $, and $ \var( Y_{1, i}) \geq c_2^2 $, it is easy to see that $ 2 c_2 ^2 p \leq B_X \leq 2 c_1^2 p $ and $ 2 c_2^{2} p \leq B_Y \leq 2 c_1^2 p $. Since $\{ X_{1, i}, 1 \leq i \leq p \}$ are $m$-dependent, we have 
	\begin{align*}
	\sum_{i = 1}^p ( \e [ X_{1, i}^2 ] )^2  & \leq   \e [ (X_1^T X_2)^2 ] 
	= \sum_{i = 1}^p \sum_{j = 1}^p (\e [ X_{1, i} X_{1, j} ] )^2  \\
	&=  \sum_{i = 1}^p \sum_{| i - j  | \leq m } ( \e [ X_{1, i} X_{1, j} ] )^2 ,
	\end{align*}
	which yields $ c_2^4 p \leq \e [ ( X_1^T X_2 )^2 ]  \leq 2(m + 1) c_1^4 p $. In the same manner, we can obtain $  L_{x, 1/2} \leq C p,  L_{y, 1/2} \leq C p, L_{x} \leq C p^2$, and $L_{y} \leq C p^2 $ with some positive constant $C$ depending on $c_1, c_2 $, and $m$. Consequently, there exist some positive constants $C_1$ and $C_2$ depending on $ c_1, c_2 $, and $m$ such that $ C_1 \leq \V^2 (X) \leq C_2 $. Similarly, we have $ C_1 \leq \V^2 (Y)  \leq C_2$. This along with \eqref{v-bound} entails that $ \Rc^2 (X, Y) \geq A p^{-1} + O(p^{ - 3/2}) $, where $A > 0$ is some constant depending on $c_1, c_2 $, and $m$. 
	
	From the above analysis, it holds that $B_X^{1-/2} B_Y^{-1/2} L_x^{1/4} L_y^{1/4} \leq A_1 $ for some positive constant $A_1$ depending on $ c_1, c_2 $, and $m$. Thus we can obtain under the assumption of $ p = o(\sqrt n) $ that
	$$ 
	n \Rc^2 (X, Y) \geq A n p^{-1} + O(n p^{ - 3/2}) \to \infty
	$$
	and 
	$$
	\sqrt{n} \V^2 (X, Y) / ( B_X^{- 1/2} B_Y^{ - 1/2} L_x^{1/4} L_y^{1/4})  \geq A \sqrt n p^{-1} + O (\sqrt n p^{- 3/2}) \to \infty.
	$$
	Finally, it follows from Theorem \ref{thm-power} that for any arbitrarily large $C > 0$, $\mathbb{P} (T_n > C ) \to 1 $  as $n \to \infty$, which concludes the proof of Theorem \ref{prop-power}.
}

\subsection{Proof of Proposition \ref{prop1}} \label{pf-prop1}
In view of the definition $  B_X = \e [ \Vert X_1 - X_2 \Vert^2  ] $, we can write
\begin{align}
d(X_1, X_2 ) &  =  ( \Vert X_1 - X_2 \Vert -  B_X^{1/2}  )  - \e [  (  \Vert X_1 - X_2 \Vert -  B_X^{1/2}  ) | X_1 ]  \nn \\
&   \quad -  \e [  (  \Vert X_1 - X_2 \Vert -  B_X^{1/2}  ) | X_2 ]  + \e ( \Vert X_1 - X_2 \Vert   - B_X^{1/2} ).   \label{d}
\end{align}
Thus it follows from Jensen's inequality that for $ \tau > 0  $,
\begin{align*}
\e [ | d ( X_1, X_2 ) |^{2 + 2 \tau } ] & \leq C_{ \tau }  \e \big[ \big| \Vert X_1 - X_2 \Vert -  B_X^{1/2} \big|^{ 2 + 2 \tau } \big] \nn \\
& = C_{ \tau }   \e \Big [   \frac {      \big|  \Vert X_1 - X_2 \Vert ^2  - B_X \big|^{ 2 + 2 \tau  }    } {  ( \Vert X_1 - X_2 \Vert +  B_X^{1/2} )^{ 2 + 2 \tau  } } \Big]  \nn \\
& \leq   C_{ \tau }  B_X^{ - ( 1 + \tau ) }  \e \big[ \big|  \Vert X_1 - X_2 \Vert^2  - B_X  \big|^{ 2 + 2 \tau } \big] .
\end{align*}
Moreover, we have 
\begin{align}
& \e  \big[ \big|  \Vert X_1 - X_2 \Vert^2  - B_X  \big|^{ 2 + 2 \tau } \big]  \nn  \\
& \leq  C_{ \tau } \Big\{  \e \big[  \big|  \Vert X_1 \Vert^2 - \e [ \Vert X_1 \Vert^2 ]  \big|^{ 2 + 2 \tau } \big] +   \e \big[  \big|  \Vert X_2 \Vert^2 - \e [ \Vert X_2 \Vert^2 ]  \big|^{ 2 + 2 \tau } \big] +  \e [ | X_1^T X_2 |^{ 2 + 2 \tau } ]  \Big\} \nn  \\
& \leq C_{ \tau } L_{ x, \tau },  \label{L_xt}
\end{align}
which completes the proof of Proposition \ref{prop1}.

\subsection{Proof of Proposition \ref{prop2}} \label{pf-prop2}

The essential idea of the proof is to conduct the Taylor expansion for function $ (1 + x)^{ 1/ 2 } $ to relate the $L_1$-norm to the $L_2$-norm. Let us define 
$$
b( X_1, X_2 ) = \Vert X_1 - X_2 \Vert - B_X^{1/2}, ~ b_1 ( X_1 ) = \e [   b ( X_1, X_2 ) | X_1 ], ~ b_1 ( X_2 ) = \e [   b ( X_1, X_2 ) | X_2 ].
$$
Since $ \V^2 ( X )  = \e [ d^2 (X_1, X_2 ) ] $, it follows from \eqref{d} that 
\begin{align*}
\V^2 ( X )  & =  \e \big \{  b ( X_1 , X_2  ) - b_1 ( X_1  )  -  b_1 ( X_2 )  +  \e [ b ( X_1 , X_2 ) ]  \big\} ^2.
\end{align*}
Then by expanding the square and the symmetry of $ X_1 $ and $ X_2  $, we can obtain 
\begin{align*}
\V^2 ( X  ) = \e [ b^2 ( X_1, X_2 ) ] - 2 \e [ b_1^2 ( X_1 ) ] + \{ \e [ b ( X_1, X_2 )  ] \}^2.
\end{align*}

Next we will bound the moments $ \e [ b^2 ( X_1, X_2 )] $, $ \e [ b_1^2 (X_1) ] $, and $ \e [ b( X_1, X_2 ) ] $ by resorting to the basic inequailties in Lemma \ref{le-cru} in Section \ref{SecC.7} of Supplementary Material. Denote by 
$$
W_{1 2} =  B_X^{ -1  } \big(  \Vert X_1 - X_2 \Vert^2 - B_X  \big) \ \text{ and } \  W_{1 3} = B_X^{ -1 } \big(  \Vert X_1 - X_3 \Vert^2 - B_X  \big).
$$
Observe that $ W_{12} \geq -1, W_{13} \geq - 1 $, and $ \e [ W_{12}  ] = \e [ W_{13} ] = 0 $.
For term $ \e [ b( X_1, X_2 ) ] $, by \eqref{cru1} and \eqref{cru2} we have
\begin{align}
\e [ b( X_1, X_2 ) ]  & = B_X^{ 1 / 2 } \Big[  \e  \big(  [ (  1 + W_{12} )^{1/2}  - 1 ] \I \{ W_{12} \leq 1  \}  \big) +   \e  \big(  [ (  1 + W_{12} )^{1/2}  - 1 ] \I \{ W_{12} > 1  \}  \big) \Big] \nn \\
& = B_X^{1/2} \Big[  \frac { 1 }  { 2 }  \e W_{12}  \I \{ W_{12} \leq 1  \}  + O_1 \e W_{12}^2  \I \{ W_{12} \leq 1  \}  + O_2 \e | W_{12} |  \I \{ W_{12} > 1  \}    \Big] \nn \\
& =    O_1 B_X^{1/2} \e W_{12}^2  \I \{ W_{12} \leq 1  \}  +  O_3 B_X^{1/ 2 } \e | W_{12} |  \I \{ W_{12} > 1  \}   \big) ,   \label{de2}
\end{align}
where $ \I \{ \cdot \} $ denotes the indicator function and $ O_1, O_2, O_3$ are bounded quantities such that $ | O_1 | \leq 1 /2$, $| O_2 | \leq 1 $, and $| O_3 | \leq 3/2 $. Thus it follows that 
\begin{align}
\{ \e [ b( X_1, X_2 ) ]  \}^2 & \leq  B_X  \big( \frac 1 2 \e W_{12}^2  \I \{ W_{12} \leq 1  \}  + \frac 3 2  \e | W_{12} |  \I \{ W_{12} > 1  \}   \big) ^2 \nn \\
& \leq      B_X  \Big(  \frac {1} { 4 } \e | W_{12} |^3  \I \{ W_{12} \leq 1  \} +  \frac {15}  {4} \e W_{12}^2  \I \{ W_{12} >  1  \}    \Big).  \label{delta-square}
\end{align}
If $ \e [ \Vert X \Vert^{4 + 4 \tau }  < \infty $ for some $ 0 < \tau \leq 1/2 $, then it holds that 
\begin{align}
\{ \e [ b( X_1, X_2 ) ]  \}^2  \leq \frac { 15 } { 4 } B_X \e [ | W_{12} |^{ 2 + 2 \tau }  ]. \label{part3}
\end{align} 

Similarly, by \eqref{cru1} and \eqref{cru2} again, for $ 0 < \tau \leq 1/ 2 $ we have 
\begin{align*}
\e [ b^2 ( X_1, X_2 ) ] & =  B_X \e \big(   (  1 + W_{12} )^{1/2}  - 1 \big)^2 \nn \\
& = B_X  \Big(  \e \big[ \big(  \frac 12 W_{12} + O_5 W_{12}^2 \big) \I \{  W_{12} \leq 1 \}  \big]^2  +  O_4 \e W_{12}^2 \I \{  W_{12}  > 1 \}  \Big) ,
\end{align*}
where $ | O_4 | \leq 1 $ and $ | O_5 | \leq 1/2 $. Hence for $ 0 < \tau \leq 1 /2  $, it holds that 
\begin{align}
\Big |  \e [ b^2 ( X_1, X_2 ) ] - \frac 14 B_X \e [ W_{12}^2  ] \Big| & \leq B_X \Big( \frac 5 4 \e [ W_{12}^2 \I \{ W_{12}  > 1  \}    ]  +  \frac  34 \e [ | W_{ 12 } |^3 \I\{ W_{12} \leq 1   \} ]  \Big) \nn  \\
& \leq \frac 5  4  B_X \e [  | W_{ 12 }  |^{ 2 + 2 \tau }  ] .      \label{part1}
\end{align}
Again it follows from \eqref{cru1} and \eqref{cru2} that for $ 0 < \tau \leq 1 /2  $, we have 
\begin{align}
\e [ b_1^2 ( X_1 ) ] & = \e [ b ( X_1, X_2 ) b ( X_1, X_3 ) ] \nn \\
& = B_X   \e  \big\{ [  ( 1 + W_{12} )^{1/2} - 1   ] [  ( 1 + W_{13} )^{1/2} - 1  ] \I \{  \max ( W_{12}, W_{13} ) \leq 1 \} \big\} \nn \\
& \quad + B_X  \e  \big\{ [  ( 1 + W_{12} )^{1/2} - 1   ] [  ( 1 + W_{13} )^{1/2} - 1  ]  \I \{  \max ( W_{12}, W_{13} ) >  1 \} \big\} \nn \\
& = B_X \Big( \frac 1 4 \e [ W_{12} W_{13} ] +   O_7   \e [ W_{12}^2 W_{13} \I \{  \max ( W_{12}, W_{13} ) \leq 1 \} ]  \nn \\
& \hspace{3cm}   + O_8  \e [ W_{12} W_{13} \I \{  \max ( W_{12}, W_{13} ) >  1 \} ]  \Big) \nn \\
& = \frac 1 4 B_X \e [ W_{12} W_{13} ] + O_9 B_X \e [ | W_{12}  |^{2 + 2 \tau} ],  \label{part2}
\end{align}
where $ O_7$, $O_8$, and $O_9 $ are bounded quantities satisfying $ | O_7 | \leq 3/4 $, $ | O_8 | \leq 5/4  $, and $  | O_9 | \leq   4 $.

Finally by combining \eqref{part3}--\eqref{part2} we can deduce 
\begin{equation}
\Big|  \V^2 ( X , X )  -   \frac {  B_X } { 4 } ( \e [ W_{12}^2 ] - 2 \e [ W_{12} W_{13} ] )  \Big| \leq  9  B_X  \e [ | W_{12}  |^{2 + 2 \tau} ].   \label{V2}
\end{equation}
Moreover, Lemma \ref{le-v2} in Section \ref{SecC.8} of Supplementary Material yields
\begin{align*}
\frac   { B_X } { 4  } \e [ W_{ 12 }^2 - 2 \e [ W_{ 12 } W_{13} ] ]   =  B_X^{ -1  } \e [ (X_1^T X_2)^2 ].
\end{align*}
It follows from \eqref{L_xt} that 
\begin{align*}
B_X \e [ | W_{ 12 } |^{ 2 + 2 \tau  } ] \leq B_X^{ - (1 + 2 \tau  ) } L_{ x, \tau }.
\end{align*}
Thus the desired result \eqref{V} can be derived. This concludes the proof of Proposition \ref{prop2}.

\subsection{Proof of Proposition \ref{prop3}}  \label{pf-prop3}
Similar to the proof of Proposition \ref{prop2}, the main idea of the proof is to conduct the Taylor expansion to relate the $ L_1$-norm to the $ L_2 $-norm. Denote by $ \Delta =   \e [ b ( X_1 , X_2 )] = \e [ \Vert X_1 - X_2 \Vert - B_X^{1/2}  ]$. In light of \eqref{d}, we have
\begin{align*}
& \e [ g( X_1, X_2, X_3, X_4 ) ]  \\
&  = \e \big[   \big(  b ( X_1 , X_2  ) - b_1 ( X_1  )  -  b_1 ( X_2 )  + \Delta ]   \big)  \big( b ( X_1 , X_3  ) - b_1 ( X_1  )  -  b_1 ( X_3 )  + \Delta  \big)  \\
& \qquad  \times  \big( b ( X_2 , X_4  ) - b_1 ( X_2  )  -  b_1 ( X_4 )  +  \Delta \big)  \big(  b ( X_3 , X_4  ) - b_1 ( X_3  )  -  b_1 ( X_4 ) + \Delta \big) \big].
\end{align*}

\noi Expanding the products and noting that $ X_1, X_2, X_3, X_4 $ are i.i.d. random variables, we can deduce
\begin{align}
\e [ g( X_1, X_2, X_3, X_4 ) ] =  G_1   - 4 G_2 + 2 G_3^2 + 4 \Delta G_4 - 4 \Delta^2 G_3 + \Delta^4,    \label{g-expan}
\end{align}
where 
\begin{align*}
G_1  & =   \e [ b ( X_1, X_2 ) b ( X_1, X_3 ) b ( X_2, X_4 ) b ( X_3, X_4 )  ], \\
G_2 & = \e [ b ( X_1, X_2 )  b_1 ( X_1, X_3 ) b( X_2, X_4 ) b ( X_4, X_5 ) ], \\
G_3 & =   \e [  b ( X_1, X_2 ) b (X_1, X_3 ) ],    \\
G_4 & =    \e [ b ( X_1, X_2 ) b ( X_1, X_3  ) b( X_2, X_4 ) ].
\end{align*}

Next we will analyze the six terms on the right hand side of (\ref{g-expan}) separately.  The same technique as in the proof of Proposition \ref{prop2} will be used. For any $ i \neq j $, let us define 
\begin{align*}
W_{ij}  =  B_X^{ -1 }  \big(\Vert X_i - X_j \Vert^2 - B_X \big).  
\end{align*}
First for term $ G_1 $, by definition it holds that 
\begin{align*}
G_1   =  B_X^2    \e  \big[  \{  ( 1 + W_{12} )^{1/2} - 1 \}    \{ ( 1 + W_{13} )^{1/2} - 1  \}  \{ ( 1 + W_{24} )^{1/2} - 1  \}  \{ ( 1 + W_{34} )^{1/2} - 1  \}  \big].
\end{align*}
Denote by 
\begin{align*}
D_1 = \{  \max ( W_{12} , W_{ 13} , W_{24}, W_{34} ) \leq 1 \} 
\end{align*} 
and $ D_1^{ c } $ the complement of $ D_1 $.  
By separating the integration region into $ D_1 $ and $ D_1^c $ and applying \eqref{cru1} and \eqref{cru2}, we can deduce
\begin{align*}
G_1  & =  B_X^2 \e \Big(  \big[ \frac 12 W_{12} + O(1) ( W_{12}^2 )  \big]  \big[ \frac 12 W_{13} + O(1) ( W_{13}^2 )  \big] \nn \\
& \hspace{4cm} \times  \big[ \frac 12 W_{24} + O(1) ( W_{24}^2 )  \big]  \big[ \frac 12 W_{34} + O(1) ( W_{34}^2 )  \big] \I \{ D_1 \} \Big) \nn \\
& \hspace{1cm} +  O(1) B_X^2  \e [ | W_{12} W_{13} W_{24} W_{34} | \I \{ D_1^{c} \} ], 
\end{align*}
where $ O(1) $ represents a bounded quantity satisfying $ | O(1) | \leq C $ for some absolute positive constant $C$. It follows from expanding the products and Chebyshev's inequality that if $ \e [  \Vert X \Vert^{ 4 + 4 \tau } ] < \infty $ for some $ 0 < \tau \leq 1/2 $, then we have 
\begin{align*}
\Big | G_1 -  \frac { B_X^2  } { 16 } \e [ W_{12} W_{13}  W_{24} W_{34} ]  \Big|  \leq C  B_X^2  \e [ | W_{12} |^{1 + 2 \tau }  | W_{13} | | W_{24}  | | W_{34} | ].
\end{align*}

Further, by conditioning on $ X_2, X_3 $, applying the Cauchy--Schwarz inequality, and noting that $ X_1, X_2, X_3, X_4$ are i.i.d. random variables, it holds that 
\begin{align}
\e [  | W_{12} |^{1 + 2 \tau }  | W_{13} | | W_{24}  | | W_{34} |  ] & =  \e \big\{ \e \big( | W_{12} | ^{ 1 + 2 \tau } | W_{13} | \big\vert X_2, X_3 \big) \e \big( | W_{24}   W_{34} | \big\vert X_2, X_3 \big)   \big\}  \nn \\
& \leq   \e \big\{ \big ( \e [ | W_{12} | ^{ 2 + 2 \tau } \vert X_2 ]  \big)^{ \frac { 1 + 2 \tau } { 2 + 2 \tau  } } \big( \e [ | W_{13} |^{ 2 + 2 \tau } \vert X_3 ] \big)^{ \frac {1}   { 2 + 2 \tau  } } \nn \\
& \qquad \times  \big( \e [ | W_{24} |^{ 2 + 2 \tau } \vert X_2 ] \big)^{ \frac {1}   { 2 + 2 \tau  } }    \big( \e [ | W_{3 4} |^{ 2 + 2 \tau } \vert X_3 ] \big)^{ \frac {1}   { 2 + 2 \tau  } }  \big\} \nn  \\
&  =   \e \big\{  \e [ | W_{12} | ^{ 2 + 2 \tau } \vert X_2 ]  \big\} \times  \e \big\{ \big( \e [ | W_{13} |^{ 2 + 2 \tau } \vert X_3 ] \big)^{ \frac { 1 }   { 1 + \tau  } }  \big\} \nn \\
&  \leq \big(  \e [ | W_{12} |^{ 2 + 2 \tau } ] \big)^{ \frac {  2 + \tau } { 1 +  \tau  } }.   \label{er-bdd}
\end{align}
Consequently, we have 
\begin{align}
\Big | G_1 -  \frac { B_X^2  } { 16 } \e [ W_{12} W_{13}  W_{24} W_{34} ]  \Big|   \leq  C B_X^2  \big(  \e [ | W_{12} |^{ 2 + 2 \tau } ] \big)^{ \frac {  2 + \tau } { 1 +  \tau  } }.  \label{g-r-p1}
\end{align}
An application of the similar argument as for the proof of \eqref{g-r-p1} yields 
\begin{align}
G_2 =  \frac { B_X^2  } { 16 } \e [ W_{12} W_{13} W_{24} W_{45} ] + O(1) B_X^2 \big(  \e [ | W_{12} |^{ 2 + 2 \tau } ] \big)^{ \frac {  2 + \tau } { 1 +  \tau  } }.  \label{g-r-p2}
\end{align}

As for term $ G_3^2  $, by the same token we can deduce 
\begin{align*}
G_3  & =  B_X \e  \big\{  [ ( 1 + W_{12} )^{1/2} - 1 ] [ ( 1 + W_{13} )^{1/2} - 1 ]  \I \{  \max ( W_{12}, W_{13} ) \leq 1 \} \big\} \nn \\
&  \quad + B_X \e  \big\{  [ ( 1 + W_{12} )^{1/2} - 1 ] [ ( 1 + W_{13} )^{1/2} - 1 ]  \I \{  \max ( W_{12}, W_{13} ) >  1 \} \big\} \nn \\
& = \frac { B_X } { 4 } \e [ W_{12} W_{13} ] + O(1) B_X \delta_1,   
\end{align*}
where $ \delta_1  =  \e [  W_{12}^2 | W_{13} | \I \{  \max ( W_{12}, W_{13}  ) \leq 1 \} ] + \e [ | W_{12} W_{13}  | \I \{  \max ( W_{12}, W_{13}  ) > 1   \} ]  $. Observe that when $ 0 < \tau \leq 1 / 2 $, we have 
\begin{align*}
\delta_1  \cdot  | \e [ W_{12} W_{13} ] |  & \leq  2 \e [ | W_{12} |^{ 1 + 2 \tau } | W_{13} |  ] \e [ | W_{12} W_{13} | ]   \leq 2 \big(  \e [ | W_{12} |^{ 2 + 2 \tau } ] \big)^{  \frac { 2 + \tau  } { 1 + \tau  } }
\end{align*}
and 
\begin{align*}
\delta_1^2 &    \leq   4 ( \e [ | W_{12} |^{1 + \tau } | W_{13} | ] )^2  \leq  4 \big( \e [ | W_{12} |^{ 2 + 2 \tau  } ] \big)^{ \frac { 2 + \tau  } { 1 + \tau  } }.
\end{align*}
As a consequence, it holds that 
\begin{align}
\Big| G_3^2  -  \frac { B_X^2  } { 16 }  \big(  \e [ W_{12} W_{13} ] \big)^2 \Big| \leq  C B_X^2  \big( \e [ | W_{12} |^{ 2 + 2 \tau  } ] \big)^{ \frac { 2 + \tau  } { 1 + \tau  } }.      \label{g-r-p3}
\end{align}

We next deal with term $\Delta  G_4 $. It follows from \eqref{cru1} and the Cauchy--Schwarz inequality that
\begin{align}
|G_4 | &  =   B_X^{3/2}  \big| \e  \big\{  [ ( 1 + W_{12} )^{ 1 / 2 }  - 1  ] [ ( 1 + W_{13} )^{ 1 / 2 }  - 1 ]  [ ( 1 + W_{24} )^{ 1 / 2 }  - 1 ] \big\} \big| \nn \\
& \leq   B_X^{3/2}   \e [ | W_{12} W_{13} W_{24} | ] =   B_X^{3/2}   \e \{  \e ( | W_{12} W_{13} | \vert X_2, X_3 )  \e ( | W_{24} | \vert X_2 ) \}  \nn \\
& \leq B_X^{3/2}   \e \{ ( \e [ W_{12}^2 \vert X_2 ] )^{ 1 / 2 }   ( \e [ W_{13}^2 \vert X_3 ] )^{ 1 / 2 }  ( \e [ W_{24}^2 \vert X_2 ]   )^{ 1 / 2 } \} \nn \\
&  = B_X^{3/2}  \e \{  \e [ W_{12}^2 \vert X_2 ]  ( \e [ W_{13}^2 \vert X_3 ] )^{ 1 / 2 } \}  \nn \\
&  \leq B_X^{3/2}   (  \e [ W_{12}^2  ] )^{3/2}.  \label{g-p4-1}
\end{align}
Moreover, \eqref{de2} entails that for $ 0 < \tau \leq 1/2   $, we have 
\begin{align*}
|  \Delta |  = |  \e [ b ( X_1, X_2 ) ]  |  \leq C  B_X^{1/2} \e [ | W_{12} |^{ 1 + 2 \tau  } ].
\end{align*}
As a result, it follows that 
\begin{align}
|  \Delta G_4  | &  \leq  C B_X^2  (  \e [ W_{12}^2  ] )^{3/2}   \e [ | W_{12} |^{ 1 + 2 \tau  } ] \nn \\
& \leq  C B_X^2  \big( \e [ | W_{12} |^{ 2 + 2 \tau  } ] \big)^{ \frac { 2 + \tau  } { 1 + \tau  } }.    \label{g-r-p4}
\end{align}

As for term $ \Delta^2  G_3 $, note that \eqref{cru1} leads to 
\begin{align*}
|  G_3 | & = B_X \big| \e  \big\{  [ ( 1 + W_{12} )^{ 1 / 2 }  - 1  ] [ ( 1 + W_{13} )^{ 1 / 2 }  - 1 ]  \big\} \big| \nn \\
& \leq B_X \e [ | W_{12} W_{13} | ] \leq B_X \e [ W_{12}^2 ].
\end{align*}
It follows from \eqref{part3} that for $ 0 < \tau \leq 1/ 2 $, we have 
\begin{align*}
\Delta^2 \leq C  B_X \e [ | W_{12} |^{ 2 + 2 \tau } ].
\end{align*}
Hence it holds that 
\begin{align}
\Delta^2 |  G_3 |  \leq C  B_X^2 \e [  | W_{12} |^{ 2 + 2 \tau } ] \e [ W_{12}^2 ] \leq C B_X^2 ( \e [ | W_{12} |^{ 2 + 2 \tau  } ] )^{ \frac { 2 + \tau  } { 1 + \tau  } } .  \label{g-r-p5}
\end{align}

Furthermore, note that \eqref{delta-square} implies that for $ 0 < \tau \leq 1/2  $, we have 
\begin{align}
\Delta^4 \leq  C B_X^2 (  \e [ | W_{12} |^{ 2 + \tau } ] )^2  \leq  C B_X^2 ( \e [ | W_{12} |^{ 2 + 2 \tau } ] )^{ \frac { 2 + \tau  } { 1 + \tau  } }.     \label{g-r-p6}
\end{align}
Therefore, by substituting \eqref{g-r-p1}--\eqref{g-r-p3} and \eqref{g-r-p4}--\eqref{g-r-p6} into \eqref{g-expan} we can obtain that if $ \e \Vert X \Vert ^{4 + 4 \tau } < \infty $ for  some $ 0 < \tau \leq 1 /2  $, then 
\begin{align*}
\e [ g ( X_1, X_2, X_3, X_4 ) ]   & = \frac { B_X^2  } { 16 } \Big\{ \e [ W_{12} W_{13} W_{24} W_{34} ] - 4 \e [ W_{12} W_{13} W_{24} W_{45}  ]  \nn \\
& \hspace{2cm} + 2 ( \e [ W_{12} W_{13} ] )^2 + O(1)   \big(  \e [ | W_{12} |^{ 2 + 2 \tau } ] \big)^{ \frac {  2 + \tau } { 1 +  \tau  } }  \Big\}.
\end{align*}
Finally, the desired result \eqref{g1} can be derived from \eqref{L_xt} and Lemma \ref{le-g-mo1} given in Section \ref{SecC.9} of Supplementary Material. This completes the proof of Proposition \ref{prop3}.

\renewcommand{\thesubsection}{B.\arabic{subsection}}
\section{Proofs of Propositions \ref{cor0}--\ref{cor2}} \label{proof}

\subsection{Proof of Proposition \ref{cor0}} \label{SecB.1}
The desired result follows from Theorem \ref{thm4}. By conditions \eqref{cor0-c1}--\eqref{cor0-c4}, it holds that 
\begin{align*}
B_{X}^{ - 2 \tau  }  L_{ x, \tau } / \e [ ( X_1^T X_2 )^2 ] \leq c_1 c_2^{ - ( 1 + 2 \tau ) } p^{ - \tau } \ \text{ and } \ 
B_{Y}^{ - 2 \tau  }  L_{ y, \tau } / \e [ ( Y_1^T Y_2 )^2 ] \leq c_1 c_2^{ - ( 1 + 2 \tau ) } q^{ - \tau }.
\end{align*}
Thus by Theorem \ref{thm4}, the fact that $ p \rightarrow \infty$ and $ q \rightarrow \infty $, and substituting the bounds in  \eqref{cor0-c1}--\eqref{cor0-c4} into \eqref{rate1}, we can obtain
\begin{align*}
\sup \limits_{ x \in \mathbb{R} } | \mathbb{ P  } ( T_n \leq x  )  - \Phi (x) | \leq  A( c_1, c_2 ) \big[  ( p q  )^{ - \tau  ( 1 + \tau ) / 2 }  +  n^{ - \tau }  \big]^{ 1/ ( 3 + 2 \tau ) },
\end{align*}
which concludes the proof of Proposition \ref{cor0}.

\subsection{Proof of Proposition \ref{cor1}}  \label{pf-cor1}
The proof is based on Theorem \ref{thm4} in Section \ref{rate} for the case of $ 0 < \tau \leq 1 /2  $ and Theorem \ref{thm8} in Section \ref{SecD.1} for the case of $ 1 /2 < \tau \leq 1  $. We need to calculate the moments involved therein. The main idea is to use the block technique to deal with the $m$-dependent structure so that the moment inequalities for independent random variables can be applied. For simplicity, assume that $ k = p / ( m_1 + 1 ) $ is an integer. For $ 1 \leq r \leq k $, we define
\begin{align*}
H_r = \{ i :  ( k - 1 ) ( m_1 + 1 ) + 1 \leq i \leq k  ( m_1 + 1 )  \}
\end{align*}
and 
\begin{align*}
S_{1, r} =   \sum \limits_{  i \in H_r }  ( X_{1, i}^2 - \e [ X_{1, i}^2 ] ), \quad S_{2, r} =   \sum \limits_{  i \in H_r } X_{1, i} X_{2, i}.
\end{align*}
By the $ m_1 $-dependent component structure of random vector $ X $, the odd blocks are mutually independent and so are the even blocks. Hence $ \{  S_{1, r} ,  r ~ \text{is odd} \} $, $ \{  S_{1, r}, r ~ \text{is even} \} $, $ \{  S_{2, r} ,  r ~\text{is odd} \} $, and $ \{  S_{2, r}, r ~ \text{is even} \} $ are sequences of independent random variables with zero mean, respectively. 

Let us first analyze term $ L_{ x, \tau } $. It holds that 
\begin{align*}
\e  (  | \Vert X \Vert^2 - \e \Vert X \Vert^{ 2} |^{ 2 + 2 \tau  } ) & =  \e \Big( \Big|  \sum_{ r = 1}^k S_{1, r} \Big|^{ 2 + 2 \tau } \Big) \nn \\
& \leq C  \Big( \e \Big[ \Big|   \sum_{r : \,odd} S_{1, r} \Big|^{ 2 + 2 \tau } \Big] +  \e \Big[ \Big|   \sum_{r : \, even} S_{1, r} \Big|^{ 2 + 2 \tau } \Big]  \Big). 
\end{align*}
Then it follows from Rosenthal's inequality that 
\begin{align*}
& \e  (  | \Vert X \Vert^2 - \e \Vert X \Vert^{ 2} |^{ 2 + 2 \tau } )  \nn\\
& \leq  C   \Big\{ \big(   \sum_{r : \, odd} \e [ S_{1, r}^2 ] \big)^{ 1 + \tau } + \big(   \sum_{r : \, even} \e [ S_{1, r}^2 ] \big)^{ 1 + \tau }  +     \sum_{r = 1}^k \e [ | S_{1, r} |^{ 2 + 2 \tau } ]  \Big\}. 
\end{align*}
Note that for positive numbers $ s >1  $ and $ t > 1 $ with $  s^{-1} + t^{ -1} = 1$, we have 
\begin{align}
\Big| \sum_{i = 1}^ n a_i b_i \Big| \leq \Big(\sum_{i = 1}^n  |a_i|^{ s } \Big)^{  1/s } \Big( \sum_{i = 1}^n b_i^{ t} \Big)^{ 1/t}. \label{ineq-m}
\end{align}
Thus we can deduce 
\begin{align*}
\e [ S_{1, r}^2 ] & = \e \Big[   \sum_{i \in H_r} ( X_{1, i}^2 - \e [ X_{1, i}^2 ] ) \Big]^2 \nn \\
& \leq ( m_1 + 1 )   \sum_{ i \in H_r }   \e  [ ( X_{1, i}^2 - \e [ X_{1, i}^2 ] )^2 ] \leq ( m_1 + 1 )  \sum_{ i \in H_r }  \e  [  X_{1, i}^4 ]     \nn \\  
\intertext{and}
\e [ | S_{1, r} |^{ 2 + 2 \tau } ] & \leq ( m_1 + 1 )^{ 1 + 2 \tau }  \sum_{i \in H_r}  \e [ | X_{1, i}^2 - \e [ X_{1, i}^2 ] |^{ 2 + 2 \tau }  ]          \nn \\
& \leq C   ( m_1 + 1 )^{ 1 + 2 \tau }   \sum_{i \in H_r}  \e [ | X_{1, i} |^{ 4 + 4 \tau } ] .              
\end{align*}
By plugging in the above bounds and applying \eqref{ineq-m}, it follows that 
\begin{align}
\e  (  | \Vert X \Vert^2 - \e \Vert X \Vert^{ 2} |^{ 2 + 2 \tau  } ) & \leq C  \big\{  ( m_1 + 1 )^{ 1 + \tau  } ( p / 2 )^{ \tau } \sum_{r :\, odd}   \sum_{ i \in H_r}  ( \e [  X_{1, i}^4 ] )^{ 1 + \tau  }  \nn \\
& \qquad +  ( m_1 + 1 )^{  1 + \tau } ( p / 2 )^{ \tau } \sum_{r :\, even}   \sum_{ i \in H_r}  ( \e [  X_{1, i}^4 ] )^{  1 + \tau }  \nn \\ 
& \qquad + ( m_1 + 1 )^{ 1 + 2 \tau }   \sum_{ r = 1}^k \sum_{ i \in H_r} \e [ | X_{1, i} |^{ 4 + 4 \tau } ] \big\} \nn \\
& \leq C  ( m_1 + 1 )^{ 1 + \tau } p ^{ \tau }      \sum_{i = 1}^{p} \e [  | X_{1, i} |^{ 4 + 4 \tau } ].  \label{L1} 
\end{align}

In a similar fashion, we have
\begin{align*}
\e [ | X_1^T X_2 |^{ 2 + 2 \tau } ] & =  \e \Big[ \Big|   \sum_{r = 1}^k S_{2, r} \Big|^{ 2 + 2 \tau  } \Big] \nn \\
& \leq C  \Big\{ \Big(  \sum_{r : \, odd} \e [ S_{2, r}^2 ] \Big)^{ 1 + \tau  } + \Big(   \sum_{r : \, even} \e [ S_{2, r}^2 ] \Big)^{ 1 + \tau  } \nn \\
& \qquad  +    \sum_{r = 1}^k \e [ | S_{2, r} |^{ 2 + 2 \tau  } ]  \Big\}. 
\end{align*}
In addition, it follows from the basic inequality \eqref{ineq-m} that 
\begin{align*}
\e [ S_{2, r}^2 ]   &  \leq ( m_1 + 1 )    \sum_{ i \in H_r}  \e  [ X_{1, i}^2 X_{2, i}^2 ] \leq  ( m_1 + 1 ) \sum_{ i \in H_r }  \e [ X_{1, i}^4 ] , \nn \\ 
\e [ | S_{2, r} |^{ 2 + 2 \tau }  ] &  \leq  ( m_1 + 1 )^{ 1 + 2 \tau  }  \sum_{ i \in H_r }  \e [ | X_{1, i} X_{2, i} |^{2 + 2 \tau} ] \\
& \leq  ( m_1 + 1 )^{ 1 + 2 \tau }   \sum_{ i \in H_r}  \e [  | X_{1, i} |^{ 4 + 4 \tau }  ].
\end{align*}
Thus an application of the same argument as in \eqref{L1} results in 
\begin{align}
\e [ | X_1^T X_2 |^{ 2 + 2 \tau } ]  & \leq  C ( m_1 + 1 )^{ 1 + \tau } p ^{ \tau }    \sum_{i = 1}^{p} \e [  | X_{1, i} |^{ 4 + 4 \tau } ], \label{L2}
\end{align}
wich together with \eqref{L1} entails that under condition \eqref{cor1-c1}, we have 
\begin{align}
L_{x, \tau}   =      \e  (  | \Vert X \Vert^2 - \e \Vert X \Vert^{ 2} |^{ 2 + 2 \tau   } ) + \e ( | X_1^T X_2 |^{2 + 2 \tau  }  )     \leq C  \kappa_1  ( m_1 + 1 )^{ 1 + \tau } p^{ 1 + \tau }.   \label{L-tau}
\end{align}

Next we deal with term $  \e [  ( X_1^T \Sigma_x X_2 )^2  ] $. Denote by $ \sigma_{i j} $ the $ (i, j) $th entry of matrix $ \Sigma_x $. By \eqref{cor1-c4} and the $ m_1 $-dependent structure, it holds that 
\begin{align}
\e [ ( X_1^T \Sigma_x X_2 )^2 ] &  =   \e \Big[ \Big(  \sum_{i =1 }^{p}  \sum_{ | j - i | \leq m_1  }  \sigma_{i j}  X_{1, i} X_{2, j}   \Big)^2  \Big]  \nn \\
& =    \sum_{i =1 }^{p}   \sum_{ u = 1}^{p}  \sum_{ | j - i | \leq m_1  }   \sum_{ | v - u | \leq m_1  }  \sigma_{i j} \sigma_{u v}  \e ( X_{1, i} X_{1, u} ) \e ( X_{2, j} X_{2, v} )   \nn \\
& \leq  \kappa_4^2  \sum_{i =1 }^{p}   \sum_{ | u - i | \leq m_1 }  \sum_{ | j - i | \leq m_1  }   \sum_{ | v - u | \leq m_1  }   | \e ( X_{1, i} X_{1, u} ) | | \e ( X_{2, j} X_{2, v} )  |  \nn \\ 
& \leq C  \kappa_4^2 ( m_1 + 1  )^3  p \kappa_{4}^2 = C  \kappa_4^4 ( m_1 + 1  )^3 p. \label{G0}  
\end{align}
Similar results as in \eqref{L-tau} and \eqref{G0} also hold for $Y$. That is, 
\begin{gather}
L_{ y, \tau } \leq   C\kappa_1 ( m_2 + 1  )^{ 1 + \tau } q^{ 1 + \tau },  \label{L-y} \\
\e [ ( Y_1^T \Sigma_y Y_2 )^2 ] \leq C \kappa_4^4 ( m_2 + 1  )^3 q . \label{G0y}
\end{gather}

As a consequence, under conditions \eqref{cor1-c1}--\eqref{cor1-c4} there exists some positive constant $ C_{\kappa} $ depending on $  \kappa_1 , \kappa_2, \kappa_3$, and $ \kappa_4 $ such that 
\begin{align*}
B_X^{ - 2 \tau  } L_{ x, \tau } / \e [ ( X_1^T X_2 )^2 ]   \leq  \frac { C_{\kappa} ( m_1 + 1 )^{ 1 + \tau } } {    p^{ \tau } } \rightarrow 0, 
\end{align*}
\begin{align*}
B_Y^{ - 2 \tau  } L_{ y, \tau  } / \e [ ( Y_1^T Y_2 )^2 ] \leq   \frac {  C_{\kappa} ( m _2 + 1  )^{ 1 + \tau  } }  {   q^{ \tau } }  \rightarrow 0 ,
\end{align*}
and 
\begin{align*}
\frac {  n^{ - \tau } L_{ x, \tau } L_{ y, \tau }  } {  \big(  \e [ ( X_1^T X_2 )^2  ]  \e [ ( Y_1^T Y_2 )^2  ]   \big)^{ 1 + \tau }  } \leq  \frac {  C_{\kappa}  ( m_1 + 1  )^{ 1 + \tau }  ( m_2 + 1 )^{ 1 + \tau }  } {     n^{ \tau }  },
\end{align*}  
\begin{align*}
\frac { \e [ ( X_1^T \Sigma_x X_2 )^2 ]  + B_X^{ - 2 \tau  } L_{x, \tau}^{ (2 + \tau) / ( 1 + \tau  ) }  } {  (  \e [ ( X_1^T X_2 )^2 ]  )^2 }  \leq C_{\kappa}  ( m_1 + 1  )^{ 2 + \tau }  p^{ - \tau },
\end{align*}
\begin{align*}
\frac { \e [ ( Y_1^T \Sigma_y Y_2 )^2 ]  + B_Y^{ - 2 \tau  } L_{y, \tau}^{ (2 + \tau) / ( 1 + \tau  ) }  } {  (  \e [ ( Y_1^T Y_2 )^2 ]  )^2 }  \leq C_{\kappa}   ( m_2 + 1  )^{ 2 + \tau }  q^{ - \tau }.
\end{align*}
Hence by Theorem \ref{thm4}, we see that \eqref{rate-dep} holds for $ 0 < \tau \leq 1 / 2 $.

We next prove the result for the case of $ 1/ 2 < \tau \leq 1 $. By the  previous analysis, it holds that 
\begin{align*}
B_X^{ -1  } L_{ x, 1/2 } / \e [ ( X_1^T X_2 )^2 ] \leq \frac {  C_{\kappa}  ( m_1 + 1  )^{ 3/2 }  } {   p^{ 1/2 } } \rightarrow  0, \\
B_Y^{ -1  } L_{ y, 1/2 } / \e [ ( Y_1^T Y_2 )^2 ] \leq \frac { C_{\kappa} ( m_2 + 1  )^{ 3/2 }  } {  q^{ 1/2 }  } \rightarrow 0  , 
\end{align*}
where the convergence to zero is by the assumption of $ m_1 = o ( p^{ \tau / ( 2 + \tau  ) } ) $ and $ m_2 = o( q^{ \tau / ( 2 + \tau  ) } ) $.
In view of Theorem \ref{thm8} in Section \ref{SecD.1}, it suffices to calculate $ \sum_{ i = 1 }^3  \mathscr{G}_i (X) $ and $ \sum_{ i = 1 }^3 \mathscr{G}_i (Y) $, where
\begin{align*}
\mathscr{G}_1 (X)    &=     \big| \e [ ( X_1^T X_2 )^2  X_1^T \Sigma_x^2  X_2  ] \big|  ,   \\
\mathscr{G}_2 ( X )     &=   \e [  \Vert X_1 \Vert^2 ( X_1^T \Sigma_x X_2 )^2   ],   \\   
\mathscr{G}_3 ( X )     & = \e [ X^T X X^T  ] \Sigma_x^2 \e [ X X^T X ] .            
\end{align*}

Let us begin with considering term $  \mathscr{G}_1 (X) $. Note that
\begin{align*}
\mathscr{G}_1 (X)  & = \big | \e [ ( X_1^T X_2 )^2 X_1^T \Sigma_x^2  X_2 ] \big| \nn \\
&  \leq  ( \e [ | X_1^T X_2 |^{ 2 + 2 \tau } ] )^{  1 / ( 1 +  \tau  ) }  ( \e [ | X_1 ^T \Sigma_x^2 X_2 |^{  ( 1 + \tau ) / \tau  } ] )^{ \tau / ( 1 +  \tau  ) }.
\end{align*}
It follows from \eqref{L2} and assumption \eqref{cor1-c1} that
\begin{align}
( \e [ | X_1^T X_2 |^{ 2 + 2 \tau } ] )^{  1 / ( 1 +  \tau  ) }  
& \leq C  \kappa_1^{ 1/ ( 1 + \tau ) } ( m_1 + 1 ) p.  \label{G1-P1}
\end{align}
Then we analyze term $ \e [ | X_1^T \Sigma_x^2 X_2 |^{ (1 + \tau ) / \tau } ] $. Denote by $ X_{1, H_r} $ the $ r $th block of $ X_1 $ for $ 1 \leq r \leq k $,  and $ \Sigma_{i, j} $ the $ (i, j) $th block of $ \Sigma_x $ for $ 1 \leq i, j  \leq k $. In particular, let $ \Sigma_{1, 0} $ and $  \Sigma_{k, k+1} $ be zero matrices. By the $ m_1 $-dependent structure, $ \Sigma_x $ is a tridiagonal block matrix and thus 
\begin{align*}
\e [ | X_1 ^T \Sigma_x^2 X_2 |^{  ( 1 + \tau )/ \tau  } ]   & = \e \Big[ \Big|    \sum_{r = 1}^k S_{3, r} \Big|^{ ( 1 + \tau )/ \tau } \Big],
\end{align*}
where
\begin{align*}
S_{3, r}  & =   ( \Sigma_{r, r-1} X_{1, H_{r - 1} } + \Sigma_{r, r} X_{1, H_{r } } + \Sigma_{r, r+1} X_{1, H_{r + 1} }  )^T \\
& \hspace{2cm} \cdot ( \Sigma_{r, r-1} X_{2, H_{r - 1} } + \Sigma_{r, r} X_{2, H_{r } } + \Sigma_{r, r+1} X_{2, H_{r + 1} }  )  .
\end{align*}

In addition, $ \{  S_{3, r} , 1 \leq r \leq k \} $ is a $ 3 $-dependent sequence. For simplicity, assume that $ k/ 8 $ is an integer. Then it is easy to see that $ \{  \sum_{ r = 8 ( l - 1 ) + 1 }^{ 8 ( l - 1 ) + 4 } S_{3, r} , 1 \leq l \leq k/8 \} $ and $ \{  \sum_{ r = 8 ( l - 1 ) + 5 }^{ 8 l } S_{3, r} , 1 \leq l \leq k/8 \} $ are sequences of independent random variables. Since $ 2 \leq ( 1 + \tau ) / \tau <  3 $ when $ 1/ 2 < \tau \leq 1  $, it follows from Rosenthal's inequality that 
\begin{align*}
\e [ | X_1 ^T \Sigma_x^2 X_2 |^{  \frac { 1 + \tau } { \tau }  } ]    & \leq C  \e \Big( \Big|   \sum_{ l = 1}^{ k / 8 }  \sum_{  r = 8 (l - 1) + 1   }^{ 8 (l - 1) + 4  }  S_{3, r}  \Big|^{ \frac { 1 + \tau  } { \tau } } \Big)+ C   \e \Big( \Big|  \sum_{ l = 1}^{ k / 8 }  \sum_{  r = 8 (l - 1) + 5   }^{ 8 l  }  S_{3, r}  \Big|^{  \frac { 1 + \tau  } { \tau } } \Big) \nn \\
& \leq C  \bigg\{  \Big( \e \Big[  \Big(  \sum_{ l = 1}^{ k / 8 }  \sum_{  r = 8 (l - 1) + 1   }^{ 8 (l - 1) + 4  }  S_{3, r} \Big)^2 \Big]  \Big)^{ \frac { 1 + \tau } { 2 \tau } } +     \sum_{ l = 1}^{ k / 8 }  \e \Big[ \Big|   \sum_{  r = 8 (l - 1) + 1   }^{ 8 (l - 1) + 4  }  S_{3, r} \Big|^{ \frac { 1 + \tau  } { \tau } } \Big] \nn \\
& \quad +  \Big( \e \Big[  \Big(    \sum_{ l = 1}^{ k / 8 }  \sum_{  r = 8 (l - 1) + 5   }^{ 8 l  }  S_{3, r} \Big)^2 \Big]  \Big)^{  \frac { 1 + \tau } { 2 \tau } } +    \sum_{ l = 1}^{ k / 8 }  \e \Big[ \Big|   \sum_{  r = 8 (l - 1) + 5   }^{ 8 l   }  S_{3, r} \Big|^{ \frac { 1 + \tau  } { \tau } } \Big] \bigg\}.
\end{align*}
Then by inequality \eqref{ineq-m}, we can obtain 
\begin{align*}
\e [ | X_1 ^T \Sigma_x^2 X_2 |^{  \frac { 1 + \tau } { \tau }  } ]  & \leq  C  \bigg\{   \Big(  \sum_{ l = 1}^{ k / 8 }    \sum_{  r = 8 (l - 1) + 1   }^{ 8 (l - 1) + 4  }  \e [  S_{3, r} ^2 ]  \Big)^{  \frac { 1 + \tau } { 2 \tau }  } + \Big(   \sum_{ l = 1}^{ k / 8 }     \sum_{  r = 8 (l - 1) + 5   }^{ 8 l  }  \e [  S_{3, r} ^2 ]  \Big)^{  \frac { 1 + \tau } { 2 \tau }  } \nn \\
& \hspace{6cm}
+   \sum_{r = 1}^k \e [ | S_{3, r} |^{ \frac { 1 + \tau  } { \tau } } ] \bigg\} \nn \\
& \leq C   k^{ \frac { 1 - \tau } { 2 \tau }  } \bigg\{    \sum_{ l = 1}^{ k / 8 }    \sum_{  r = 8 (l - 1) + 1   }^{ 8 (l - 1) + 4  }   \e [ | S_{3, r} |^{ \frac { 1 + \tau  } { \tau } } ] +    \sum_{ l = 1}^{ k / 8 }    \sum_{  r = 8 (l - 1) + 5   }^{ 8 l  }   \e [ | S_{3, r} |^{\frac { 1 + \tau  } { \tau } } ] \bigg\} \nn \\
& \leq C  [ p / ( m_1 + 1 ) ]^{  \frac { 1 - \tau } { 2 \tau }     }  \sum_{r = 1}^k \e [ | S_{3, r} |^{ \frac { 1 + \tau  } { \tau } } ].
\end{align*}

Furthermore, it holds that 
\begin{align*}
\e [ | S_{3, r} |^{  \frac { 1 + \tau  } { \tau }  } ] & \leq C  \e \big[ \big( \Vert  \Sigma_{r, r-1} X_{1, H_{r - 1} } \Vert^2 +  \Vert  \Sigma_{r, r} X_{1, H_{r } } \Vert^2  +  \Vert \Sigma_{r, r+1} X_{1, H_{r + 1} }  \Vert^2   \big)^{    \frac { 1 + \tau  } { \tau }  } \big] \nn \\
& \leq C  \big( \e [  \Vert  \Sigma_{r, r-1} X_{1, H_{r - 1} } \Vert^{  \frac { 2 + 2 \tau } { \tau }  } ] + \e [  \Vert  \Sigma_{r, r} X_{1, H_{r } } \Vert^{  \frac { 2 + 2 \tau } { \tau }  } ] + \e [  \Vert \Sigma_{r, r+1} X_{1, H_{r + 1} }  \Vert^{  \frac { 2 + 2 \tau } { \tau }  } ] \big).   
\end{align*}
For $ 1 \leq i, j \leq m_1 + 1 $, denote by $\Sigma_{r, r- 1}^{ (i, j) }  $ the $ (i, j) $th entry of $  \Sigma_{r, r-1} $ and $  X_{1, H_{r - 1} }^{ (i) } $ the $ i $th component of $ X_{1, H_{r - 1} } $. Observe that by  assumption \eqref{cor1-c4}, we have 
\begin{align*}
& \e [  \Vert  \Sigma_{r, r-1} X_{1, H_{r - 1} } \Vert^ {  \frac { 2 + 2 \tau } { \tau }  } ] \nn \\
& = \e \Big[ \Big| \sum_{i = 1}^{ m_1 + 1 } \sum_{ j = 1 }^{ m_1 + 1} \sum_{ l = 1 }^{ m + 1} \Sigma_{r, r- 1}^{ (i, l) } \Sigma_{r- 1, r}^{ ( l, j ) } X_{1, H_{r - 1} }^{ (i) } X_{1, H_{r - 1} }^{ (j) } \Big|^{  \frac { 1 +  \tau } { \tau }  } \Big] \nn \\
& \leq  \e \Big[ \Big\{  \sum_{i = 1}^{ m_1 + 1 } \sum_{j = 1}^{ m_1 + 1 } \sum_{ l = 1 }^{ m_1 + 1}  [  \e  ( X_{1, H_r }^{ ( i ) } ) ^2 ] ^{1/2}  \e [  ( X_{1, H_{ r  - 1 } }^{ ( l ) } ) ^2 ]  [  \e  ( X_{1, H_r }^{ ( j ) } ) ^2 ] ^{1/2}  | X_{1, H_{r - 1} }^{ (i) } X_{1, H_{r - 1} }^{ (j) } | \Big\}^{  \frac { 1 +  \tau } { \tau }  } \Big]  \\
& \leq  ( m_1 + 1 )^{ \frac { 1 + \tau } { \tau } } \kappa_4^{ \frac {  ( 3 - 2 \tau ) ( 1 + \tau )  } { \tau } }  \e \Big[ \Big\{  \sum_{i = 1}^{ m_1 + 1 } \sum_{j = 1}^{m_1 + 1}   [  \e  ( X_{1, H_r }^{ ( i ) } ) ^2 ] ^{ \tau - \frac 1 2 }  \nn \\
& \hspace{7cm} \times  [  \e  ( X_{1, H_r }^{ ( j ) } ) ^2 ] ^{ \tau - \frac 1 2 }  | X_{1, H_{r - 1} }^{ (i) } X_{1, H_{r - 1} }^{ (j) } | \Big\}^{  \frac { 1 +  \tau } { \tau }  } \Big]. 
\end{align*}
Moreover, it follows from \eqref{ineq-m} that 
\begin{align*}
& \e [    \Vert  \Sigma_{r, r-1} X_{1, H_{r - 1} } \Vert^ {  \frac { 2 + 2 \tau } { \tau }  } ] \nn \\
& \leq    ( m_1 + 1 )^{ \frac { 1 + \tau } { \tau } } \kappa_4^{ \frac {  ( 3 - 2 \tau ) ( 1 + \tau )  } { \tau } }  \e \Big[ \Big(   \sum_{i = 1}^{ m_1 + 1 }  [  \e  ( X_{1, H_r }^{ ( i ) } ) ^2 ] ^{ \tau - \frac 1 2 }    | X_{1, H_{r - 1} }^{ (i) } |   \Big) ^{ \frac { 2 + 2 \tau } { \tau } } \Big]  \nn \\
& \leq   ( m_1 + 1 )^{ \frac { 3 + 2 \tau } { \tau } } \kappa_4^{ \frac {  ( 3 - 2 \tau ) ( 1 + \tau )  } { \tau } }      \sum_{i = 1}^{ m_1 + 1 }   [  \e  ( X_{1, H_r }^{ ( i ) } ) ^2 ] ^{  ( \tau - \frac 1 2 ) ( 2  +  2 \tau ) / \tau }  \e [ | X_{1, H_{r - 1} }^{ (i) } |^{ ( 2  + 2 \tau  ) / \tau }  ] \nn \\
& \leq  ( m_1 + 1 )^{ \frac { 3 + 2 \tau } { \tau } } \kappa_4^{ \frac {  ( 3 - 2 \tau ) ( 1 + \tau )  } { \tau } }     \sum_{i = 1}^{ m_1 + 1 }   [  \e  ( X_{1, H_r }^{ ( i ) } ) ^{ 4 + 4 \tau  } ]^{ 1 - \frac { 1 } { 2 \tau } }  [  \e  ( X_{1, H_{ r - 1 } }^{ ( i ) } ) ^{ 4 + 4 \tau  } ]^{ \frac { 1 } { 2 \tau }   }.
\end{align*}

Note that for any $ a > 0$, $b > 0$, and $0 < \alpha < 1$, we have 
\begin{equation}  
a^{ 1 - \alpha} b ^{ \alpha} \leq a + b . \label{basic}
\end{equation} 
Thus it holds that 
\begin{align*}
& \e [  \Vert  \Sigma_{r, r-1} X_{1, H_{r - 1} } \Vert^ {  \frac { 2 + 2 \tau } { \tau }  } ] \nn \\
& \leq    ( m_1 + 1 )^{ \frac { 3 + 2 \tau } { \tau } } \kappa_4^{ \frac {  ( 3 - 2 \tau ) ( 1 + \tau )  } { \tau } }  \Big(  \sum_{ i \in H_r } \e [ | X_{1, i} |^{ 4 + 4 \tau } ] + \sum_{ i \in H_{ r - 1 } } \e [ | X_{1, i} |^{ 4 + 4 \tau } ] \Big).
\end{align*} 
In the same manner, we can deduce 
\begin{align*}
\e [  \Vert  \Sigma_{r, r} X_{1, H_{r } } \Vert^ {  \frac { 2 + 2 \tau } { \tau }  } ] & \leq   2   (  m_1 + 1  )^{ \frac { 3 + 2 \tau } { \tau } } \kappa_4^{ \frac {  ( 3 - 2 \tau ) ( 1 + \tau )  } { \tau } }    \sum_{ i \in H_r } \e [ | X_{1, i} |^{ 4 + 4 \tau } ]   \nn \\
\e [  \Vert \Sigma_{r, r+1} X_{1, H_{r + 1} }  \Vert^ {  \frac { 2 + 2 \tau } { \tau }  }  ] & \leq    ( m_1 + 1 )^{ \frac { 3 + 2 \tau } { \tau } } \kappa_4^{ \frac {  ( 3 - 2 \tau ) ( 1 + \tau )  } { \tau } }  \Big(   \sum_{ i \in H_r } \e [ | X_{1, i} |^{ 4 + 4 \tau } ] + \sum_{ i \in H_{ r + 1 } } \e [ | X_{1, i} |^{ 4 + 4 \tau } ] \Big). 
\end{align*}
Thus by \eqref{cor1-c1}, it holds that 
\begin{align*}
\e [ | X_1^T \Sigma_x^2 X_2 |^{ \frac { 1 + \tau } { \tau } }  ] & \leq  C  [ p / ( m_1 + 1 ) ]^{  \frac { 1 - \tau } { 2 \tau }     }    ( m_1 + 1 )^{ \frac { 3 + 2 \tau } { \tau } } \kappa_4^{ \frac {  ( 3 - 2 \tau ) ( 1 + \tau )  } { \tau } }  \sum_{r = 1}^k  \sum_{ i \in H_r } \e [ X_{1, H_{r} } ^{ (i) } ]^{ 4 + 4 \tau } \nn \\
& =  C   \kappa_4^{ \frac {  ( 3 - 2 \tau ) ( 1 + \tau )  } { \tau } }   ( m_1 + 1 )^{ \frac { 5 + 5 \tau } { 2 \tau }  }   p^{ \frac { 1 - \tau } {  2 \tau } }  \sum_{i = 1}^{ p } \e [ | X_{1, i} |^{ 4 + 4 \tau }  ]  \nn \\
& \leq   C  \kappa_1  \kappa_4^{ \frac {  ( 3 - 2 \tau ) ( 1 + \tau )  } { \tau } }   ( m_1 + 1 )^{ \frac { 5 + 5 \tau } { 2 \tau }  }   p^{ \frac { 1 + \tau } {  2 \tau } }, 
\end{align*}
which together with \eqref{G1-P1} leads to 
\begin{align}
\mathscr{G}_1 (X) & \leq C  \kappa_1 \kappa_{4}^{ 3 - 2 \tau } ( m_1 + 1 )^{ 7/2 } p^{ 3/2 } .   \label{G1}
\end{align}

We proceed with bounding term $ \mathscr{G}_2 (X) $. Denote by $ \sigma_{i, j} $ the $ (i, j) $th entry of matrix $ \Sigma_x $. Under the $ m_1 $-dependent structure, we have 
\begin{align*}
\mathscr{G}_2 (X)
& =   \sum_{l = 1}^{p}   \sum_{ i= 1 }^{p}  \sum_{ u= 1 }^{p }  \sum_{ | j - i | \leq m_1 } \sum_{ | v - u | \leq m_1 }  \sigma_{i, j} \sigma_{u, v} \e [ X_{1, l}^2 X_{1, i} X_{1, u} ] \e [ X_{2, j} X_{2, v} ]. 
\end{align*}
Observe that $ \e [ X_{2, j} X_{2, v} ] = 0 $ if $ | j - v | > m_1 $. Thus it follows that 
\begin{align*}
\mathscr{G}_{ 2 } (X) & \leq   \sum_{l = 1}^{ p }  \sum_{ i= 1 }^{ p }   \sum_{ | u - i | \leq 3 m_1 }   \sum_{ | j - i | \leq m_1 } \sum_{ | v - i | \leq 2 m_1 }  \sigma_{i, j} \sigma_{u, v} \e [ X_{1, l}^2 X_{1, i} X_{1, u} ] \e [ X_{2, j} X_{2, v} ] .
\end{align*}
By the Cauchy--Schwarz inequality, we can obtain
\begin{align*}
\mathscr{G}_{ 2 } (X)  & \leq  \kappa_{4}^{ 3 - 2 \tau }   \sum_{l = 1}^{ p }  \sum_{ i= 1 }^{ p }   \sum_{ | u - i | \leq 3 m_1 }   \sum_{ | j - i | \leq m_1 } \sum_{ | v - i | \leq 2 m_1 } \Big[ (  \e [ | X_{1, l} |^{ 4 + 4 \tau } ] )^{ \frac {1} { 2 + 2 \tau } }   ( \e [ | X_{1, i} |^{ 4 + 4 \tau }  ] ) ^{ \frac { \tau } { 2 + 2 \tau }  } \nn \\
& \hspace{3cm} \times ( \e [ | X_{1, u} |^{ 4 + 4 \tau }  ] ) ^{ \frac { \tau } { 2 + 2 \tau }  }     ( \e [ | X_{1, j} |^{ 4 + 4 \tau }  ] ) ^{ \frac { 1 } { 4 + 4 \tau }  }  ( \e [ | X_{1, v} |^{ 4 + 4 \tau }  ] ) ^{ \frac { 1 } { 4 + 4 \tau }  } \Big]      \nn \\
& \leq  \kappa_{4}^{ 3 - 2 \tau }  \Big(  \sum_{l = 1}^{ p } (  \e [ | X_{1, l} |^{ 4 + 4 \tau } ] )^{ \frac {1} { 2 + 2 \tau } }  \Big)  \Big\{  \sum_{ i= 1 }^{ p }   ( \e [ | X_{1, i} |^{ 4 + 4 \tau }  ] ) ^{ \frac { \tau } { 2 + 2 \tau }  } \nn \\
& \hspace{2cm} \times   \Big(  \sum_{ | j - i | \leq 2 m_1 } ( \e [ | X_{1, j} |^{ 4 + 4 \tau }  ] ) ^{ \frac { 1 } { 4 + 4 \tau }  }   \Big)^2  \Big(  \sum_{ | u - i | \leq 3 m_1 } ( \e [ | X_{1, u} |^{ 4 + 4 \tau }  ] ) ^{ \frac { \tau } { 2 + 2 \tau }  }     \Big) \Big\}.
\end{align*}

Further, by the basic inequality \eqref{ineq-m} it holds that 
\begin{align*}
& \mathscr{G}_{ 2 } (X)   \nn \\
& \leq   \kappa_{4}^{ 3 - 2 \tau }  \Big(  \sum_{l = 1}^{p} (  \e [ | X_{1, l} |^{ 4 + 4 \tau } ] )^{ \frac {1} { 2 + 2 \tau } }  \Big)  \Big( \sum_{i = 1}^{p}   \e [ | X_{1, i} |^{ 4 + 4 \tau } ]   \Big)^{ \frac { \tau } {  2 + 2 \tau } }  \nn \\                      
& \quad   \times  \Big\{  \sum_{ i= 1 }^{p}  \Big(   \sum_{ | j - i | \leq 2 m_1 } ( \e [ | X_{1, j} |^{ 4 + 4 \tau }  ] ) ^{ \frac { 1 } { 4 + 4 \tau }  }   \Big)^{ \frac { 4 + 4 \tau } { 2 + \tau } }  \Big(  \sum_{ | u - i | \leq 3 m_1  } ( \e [ | X_{1, u} |^{ 4 + 4 \tau }  ] ) ^{ \frac { \tau } { 2 + 2 \tau }  }     \Big)^{ \frac { 2 + 2 \tau } { 2 + \tau } } \Big\}^{ \frac { 2 + \tau } { 2 + 2 \tau } }  \nn \\      
& \leq   C   \kappa_{4}^{ 3 - 2 \tau } ( m_1  + 1 )^{ \frac { 1 + 2 \tau } { 1 + \tau } } p^{ \frac { 1 + 2 \tau } { 2 + 2 \tau } } \Big(  \sum_{i = 1}^{p}   \e [ | X_{1, i} |^{ 4 + 4 \tau } ]   \Big)^{  1/2 }  \Big\{  \sum_{ i= 1 }^{p}  \Big(   \sum_{ | j - i | \leq 2 m_1 } ( \e [ | X_{1, j} |^{ 4 + 4 \tau }  ] ) ^{ \frac { 1 } { 2 + \tau } }     \Big) \nn \\
& \hspace{6cm}  \times \Big(    \sum_{ | u - i | \leq 3 m_1 } ( \e [ | X_{1, u} |^{ 4 + 4 \tau }  ] ) ^{ \frac { \tau } { 2 + \tau } }     \Big) \Big\}^{ \frac { 2 + \tau } { 2 + 2 \tau } }.
\end{align*} 
Hence it follows from the basic inequality \eqref{basic} that
\begin{align*}
( \e [ | X_{1, j} |^{ 4 + 4 \tau }  ] ) ^{ \frac { 1 } { 2 + \tau } } ( \e [ | X_{1, u} |^{ 4 + 4 \tau }  ] )^{ \frac { \tau } { 2 + \tau } }  & \leq  \big( \e [ | X_{1, j} |^{ 4 + 4 \tau }  ] +   \e [ | X_{1, u} |^{ 4 + 4 \tau }  ] \big)^{ \frac { 1 + \tau } { 2 + \tau } } \nn \\
& \leq C    \big(   \e [ | X_{1, j} |^{ 4 + 4 \tau }  ] \big)^{  \frac { 1 + \tau } { 2 + \tau }  } +   C   \big(   \e [ | X_{1, u} |^{ 4 + 4 \tau }  ] \big)^{  \frac { 1 + \tau } { 2 + \tau }  },
\end{align*}
which together with \eqref{ineq-m} and assumption \eqref{cor1-c1} yields
\begin{align}
\mathscr{G}_2 (X) & \leq  C  \kappa_{2}^{ 3 - 2 \tau } ( m_1  + 1 )^{ \frac { 4 + 5 \tau } { 2 + 2 \tau } } p^{ \frac { 1 + 2 \tau } { 2 + 2 \tau } } \Big(   \sum_{i = 1}^{p}   \e [ | X_{1, i} |^{ 4 + 4 \tau } ]   \Big)^{  1/2 } \nn \\
& \qquad \times   \Big\{  \sum_{ i= 1 }^{p}  \sum_{ | j - i | \leq 3 m_1 }   \big(   \e [ | X_{1, j} |^{ 4 + 4 \tau }  ] \big)^{  \frac { 1 + \tau } { 2 + \tau }  }  \Big\}^{ \frac { 2 + \tau } { 2 + 2 \tau }  }  \nn \\
& \leq C \kappa_{2}^{ 3 - 2 \tau } ( m_1  + 1 )^{  3  } p  \Big(   \sum_{i = 1}^{p}   \e [ | X_{1, i} |^{ 4 + 4 \tau } ]   \Big)^{  1/2 }   \Big(  \sum_{ j = 1 }^{p}     \e [ | X_{1, j} |^{ 4 + 4 \tau }  ]   \Big)^ { 1/ 2 } \nn \\
& = C  \kappa_{2}^{ 3 - 2 \tau } ( m_1  + 1 )^{  3  } p   \sum_{i = 1}^{p}   \e [ | X_{1, i} |^{ 4 + 4 \tau } ] \leq C  \kappa_1 \kappa_{4}^{ 3 - 2 \tau } ( m + 1 )^3 p^2. \label{G2}
\end{align}

As for term $ \mathscr{G}_3 ( X ) $, we exploit similar arguments. It is easy to see that the $ r $th block of $ \e [ X_1^T X_1 X_1^T ] $ is given by 
\begin{align*}
( \e [ X_1^T X_1 X_1^T ] )^{ (r) } =  \e \big[  X_{1, H_r} ( \Vert X_{ 1, H_{r - 1} } \Vert^2 +  \Vert X_{ 1, H_{r } } \Vert^2 +  \Vert X_{ 1, H_{r + 1} } \Vert^2 ) \big]. 
\end{align*}
Thus the $ r $th block of $ \Sigma_x \e [ X_1 X_1^T X_1 ] $ is 
\begin{align*}
( \Sigma_x \e [ X_1 X_1^T X_1 ] )^{ (r) } & = \Sigma_{r, r - 1 } \e \big[  X_{1, H_{r - 1} } ( \Vert X_{1 , H_{ r - 2} } \Vert^2  +  \Vert X_{1 , H_{ r - 1} } \Vert^2 +  \Vert X_{1 , H_{ r } } \Vert^2  ) \big] \nn \\
& \quad + \Sigma_{r, r  } \e \big[ X_{1, H_{r } } ( \Vert X_{1 , H_{ r - 1} } \Vert^2  +  \Vert X_{1 , H_{ r } } \Vert^2 +  \Vert X_{1 , H_{ r + 1 } } \Vert^2  ) \big]  \nn \\
& \quad + \Sigma_{r, r + 1 }  \e \big[ X_{1, H_{r + 1} } ( \Vert X_{1 , H_{ r } } \Vert^2  +  \Vert X_{1 , H_{ r + 1} } \Vert^2 +  \Vert X_{1 , H_{ r + 2 } } \Vert^2  ) \big]. 
\end{align*}
Then it follows that 
\begin{align}
\mathscr{G}_3 ( X ) & = \Vert \Sigma_x \e [ X_1 X_1^T X_1 ]  \Vert^2  =  \sum_{ r = 1}^k \Vert ( \Sigma_x \e [ X_1 X_1^T X_1 ] )^{ (r) } \Vert^2      \nn \\
& \leq C  \sum_{ r = 1}^k \sum_{u \in \{ r-1, r, r+1 \}} \Big\{   \Vert \Sigma_{r, u } \e (  X_{1, H_{ u } }  \Vert X_{1 , H_{ u - 1 } } \Vert^2 )  \Vert^2   \nn \\
& \hspace{1cm}  +   \Vert \Sigma_{r, u } \e (  X_{1, H_{ u } }  \Vert X_{1 , H_{ u } } \Vert^2 )  \Vert^2 +  \Vert \Sigma_{r, u } \e (  X_{1, H_{ u } }  \Vert X_{1 , H_{ u + 1 } } \Vert^2 )  \Vert^2  \Big\}. \label{G4-med}
\end{align}
In fact, the terms on the right hand side of the above inequality share the same bounds. Thus we show the analysis only for the first term. 

Observe that 
\begin{align*}
&   \sum_{ r = 1}^k  \Vert \Sigma_{r, r - 1 } \e (  X_{1, H_{ r -1 } }  \Vert X_{1 , H_{ r - 2 } } \Vert^2 ) \Vert^2 \nn \\
& =   \sum_{ r = 1}^k \sum_{ i \in H_r} \sum_{ j  \in H_{r} } \sum_{ l  \in H_{r - 1 }} \Sigma_{r, r - 1}^{ ( i, l ) } \Sigma_{r - 1, r}^{ ( l ,j ) } \e [ X_{1, H_{r - 1} }^{ (i) } \Vert X_{1 , H_{ r - 2 } } \Vert^2 ]  \e [ X_{1, H_{r - 1} }^{ (j) } \Vert X_{1 , H_{ r - 2 } } \Vert^2 ] .
\end{align*}
Then it follows from the Cauchy--Schwarz inequality, assumption \eqref{cor1-c4}, and the basic inequality \eqref{ineq-m} that 
\begin{align*}
&   \sum_{ r = 1}^k  \Vert \Sigma_{r, r - 1 } \e (  X_{1, H_{ r -1 } }  \Vert X_{1 , H_{ r - 2 } } \Vert^2 ) \Vert^2 \nn \\
& \leq \kappa_4^{3 - 2 \tau }  \sum_{ r = 1}^k \Big( \sum_{ l \in H_{r - 1} } (\e [ X_{1, l}^2 ] )^{ \frac { 2 \tau - 1  } { 2 } } \Big) \Big( \sum_{ i \in H_{r} }  \big| \e [ X_{1, i } \Vert X_{1 , H_{ r - 2 } } \Vert^2 ]  \big| \Big)  \nn \\
& \hspace{5cm}  \times \Big(    \sum_{ j \in H_r }  ( \e [ X_{1, j}^2 ] )^{ \frac { 2 \tau - 1 } { 2 } } \big|   \e [ X_{1, j} \Vert X_{1, H_{ r - 2 } } \Vert^2  ] \big|   \Vert \Big) \nn \\
& \leq \kappa_4^{3 - 2 \tau } ( m_1 + 1 )    \sum_{ r = 1}^k \Big( \sum_{ l \in H_{r - 1} } (\e [ X_{1, l}^2 ] )^{  2 \tau  - 1   } \Big)^{ \frac 1 2  }  \Big( \sum_{ j \in H_r }  ( \e [ X_{ 1, j }^2 ] )^{ 2 \tau - 1 } \Big)^{ \frac 1 2 }    \e [  \Vert X_{1, H_{r}} \Vert^2 \Vert X_{1 , H_{ r - 2 } } \Vert^4 ] \nn \\
& \leq   \kappa_4^{3 - 2 \tau } ( m_1 + 1 )     \sum_{ r = 1}^k  \Big\{  \sum_{ l \in H_{r - 1}  \cup H_{r} } (\e [ X_{1, l}^2 ] )^{  2 \tau  - 1   }    \Big\}   \e [  \Vert X_{1, H_{r}} \Vert^2 \Vert X_{1 , H_{ r - 2 } } \Vert^4 ] .
\end{align*}
Moreover, note that for any $ a, b > 0$, we have 
\begin{align*}
a b^2 \leq a^3 + b^3.
\end{align*}
Thus in light of \eqref{ineq-m}, we can obtain
\begin{align*}
\e [  \Vert X_{1, H_{r}} \Vert^2 \Vert X_{1 , H_{ r - 2 } } \Vert^4 ]  & \leq  \e[ \Vert X_{1 , H_{ r  } } \Vert^6 ] + \e[ \Vert X_{1 , H_{ r - 2 } } \Vert^6 ]    \nn \\ 
& \leq   ( m_1 + 1 )^2  \Big(  \sum_{ i \in H_{r} } \e [ X_{1, i }^6 ] + \sum_{ i \in H_{r - 2} }    \e [ X_{1, i }^6 ] \Big) .  
\end{align*}

Furthermore, it follows from \eqref{ineq-m} that 
\begin{align*}
& \Big(   \sum_{ l \in H_{r - 1} \cup H_r } (\e [ X_{1, l}^2 ] )^{ 2 \tau - 1 } \Big) \Big( \sum_{ i \in H_{r} }    \e [ X_{1, i }^6 ] \Big) \nn \\
& \leq  2  ( m_1 + 1 ) \big(   \sum_{ l \in H_{ r - 1} \cup H_r } \e [ | X_{1, l} |^{ 4 + 4 \tau } ] \big)^ { \frac { 2 \tau  - 1 } { 2 + 2 \tau } }  \big(   \sum_{ i \in H_{ r  } } \e [ | X_{1, i} |^{ 4 + 4 \tau } ] \big)^ {  \frac { 3 } { 2 +  2 \tau } } \nn \\
& \leq 4 ( m_1 + 1 ) \big\{   \sum_{ l \in H_{ r - 1} } \e [ | X_{1, l} |^{ 4 + 4 \tau } ] + \sum_{ i \in H_{ r  } } \e [ | X_{1, i} |^{ 4 + 4 \tau } ]    \big\}.
\end{align*}
Similarly, we can deduce 
\begin{align*}
& \big(  \sum_{ l \in H_{r  - 1} \cup H_r } (\e [ X_{1, l}^2 ] )^{ 2 \tau - 1 } \big) \big( \sum_{ i \in H_{r - 2} }    \e [ X_{1, i }^6 ]   \big) \nn \\
& \leq 2 ( m_1 + 1 ) \Big\{    \sum_{ l \in H_{ r } } \e [ | X_{1, l} |^{ 4 + 4 \tau } ] + \sum_{ l \in H_{ r - 1} } \e [ | X_{1, l} |^{ 4 + 4 \tau } ] + \sum_{ i \in H_{ r -2 } } \e [ | X_{1, i} |^{ 4 + 4 \tau } ]    \Big\}.
\end{align*}
Consequently, it holds that 
\begin{align*}
&    \sum_{ r = 1}^k  \Vert \Sigma_{r, r - 1 } \e (  X_{1, H_{ r -1 } }  \Vert X_{1 , H_{ r - 2 } } \Vert^2 )  \Vert^2  \nn \\
& \leq  C \kappa_4^{3 - 2 \tau } ( m_1 + 1 )^4  \sum_{ r = 1}^k \sum_{ i \in H_r } \e [ | X_{1, i} |^{ 4 + 4 \tau } ]  \nn \\
& =  C \kappa_4^{ 3 - 2 \tau } ( m_1 + 1 )^4   \sum_{ i = 1}^{p} \e [ | X_{1, i} |^{ 4 + 4 \tau } ] \leq C \kappa_1 \kappa_4^{ 3 - 2 \tau } ( m_1 + 1 )^4 p.
\end{align*}

For the other terms on the right hand side of \eqref{G4-med}, the same bound can be derived in a similar way. Thus we have 
\begin{align}
\mathscr{G}_3 ( X ) &  \leq  C \kappa_1 \kappa_4^{ 3 - 2 \tau  } ( m_1 + 1 )^4 p. \label{G3}
\end{align}
Combining \eqref{G1}, \eqref{G2}, and \eqref{G3}, and noting that $ m_1 + 1 \leq p $, we can obtain
\begin{align}
\mathscr{G}_1 (X) + \mathscr{G}_2 (X) + \mathscr{G}_3 (X) \leq C_{\tau} \kappa_1 \kappa_4^{ 3 - 2 \tau } ( m_1 + 1 )^3 p^2. \label{Gx}
\end{align}
In the same manner, we can also show that 
\begin{align}
\mathscr{G}_1 (Y) + \mathscr{G}_2 (Y) + \mathscr{G}_3 (Y) \leq C_{\tau} \kappa_1 \kappa_4^{ 3 - 2 \tau } ( m_2 + 1 )^3 q^2. \label{Gy}
\end{align}
Hence \eqref{rate-dep}  follows from substituting \eqref{L-tau}--\eqref{G0y} and \eqref{Gx}--\eqref{Gy} into Theorem \ref{thm8}. Then we can see that when $ m_1 $ and $ m_2 $ satisfy \eqref{range-m}, $ T_n \stackrel { \mathscr{D} } { \rightarrow } N( 0, 1 )  $.  This completes the proof of Proposition \ref{cor1}.

\subsection{Proof of Proposition \ref{cor2}} \label{SecB.3}
Assume that $ \Sigma_x = \Gamma_1^T \diag ( \lambda_1^X, \ldots, \lambda_{ p }^X  )  \Gamma_1 $ and $ \Sigma_y = \Gamma_2^T \diag ( \lambda_1^Y , \ldots, \lambda_{ q }^Y ) \Gamma_2  $ for some orthogonal matrices $ \Gamma_1 $ and $ \Gamma_2  $. A useful fact is that the Euclidean norm is invariant to orthogonal transformations. Thus $ X $ and $ Y $ can be replaced with the transformed random vectors $ \breve{X} = \Gamma_1 X $ and $ \breve{Y}  = \Gamma_2 Y $, respectively. Clearly the transformed random vectors are distributed as 
\begin{align*}
\breve{X} \sim  N (0,  \diag( \lambda_1^X , \ldots, \lambda_{ p }^X  ) ) \ \text{ and } \   \breve{Y}   \sim   N(0, \diag (  \lambda_1^Y , \ldots, \lambda_{ q }^Y  ) ).
\end{align*}
It is equivalent to analyze the distance correlation between the new  multivariate normal random variables $ \breve{X} $ and $ \breve{Y} $. Ii is easy to show that 
\begin{gather*}
\max \limits_{  1 \leq i \leq p  } \e [ \breve{X}_{ 1, i }^2 ] \leq a_2   ,  \quad  p^{ - 1 }    \textstyle \sum_{ i = 1 }^{ p } \e [ | \breve{X}_{ 1, i } |^{ 8  }]  \leq C a_2^{ 4 }, \\
p^{ - 1 } \e [ ( \breve{X}_1^T \breve{X}_2 )^2  ]   \geq  a_1^2 , \quad p^{ - 1 } B_X \geq a_1.
\end{gather*}
Similar bounds also hold for $ Y $. Then the conditions of Proposition \ref{cor1} are satisfied and the independence of coordinates entails that $ m_1 = m_2 = 0 $. Therefore, the desired result can be derived by applying Proposition \ref{cor1} with $ \tau = 1 $ and $ m_1 = m_2 = 0 $. This concludes the proof of Proposition \ref{cor2}.

\renewcommand{\thesubsection}{C.\arabic{subsection}}
\section{Some key lemmas and their proofs}
\label{lemmas}

\subsection{Lemma \ref{le-consist} and its proof} \label{SecC.1}

\begin{lemma} \label{le-consist}
	Under condition \eqref{cond1}, we have
	\begin{equation} 
	\V_n^*( X ) / \V^2 ( X )  \longrightarrow  1 \quad \mbox{in probability}  \label{X-con}
	\end{equation}
	and 
	\begin{equation} 
	\V_n^*( Y )  /  \V^2 ( Y )  \longrightarrow    1 \quad \mbox{in probability} \label{Y-con}
	\end{equation}
	as $ n \rightarrow \infty $.
\end{lemma}

\noindent \textit{Proof}. For any $X$ and $Y$, since $\V_n^* ( X, Y) $ is a U-statistic and noting that $ \e [ \V_n^* (X, Y) ]  =  \V^2 (X, Y) $ by \eqref{U-stat}, it follows from the moment inequality of U-statistics \citep[p.~72]{book1994} and conditional Jensen's inequality that for $ 0 < \tau \leq 1 $,
\begin{align}
& \e \big[  \big| \V_n^* ( X, Y )    - \V^2 (X, Y) \big|^{1 + \tau} \big] \nn \\ 
&   \leq  C    \sum_{i = 1}^4 { 4 \choose i }^{1 + \tau} { n \choose i}^{ - \tau }   \e \big[ \big| h ( ( X_1, Y_1 ), ( X_2, Y_2 ), ( X_3, Y_3 ), ( X_4, Y_4 ) )  \big|^{ 1 + \tau } \big] \nn \\
&  \leq   C n^{ - \tau }  \e \big[ \big| h ( ( X_1, Y_1 ), ( X_2, Y_2 ), ( X_3, Y_3 ), ( X_4, Y_4 ) )  \big|^{ 1 + \tau } \big]   .  \label{h-b1}
\end{align}
In fact, the moment of $  h ( ( X_1, Y_1 ), ( X_2, Y_2 ), ( X_3, Y_3 ), ( X_4, Y_4 ) )  $ can be dominated by that of $ d ( X_1, X_2 ) d (Y_1, Y_2) $ based on the expression given in Lemma \ref{h-d} in Section \ref{SecC.5}.

By  choosing $ X = Y $ in \eqref{ex-h-d} and the Cauchy--Schwarz inequality, we can obtain that for $  0 < \tau \leq 1  $, 
\begin{align*}
\e [ | h ( ( X_1, X_1 ), ( X_2, X_2 ), ( X_3, X_3 ), ( X_4, X_4 ) )  |^{ 1 + \tau }  ] \leq C \e [ | d ( X_1, X_2 ) |^{2 + 2 \tau } ].  
\end{align*}
Thus it follows from \eqref{h-b1} that
\begin{align}
\e \big[  \big| \V_n^* ( X ) /  \V^2 ( X )  - 1 \big|^{1 + \tau} \big] \leq \frac { C \e [ | d ( X_1, X_2 ) |^{ 2 + 2 \tau } ] } { n^{ \tau } [ \V^2 ( X ) ]^{ 1 + \tau }  }.  \label{bound-1}
\end{align}
Moreover, since
$
\e [ | d ( Y_1, Y_2 ) |^{2 + 2 \tau} ] \geq    ( \e [ d^2 ( Y_1, Y_2 )  ] )^{ 1 + \tau }  = [ \V^2 ( Y  ) ]^ {1 + \tau }$,
it follows from condition \eqref{cond1} that 
\begin{align*}
\frac {  \e [ | d( X_1, X_2 ) |^{ 2 + 2 \tau } ]  } { n^{\tau}  [ \V^2 ( X ) ]^{1 + \tau}   } \rightarrow 0, 
\end{align*}
which yields the ratio consistency \eqref{X-con}. The result in  \eqref{Y-con} can be obtained similarly. This completes the proof of Lemma \ref{le-consist}.

\subsection{Lemma \ref{le-decom} and its proof} \label{SecC.2}

\begin{lemma} \label{le-decom}
	If $ \e [ \Vert X \Vert ^2  ]  + \e [ \Vert Y \Vert ^2  ]  < \infty $ and $ X $ is independent of $ Y $, then we have 
	\begin{equation*}
	\breve{T}_n  =  W_n^{(1)} (X, Y) + W_n^{(2)}(X, Y),  
	\end{equation*}
	where 
	\begin{equation}
	W_n^{(1)} (X, Y) = \sqrt{ \frac {2 } { n (n - 1)  } } \sum_{1 \leq i < j \leq n }  \frac { d (X_i, X_j) d (Y_i, Y_j)} { \sqrt{ \V^2( X ) \V^2 ( Y )  } }     \label{W_n_1}
	\end{equation}
	and $ W_n^{(2)} (X, Y) $ satisfies
	$
	\e (  [ W_n^{(2)} (X, Y) ]^2 )  \leq  C n^{-1} $. 
\end{lemma}

\noindent \textit{Proof}. Recall that $  	\V_n^* (X,  Y)  $ is a U-statistic and 
$$
\V_n^* (X,  Y) = { n \choose 4 } ^{-1} \sum_{1 \leq i_1 < i_2 < i_3 < i_4 \leq n} h ( ( X_{i_1}, Y_{i_1} ), \ldots, ( X_{i_4}, Y_{i_4} )  ).
$$ 
It has been shown in \citet{HH2017} that under the independence of $ X $ and $ Y $,
\begin{align*}
\e \big[ h ( (X_1, Y_1), (X_2, Y_2), (X_3, Y_3), (X_4, Y_4) ) \big\vert ( X_1, Y_1 ) \big]  & = 0  , \\
\e \big[ h ( (X_1, Y_1), (X_2, Y_2), (X_3, Y_3), (X_4, Y_4) ) \big\vert ( X_1, Y_1 ), (X_2, Y_2) \big]  
&  =  \frac {1} {6} d(X_1, X_2) d(Y_1, Y_2). 
\end{align*}
Thus by the Hoeffding decomposition \citep[e.g.][p.~23]{book1994} and dispersion for U-statistics \citep[p.~31]{book1994}, when $ X $ is independent of $ Y $ we have
\begin{equation*} 
\V_n^*(X, Y) =   { n \choose 2 } ^{-1}  \sum_{ 1 \leq i < j \leq n } d(X_i, X_j) d(Y_i, Y_j) + U_{n} (X, Y),   
\end{equation*}
where
\begin{equation*}
\e [ U_n^2 (X, Y) ] \leq \frac { C } { n^3 } \e [ h ( (X_1, Y_1), (X_2, Y_2), (X_3, Y_3), (X_4, Y_4) ) ]^2.   
\end{equation*}

Furthermore, Lemma \ref{le-h2} given in Section \ref{SecC.6} yields 
\begin{align*}
\e [ U_n^2 (X, Y) ] \leq  \frac {C \V^2 (X) \V^2 (Y)}  {  2 n^3  }.
\end{align*}
Hence it follows from \eqref{N-Tn} that
\begin{equation*}
\e [ W_n^{ (2) } ( X, Y ) ]^2  =   \frac { n ( n - 1 )  \e [ U_n^2 ( X, Y ) ]  } {  2 \V^2 ( X ) \V^2 ( Y ) } \leq \frac { C }{ n },
\end{equation*}  
which concludes the proof of Lemma \ref{le-decom}.

\subsection{Lemma \ref{le-martin} and its proof} \label{SecC.3}

\begin{lemma} \label{le-martin}
	Let $ \mathscr{F}_k   = \sigma \{ (X_1, Y_1), \ldots, (X_k, Y_k) \} $ be a $\sigma$-algebra. Then $ \{ ( \zeta_{n,k}, \mathscr{F}_k ) , k \geq 1 \} $ forms a martingale difference array under the independence of $ X $ and $ Y $, where $ \zeta_{n, k} $ is defined in \eqref{mart-dif}.
\end{lemma}

\noindent \textit{Proof}. It is easy to see that $ \zeta_{n, k} \in \mathscr{F}_k $ and  when $ X $ is independent of $ Y $,
\begin{equation*}
\e \Big[ \sum_{ i = 1}^{k - 1} d(X_i , X_k) d(Y_i, Y_k) \big\vert  \mathscr{F}_{k-1}  \Big]  =   \sum_{i = 1}^{k-1} \e \big[d(X_i , X_k)\big\vert X_i \big] \e \big[  d(Y_i, Y_k)  \big\vert Y_i  \big] =0, 
\end{equation*}
where the last equality is due to $ \e \big[ d(X_1 , X_2 ) \big\vert X_1 \big] = 0  $ and $ \e \big[ d ( X_1 , X_2 ) \big\vert X_2 \big] = 0  $.

\subsection{Lemma \ref{le-MarCLT} and its proof} \label{SecC.4}

\begin{lemma} \label{le-MarCLT}
	If $ \e [ \Vert X \Vert ^{2 + 2 \tau } ] + \e [ \Vert Y \Vert ^{2 + 2 \tau } ] < \infty$ for some constant $ 0 < \tau \leq 1 $ and $ X $ is independent of $ Y $, then we have 
	\begin{align}
	\e \bigg( \bigg|  \sum_{k =1}^n \e [ \zeta_{n,k}^2 \vert \mathscr{F}_{k-1} ]  - 1  \bigg|^{1 + \tau }  \bigg) & \leq  C  \Big ( \frac  {      \e [ g (X_1,X_2, X_3, X_4 ) ]  \e [ g( Y_1, Y_2, Y_3, Y_4 ) ]   }  {  [  \V^2 ( X ) \V^2 ( Y )  ] ^2 }  \Big)^{   ( 1 + \tau  ) / 2   } \nn \\
	& \quad  + \frac { C \e [ | d (X_1, X_2 ) |^{2 + 2 \tau } ]  \e [  | d (Y_1, Y_2) |^{2 + 2 \tau } ] } { n^{ \tau  }  [ \V^2 ( X ) \V^2 ( Y ) ]^{ 1 + \tau  }    } \label{c1} 
	\end{align}
	and 
	\begin{align}
	\sum_{ k = 1}^n \e [ | \zeta_{n,k} |^{2 + 2 \tau } ]    \leq  \frac { C  \e [ | d (X_1, X_2 ) |^{2 + 2 \tau } ]  \e [  | d (Y_1, Y_2) |^{2 + 2 \tau } ] } { n^{ \tau  }  [ \V^2 ( X ) \V^2 ( Y ) ]^{ 1 + \tau   }    }.  \label{c2}
	\end{align}
\end{lemma}

\noindent \textit{Proof}. (i) We first prove \eqref{c1}. Recall the definition of $ \zeta_{n, k} $ in \eqref{mart-dif}. Note that under the independence of $ X $ and $ Y $, we have 
\begin{equation*}
\sum_{k= 1}^n \e [ \zeta_{n, k}^2 \vert \mathscr{F}_{k-1} ]  =  \frac { 2  \sum_{k =1}^n \e \Big(  \Big[ \sum_{i = 1}^{k-1} d (X_i, X_k) d(Y_i, Y_k)  \Big]^2  \Big \vert \mathscr{F}_{k - 1 }  \Big)  } { n ( n - 1 )  \V^2 ( X ) \V^2 ( Y )   }  := R_n^{(1)} + R_n^{(2)} ,  
\end{equation*}
where $ R_n^{(1)} $ is the sum of squared terms given by
$$
R_n^{(1)} =  \frac { 2  \sum_{k =1}^n  \sum_{i = 1}^{k-1}  \e \big[ d^2 (X_i, X_k) \big \vert X_i \big]  \e \big[  d^2 (Y_i, Y_k)   \big \vert  Y_i   \big]   } { n ( n - 1 )  \V^2 ( X ) \V^2 ( Y )   }
$$
and $  R_n^{(2)} $ is the sum of cross-product terms given by 
$$
R_n^{(2)} =  \frac { 4  \sum_{k =1}^n  \sum_{ 1 \leq i < j \leq k - 1 }   \e \big[ d (X_i, X_k)  d (X_j, X_k)  \big \vert  X_i , X_j  \big]   \e \big[  d (Y_i, Y_k) d (Y_j, Y_k)  \big\vert  Y_i, Y_j  \big]  } { n ( n - 1 )  \V^2 ( X ) \V^2 ( Y )   }.  
$$
Thus it holds that 
\begin{align}
\e \bigg( \bigg|  \sum_{k =1}^n \e [ \zeta_{n,k}^2 \vert \mathscr{F}_{k-1} ]  - 1  \bigg|^{1 + \tau }  \bigg) \leq C \big( \e [ | R_n^{ (1) } - 1 |^{ 1 + \tau    } ]  +  \e [ | R_n^{(2)} |^{ 1 + \tau } ]  \big). \label{condvar}
\end{align}

We first bound term  $   \e [ | R_n^{(2)} |^{ 1 + \tau } ] $. Let $ (X', Y') $ be an independent copy of $ (X, Y) $ that is independent of $ (X_1, Y_1), \ldots, ( X_n, Y_n ) $. For notational simplicity, define
$$
\eta_1 (X_i, X_j) =  \e \big[ d (X_i, X )  d (X_j, X )  \big \vert  X_i , X_j  \big] \ \text{ and } \ \eta_2 (Y_i, Y_j) =  \e \big[ d (Y_i, Y )  d (Y_j, Y )  \big \vert  Y_i , Y_j  \big].
$$
By changing the order of summation, we can obtain
\begin{align*}
\e ( [ R_n^{(2)} ] ^2 )  &   =     \frac  {  16   \e \Big[ \sum_{ 1 \leq i < j \leq n } \sum_{ k \geq j+1 }   \eta_1 ( X_i, X_j )  \eta_2 ( Y_i, Y_j )  \Big]^2  }  { n^2 ( n - 1 )^2  [  \V^2 ( X ) \V^2 ( Y )  ]^2  }  \\
& =    \frac  { 16  \e \Big[ \sum_{ 1 \leq i < j \leq n }  (n - j)    \eta_1 ( X_i, X_j )  \eta_2 ( Y_i, Y_j )  \Big]^2  }  { n^2 ( n - 1  )^2 [  \V^2 ( X ) \V^2 ( Y )  ]^2  } .
\end{align*}
In addition, for pairwisely nonequal $ i, j, l $, it holds that 
\begin{align*}
\e [ \eta_1 ( X_i, X_j ) \eta_1 ( X_i, X_l ) ]  
& =    \e \big\{ \e [ d( X_i, X  ) d(X_j, X )  d( X_i, X'  ) d(X_l, X' ) \vert X_i, X_j , X_l ]  \big\} \nn  \\
& =    \e [ d( X_i, X  ) d(X_j, X )  d( X_i, X'  ) d(X_l, X' ) ]    \nn \\
& =    \e \big\{ \e [ d( X_i, X  ) d(X_j, X )  d( X_i, X'  ) d(X_l, X' ) \vert X, X' ]  \big\}  \nn \\
& =    \e \big( \e [ d( X_i, X  ) d(X_i, X' )  \vert X, X' ]  \e [  d( X_j, X'  ) \vert X'  ]   \e [ d(X_l, X' ) \vert X' ]  \big) \nn  \\
& =    0,  \label{cross1}
\end{align*}
where we have used the fact that $ \e [  d( X_j, X'  ) \vert X'  ]  =   \e [ d(X_l, X' ) \vert X' ]  \big) = 0   $.

It is easy to see that $ \e [ \eta_1 (X_i, X_j ) ] = \e [\eta_2 ( Y_i, Y_j )  ] = 0 $ for $ i \neq j $. Thus for pairwisely nonequal $ i, j, k, l $, it holds that 
\begin{equation*}
\e [ \eta_1 ( X_i, X_j ) \eta_1 ( X_k, X_l ) ] = 0.     \label{cross2}
\end{equation*}
Then the cross-product terms in the numerator of $ \e ([ R_n^{(2)} ]^2)  $ vanish. Moreover, in view of the definition of $ g ( X_1, X_2, X_3, X_4 ) $ in \eqref{def-g}, we have
\begin{align*}
\e \big[    \eta_1( X_i, X_j )  \big]^2    
& =   \e \big\{  \e \big[ d (X_i, X )  d (X_j, X ) d (X_i, X' )  d (X_j, X' )   \big \vert  X_i , X_j  \big]  \big\} \\
& =    \e \big[ d (X_1, X_2 )  d (X_1, X_3 ) d (X_2, X_4 )  d (X_3, X_4 )  \big] = \e [ g (X_1,X_2, X_3, X_4 ) ].
\end{align*}
Consequently, it follows that 
\begin{align*}
\e ( [ R_n^{(2)} ]^2 )   &  =      \frac  { 16  \sum_{ 1 \leq i < j \leq n }  (n - j)^2   \e \big[    \eta_1( X_i, X_j )  \big]^2  \e \big[ \eta_2  ( Y_i, Y_j )  \big]^2  }  { n^2 ( n - 1 )^2 [  \V^2 ( X ) \V^2 ( Y )  ]^2  } \nn \\
& =   \frac  { 16  \sum_{j = 1}^n  (j - 1)  (n - j)^2  \e [ g (X_1,X_2, X_3, X_4 ) ]  \e [ g( Y_1, Y_2, Y_3, Y_4 ) ]  }  { n^2 ( n - 1 )^2 [  \V^2 ( X ) \V^2 ( Y )  ]^2  } \nn \\
& \leq  \frac  { C     \e [ g (X_1,X_2, X_3, X_4 ) ]  \e [ g( Y_1, Y_2, Y_3, Y_4 ) ]  }  {  [  \V^2 ( X ) \V^2 ( Y )  ]^2  }.  
\end{align*}
Hence we can obtain
\begin{align}
\e [ | R_n^{ (2) } |^{ 1 + \tau } ]   \leq  \Big(  \frac  { C     \e [ g (X_1,X_2, X_3, X_4 ) ]  \e [ g( Y_1, Y_2, Y_3, Y_4 ) ]  }  {  [  \V^2 ( X ) \V^2 ( Y )  ]^2  }  \Big)^{ (1 + \tau ) /2  } .   \label{Rn_2}
\end{align}

Next we deal with term $ \e  [ | R_n^{(1)} - 1 |^{ 1 + \tau } ] $. Since $ \e [ d^2 (X_1, X_2) ] = \V^2 ( X ) $, clearly when $ X  $ is independent of $ Y $, we have 
$$
\e [ R_n^{(1)}  ]= \frac  { 2 \sum_{k =1}^n \sum_{i=1}^{ k -1}  \V^2 ( X ) \V^2 ( Y ) }  { n (n -1)  \V^2 ( X ) \V^2 ( Y ) } = 1.
$$
For simplicity, denote by $ \eta_3 (X_i, Y_i) =   \e \big[ d^2 (X_i, X) \big \vert X_i \big]  \e \big[  d^2 (Y_i, Y )   \big \vert  Y_i   \big] $. Then by changing the order of summation, we deduce 
\begin{align*}
\e  [  | R_n^{(1)} - 1 |^{ 1 + \tau } ]    & =   \frac { \e \Big[\Big| 2  \sum_{k =1}^n  \sum_{i = 1}^{k-1}  [ \eta_3 (X_i, Y_i)  - \e   \eta_3 (X_i, Y_i) ]  \Big|^{1 + \tau} \Big]} {  [  n ( n - 1 ) \V^2 ( X ) \V^2 ( Y ) ]^{1 + \tau}    } \\
& =   \frac { \e \Big[\Big| 2  \sum_{i =1}^n  ( n - i )  [ \eta_3 (X_i, Y_i)  - \e   \eta_3 (X_i, Y_i) ]  \Big|^{1 + \tau } \Big]} {   [ n ( n -1 ) \V^2 ( X ) \V^2 ( Y ) ]^{1 + \tau}    }. 
\end{align*}
Then it follows from the von Bahr--Esseen inequality \citep[p.~100]{bookLinBai} for independent random variables that when $ 0< \tau \leq 1 $, 
\begin{align*}
\e  [  | R_n^{(1)} - 1 |^{ 1 + \tau } ]    & \leq  \frac { C   \sum_{i = 1}^n ( n - i  )^{1 + \tau} \e [ | \eta_3 ( X_i, Y_i  ) |^{1 + \tau} ]   } {   [ n^2  \V^2 ( X ) \V^2 ( Y ) ]^{1 + \tau}    } \nn  \\
& \leq   \frac { C    \e [ | \eta_3 ( X_i, Y_i  ) |^{1 + \tau} ]   } {   n^{\tau} [ \V^2 ( X ) \V^2 ( Y ) ]^{1 + \tau}    } \nn \\
& \leq  \frac { C    \e [ | d ( X_1, X_2 ) |^{2 + 2 \tau} ]  \e [ | d ( Y_1, Y_2 ) |^{2 + 2 \tau} ]   } {   n^{\tau} [ \V^2 ( X ) \V^2 ( Y ) ]^{1 + \tau}    },
\end{align*}
which along with \eqref{condvar} and \eqref{Rn_2} leads to \eqref{c1}.

(ii) We now show \eqref{c2}. Note that 
\begin{align}
\sum_{ k = 1}^n \e [ | \zeta_{n,k} |^{2 + 2 \tau } ]    = \frac { 2 ^{ 1 + \tau }   \sum_{k=1}^n   \e \Big[\big|  \sum_{i =1 }^ {k-1}  d (X_i, X_k) d(Y_i, Y_k)  \big|^{ 2 + 2 \tau } \Big]}  {   [ n ( n - 1 )  \V^2 ( X ) \V^2 ( Y ) ]^{1 + \tau}  }.   \label{mo-2}
\end{align}
Given $ ( X_k, Y_k ) $,  $ \{ d( X_i, X_k ) d( Y_i, Y_k ), 1 \leq i \leq k-1 \} $ is a sequence of independent random variables and under the independence of $ X $ and $ Y $,
$$ 
\e [ d( X_i, X_k ) d( Y_i, Y_k ) \vert X_k, Y_k ]  = \e [ d( X_i, X_k ) \vert X_k ] \e [ d( Y_i, Y_k ) \vert Y_k ] = 0.
$$
Thus it follows from Rosenthal's inequality for independent random variables that
\begin{align*}
\e \Big[\Big|  \sum_{i =1 }^ {k-1}  d (X_i, X_k) d(Y_i, Y_k)  \Big|^{ 2 + 2 \tau } \Big] & =  \e \bigg[  \e \Big( \Big| \sum_{i =1 }^ {k-1}  d (X_i, X_k) d(Y_i, Y_k)  \Big|^{ 2 + 2 \tau  }  \Big\vert (X_k, Y_k)  \Big) \bigg]  \\
& \leq  C  \e  \bigg[  \e \Big( \Big[  \sum_{i =1 }^ {k-1}  d (X_i, X_k) d(Y_i, Y_k)  \Big]^2  \Big\vert (X_k, Y_k)  \Big) \bigg]^{ 1 + \tau }   \nn \\
& \quad  +  C ( k - 1 )   \e    \big(  | d (X_1, X_2)  |^{ 2 + 2 \tau }   \big) \e \big( | d  (Y_1, Y_2) |^{ 2 + 2 \tau }    \big)  .  
\end{align*}

Since given $ ( X_k, Y_k ) $, $ \{ d( X_i, X_k ) d( Y_i, Y_k ), 1 \leq i \leq k-1 \} $ is a sequence of independent random variables with zero means under the independence of $ X $ and $ Y $, it is easy to see that 
\begin{align*}
\e \Big( \Big[  \sum_{i =1 }^ {k-1}  d (X_i, X_k) d(Y_i, Y_k)  \Big]^2  \Big\vert (X_k, Y_k)  \Big)     & = \sum_{i=1}^{k-1}  \e \big[ d^2 (X_i, X_k) d^2 (Y_i, Y_k) \big \vert (X_k, Y_k) \big], \\
& =  ( k - 1 ) \e \big[ d^2 (X, X_k) \vert X_k \big] \e \big[ d^2 (Y, Y_k) \big \vert Y_k \big].
\end{align*} 
Then it follows from the conditional Jensen's inequality that when $ X  $ is independent of $ Y $, 
\begin{align*}
&  \e \bigg[  \e \Big( \Big[  \sum_{i =1 }^ {k-1}  d (X_i, X_k) d(Y_i, Y_k)  \Big]^2  \Big\vert (X_k, Y_k)  \Big) \bigg]^{1 + \tau}  \\
& \leq  ( k - 1 )^{ 1 + \tau }  \e [ | d (X_1, X_2) |^{ 2 + 2 \tau } ]  \e  [ | d (Y_1, Y_2)  |^{2 + 2 \tau} ]. 
\end{align*}
Finally we can obtain
\begin{align*}
\sum_{ k = 1 }^n  \e \Big[\Big|  \sum_{i =1 }^ {k-1}  d (X_i, X_k) d(Y_i, Y_k)  \Big|^{ 2 + 2 \tau } \Big] 
& \leq    C  \sum_{k=1}^n  ( k - 1 )^{1 + \tau }  \e [ | d (X_1, X_2) |^{ 2 + 2 \tau } ]  \e  [ | d (Y_1, Y_2)  |^{2 + 2 \tau}  ]    \\
& \leq   C  n^{ 2 + \tau }   \e [ | d (X_1, X_2) |^{ 2 + 2 \tau } ]  \e  [ | d (Y_1, Y_2)  |^{2 + 2 \tau}  ].
\end{align*}
Substituting the above bound into \eqref{mo-2} results in \eqref{c2}. This completes the proof of Lemma \ref{le-MarCLT}.

\subsection{Lemma \ref{h-d} and its proof} \label{SecC.5}
The following lemma provides a useful representation of the kernel function $  h ( ( X_1, Y_1 ) , ( X_2, Y_2 ),  \\ ( X_3, Y_3 ), ( X_4, Y_4 ) ) $ in terms of the double-centered distance $ d ( \cdot, \cdot ) $.
\begin{lemma}  \label{h-d}
	For any random vectors $ X $ and $ Y $ with finite first moments, we have 
	\begin{align}
	& h ( ( X_1, Y_1 ) , ( X_2, Y_2 ),  ( X_3, Y_3 ), ( X_4, Y_4 ) ) \nn \\
	& =  \frac 14 \sum_{ \substack {  1 \leq i, j \leq 4,  \\ i \neq j } }  d ( X_i, X_j ) d ( Y_i, Y_j )   - \frac 1 4 \sum_{i = 1}^4  \bigg(  \sum_{ \substack { 1 \leq  j \leq 4, \\ j \neq i  } }    d ( X_i, X_j )   \sum_{ \substack {   1 \leq  j \leq 4, \\ j \neq i  } }  d ( Y_i, Y_j )   \bigg) \nn \\
	&   \quad + \frac {1} {24} \sum_{ \substack { 1 \leq i, j \leq 4,  \\ i \neq j } }  d ( X_i, X_j )   \sum_{ \substack {  1 \leq i, j \leq 4,  \\ i \neq j } }  d ( Y_i, Y_j ) .  \label{ex-h-d}
	\end{align}
\end{lemma}

\noindent \textit{Proof}.  Let us define 
\begin{gather*}
a_1 ( X_1, X_2 ) = \Vert X_1 - X_2 \Vert - \e [  \Vert X_1 - X_2 \Vert  ], \quad a_1 (Y_1, Y_2 ) = \Vert Y_1 - Y_2 \Vert  -  \e [  \Vert Y_1 - Y_2 \Vert  ] ,  \\
a_2 ( X_1  ) = \e [ a_1 ( X_1, X_2 ) | X_1 ], \qquad a_3 ( Y_1 ) = \e [ a_1 ( Y_1, Y_2  ) | Y_1 ].
\end{gather*}
We divide the proof into two steps.

\textit{Step 1}. Recall the definition of $ h ( ( X_1, Y_1 ) , ( X_2, Y_2 ),  ( X_3, Y_3 ), ( X_4, Y_4 ) ) $ given in \eqref{kernel}. It is easy to show that 
\begin{align*}
& \frac 14 \sum_{ \substack {  1 \leq i, j \leq 4,  \\ i \neq j  } }   \Vert X_i - X_j \Vert  \Vert Y_i -  Y_j \Vert  \nn \\ 
& =  \frac 14 \sum_{ \substack {   1 \leq i, j \leq 4,  \\ i \neq j } }  a_1 ( X_i, X_j ) a_1 ( Y_i, Y_j )  +  \frac 14  \e [ \Vert X_1 - X_2 \Vert  ]  \sum_{ \substack {  1 \leq i, j \leq 4,  \\ i \neq j } }  a_1 ( Y_i, Y_j )  \nn \\
& \hspace{1cm} +  \frac 14  \e [ \Vert Y_1 - Y_2 \Vert  ]  \sum_{ \substack {  1 \leq i, j \leq 4,  \\ i \neq j  } }  a_1 (  X_i, X_j )  + 3  \e [ \Vert X_1 - X_2 \Vert  ]  \e [ \Vert Y_1 - Y_2 \Vert  ],
\end{align*}
\begin{align*}
& \frac 1 4 \sum_{i = 1}^4  \bigg(  \sum_{ \substack {  1 \leq  j \leq 4, \\  j \neq i } }    \Vert X_i - X_j \Vert   \sum_{ \substack {  1 \leq  j \leq 4, \\  j \neq i  } } \Vert Y_i - Y_j \Vert   \bigg)  \nn \\
& =   \frac 1 4 \sum_{i = 1}^4  \bigg(  \sum_{ \substack {   1 \leq  j \leq 4, \\  j \neq i  } }    a_1 ( X_i, X_j )   \sum_{ \substack { 1 \leq  j \leq 4, \\  j \neq i  } }  a_1 ( Y_i, Y_j )   \bigg)  + \frac 3 4  \e [ \Vert X_1 - X_2 \Vert  ]    \sum_{ \substack { 1 \leq i, j \leq 4,  \\ i \neq j } }  a_1 ( Y_i, Y_j )  \nn \\
&  \hspace{1cm}   + \frac 3 4  \e [ \Vert Y_1 - Y_2 \Vert  ]    \sum_{ \substack { 1 \leq i, j \leq 4,  \\ i \neq j } }  a_1 ( X_i, X_j )   + 9  \e [ \Vert X_1 - X_2 \Vert  ]  \e [ \Vert Y_1 - Y_2 \Vert  ],    
\end{align*}
and
\begin{align*}
& \frac {1} {24} \sum_{ \substack {  1 \leq i, j \leq 4,  \\ i \neq j } } \Vert X_i - X_j \Vert   \sum_{ \substack {  1 \leq i, j \leq 4,  \\ i \neq j  } }   \Vert Y_i - Y_j \Vert  \\
&   =  \frac {1} {24} \sum_{ \substack { 1 \leq i, j \leq 4,  \\ i \neq j } }  a_1 ( X_i, X_j )   \sum_{ \substack { 1 \leq i, j \leq 4,  \\ i \neq j } }  a_1 ( Y_i, Y_j ) + \frac  1 2  \e [ \Vert X_1 - X_2 \Vert  ]  \sum_{ \substack { 1 \leq i, j \leq 4,  \\ i \neq j } }  a_1 ( Y_i, Y_j ) \nn \\
& \hspace{1cm}  + \frac  1 2  \e [ \Vert X_1 - X_2 \Vert  ]  \sum_{ \substack { 1 \leq i, j \leq 4,  \\ i \neq j } }  a_1 ( Y_i, Y_j ) +  6  \e [ \Vert X_1 - X_2 \Vert  ]  \e [ \Vert Y_1 - Y_2 \Vert  ].
\end{align*} 
By these equalities and \eqref{kernel}, we can obtain
\begin{align*}
&  h ( ( X_1, Y_1 ) , ( X_2, Y_2 ),  ( X_3, Y_3 ), ( X_4, Y_4 ) ) \\
& =   \frac 14 \sum_{ \substack { 1 \leq i, j \leq 4,  \\ i \neq j } }  a_1 ( X_i, X_j ) a_1 ( Y_i, Y_j )  - \frac 1 4 \sum_{i = 1}^4  \bigg(  \sum_{ \substack {  1 \leq  j \leq 4, \\ j \neq i } }    a_1 ( X_i, X_j )   \sum_{ \substack { 1 \leq  j \leq 4, \\ j \neq i  } }  a_1 ( Y_i, Y_j )   \bigg)  \nn \\
& \quad + \frac {1} {24} \sum_{ \substack {  1 \leq i, j \leq 4,  \\ i \neq j  } }  a_1 ( X_i, X_j )   \sum_{ \substack {  1 \leq i, j \leq 4,  \\ i \neq j  } }  a_1 ( Y_i, Y_j ) .  
\end{align*}

\textit{Step 2}. Since $ d ( X_1, X_2 )  = a_1 ( X_1, X_2 ) - a_2 ( X_1 ) - a_2 ( X_2 ) $, it holds that 
\begin{align*}
&  \sum_{ \substack {  1 \leq i, j \leq 4,  \\ i \neq j  } }  a_1 ( X_i, X_j ) a_1 ( Y_i, Y_j ) \nn \\
& =   \sum_{ \substack {  1 \leq i, j \leq 4,  \\ i \neq j  } }  d ( X_i, X_j ) d ( Y_i, Y_j )  + 2 \sum_{ \substack {  1 \leq i, j \leq 4,  \\ i \neq j  } }  d ( X_i, X_j ) a_3 (Y_i) + 2  \sum_{ \substack {  1 \leq i, j \leq 4,  \\ i \neq j  } }  d ( Y_i, Y_j ) a_2 (X_i)  \nn \\
& \qquad + 4   \sum_{i = 1}^4 a_2 ( X_i ) a_3 ( Y_i ) + 2   \Big( \sum_{ i = 1}^4 a_2 ( X_i )  \Big)  \Big( \sum_{ i = 1}^4 a_3 ( Y_i )  \Big) ,
\end{align*}
\begin{align*}
& \sum_{i = 1}^4  \bigg(  \sum_{ \substack {    1 \leq  j \leq 4, \\  j \neq i  } }    a_1 ( X_i, X_j )   \sum_{ \substack { 1 \leq  j \leq 4, \\  j \neq i   } }  a_1 ( Y_i, Y_j )   \bigg)  \\
& =  \sum_{i = 1}^4  \bigg(  \sum_{ \substack { 1 \leq  j \leq 4, \\  j \neq i   } }     d ( X_i, X_j )   \sum_{ \substack { 1 \leq  j \leq 4, \\  j \neq i   } }  d ( Y_i, Y_j )   \bigg)  +  \Big(   \sum_{ j = 1 }^4 a_3 ( Y_j )  \Big) \Big(   \sum_{ \substack {  1 \leq i, j \leq 4,  \\ i \neq j  } }   d ( X_i , X_j ) \Big) \nn \\
& \qquad +  \Big(   \sum_{ j = 1 }^4 a_2 ( X_j )  \Big) \Big(   \sum_{ \substack { 1 \leq i, j \leq 4,  \\ i \neq j  } }   d ( Y_i , Y_j )  \Big)  +  2  \sum_{ \substack {  1 \leq i, j \leq 4,  \\ i \neq j  } }    d ( X_i , X_j ) a_3 ( Y_i )    \nn \\
& \qquad +  2  \sum_{ \substack { 1 \leq i, j \leq 4,  \\ i \neq j  } }    d ( Y_i , Y_j ) a_2 ( X_i )  +  8 \Big( \sum_{ i = 1 }^4 a_2 ( X_i )  \Big)   \Big( \sum_{ i = 1 }^4 a_3 ( Y_i )  \Big) + 4 \sum_{i = 1}^4 a_2 ( X_i ) a_3 ( Y_i ),
\end{align*}
and 
\begin{align*}
& \sum_{ \substack { 1 \leq i, j \leq 4,  \\ i \neq j } }  a_1 ( X_i, X_j )   \sum_{ \substack { 1 \leq i, j \leq 4,  \\ i \neq j } }  a_1 ( Y_i, Y_j )  \nn \\
& = \sum_{ \substack { 1 \leq i, j \leq 4,  \\ i \neq j } }  d ( X_i, X_j )   \sum_{ \substack { 1 \leq i, j \leq 4,  \\ i \neq j } }   d  ( Y_i, Y_j )  +  6 \Big( \sum_{ i = 1 }^4  a_2 ( X_i ) \Big) \Big( \sum_{ \substack { 1 \leq i, j \leq 4,  \\ i \neq j } }  d ( Y_i, Y_j )  \Big) \nn \\
& \qquad  +  6 \Big( \sum_{ i = 1 }^4  a_3 ( Y_i ) \Big) \Big( \sum_{ \substack { 1 \leq i, j \leq 4,  \\ i \neq j } }  d ( X_i, X_j )  \Big) + 36  \Big( \sum_{ i = 1 }^4 a_2 ( X_i )  \Big) \Big( \sum_{ i = 1 }^4 a_3 ( Y_i )  \Big).
\end{align*}
Combining the above three equalities yields \eqref{ex-h-d}. This concludes the proof of Lemma \ref{h-d}.

\subsection{Lemma \ref{le-h2} and its proof} \label{SecC.6}

\begin{lemma}  \label{le-h2}
	If $ X $ is independent of $ Y $, it holds that 
	\begin{equation*}
	\e [ h ( (X_1, Y_1), (X_2, Y_2), (X_3, Y_3), (X_4, Y_4) ) ]^2 = \frac 12 \V^2 ( X ) \V^2 ( Y ).  
	\end{equation*}
\end{lemma}

\noindent \textit{Proof}. From Lemma \ref{h-d}, we can deduce
\begin{equation*}
\e [ h ( (X_1, Y_1), (X_2, Y_2), (X_3, Y_3), (X_4, Y_4) ) ]^2 =  \sum_{ k = 1}^6 I_k,
\end{equation*} 
where
\begin{align*}
I_1 & =  \frac {1} { 16 } \e \Big[\Big(  \sum_{ \substack { 1 \leq i, j \leq 4,  \\ i \neq j } }  d ( X_i, X_j ) d ( Y_i, Y_j )   \Big)^2\Big],  \nn \\
I_2 & =   \frac {1} { 16 } \e \bigg\{\bigg[  \sum_{i = 1}^4  \bigg(  \sum_{ \substack {    1 \leq  j \leq 4, \\ j \neq i  } }    d ( X_i, X_j )   \sum_{ \substack {   1 \leq  j \leq 4, \\ j \neq i  } }  d ( Y_i, Y_j )   \bigg)   \bigg]^2\bigg\},  \nn \\
I_3 & = \frac { 1 } { 576 } \e \Big[\Big(   \sum_{ \substack { 1 \leq i, j \leq 4,  \\ i \neq j } }  d ( X_i, X_j )   \sum_{ \substack { 1 \leq i, j \leq 4,  \\ i \neq j } }  d ( Y_i, Y_j ) \Big)^2\Big], \nn \\
I_4 & =  - \frac 1 8  \e \bigg\{  \Big( \sum_{ \substack { 1 \leq i, j \leq 4,  \\ i \neq j } }  d ( X_i, X_j ) d ( Y_i, Y_j )    \Big)  \bigg[  \sum_{i = 1}^4  \Big(  \sum_{ \substack { 1 \leq  j \leq 4, \\ j \neq i } }    d ( X_i, X_j )   \sum_{ \substack { 1 \leq  j \leq 4, \\ j \neq i } }  d ( Y_i, Y_j )   \Big)  \bigg] \bigg\}, \nn \\
I_5 & = \frac { 1 } { 48  } \e \bigg[ \Big( \sum_{ \substack { 1 \leq i, j \leq 4,  \\ i \neq j } }  d ( X_i, X_j ) d ( Y_i, Y_j )   \Big)  \Big(  \sum_{ \substack { 1 \leq i, j \leq 4,  \\ i \neq j } }  d ( X_i, X_j )  \Big) \Big(    \sum_{ \substack { 1 \leq i, j \leq 4,  \\ i \neq j } }  d ( Y_i, Y_j )  \Big) \bigg],   
\end{align*}
and
\begin{align*}
I_6 & = - \frac { 1 }   { 48 }  \e \bigg\{   \bigg[  \sum_{i = 1}^4  \Big(  \sum_{ \substack { 1 \leq  j \leq 4, \\ j \neq i } }    d ( X_i, X_j )   \sum_{ \substack { 1 \leq  j \leq 4, \\ j \neq i } }  d ( Y_i, Y_j )   \Big)  \bigg]  \\
& \hspace{5cm}  \times  \Big(  \sum_{ \substack { 1 \leq i, j \leq 4,  \\ i \neq j  } }  d ( X_i, X_j )  \Big) \Big(    \sum_{ \substack { 1 \leq i, j \leq 4,  \\ i \neq j } }  d ( Y_i, Y_j )  \Big)  \bigg\}.
\end{align*} 

Since $ \e [ d ( X_1, X_2 ) d( X_1, X_3 ) ] = 0 $ and $ \e [ d ( X_1, X_2 ) ] = 0 $, under the independence of $ X $ and $ Y $ we have
\begin{align*}
I_1 & =  \frac 32   \e [ d^2 ( X_1, X_2 ) ] \e [ d^2 ( Y_1, Y_2 ) ]
\end{align*}
and
\begin{align*}
I_2 & =  \frac { 1 } { 16 } \sum_{i = 1}^4  \bigg[ \e \Big(   \sum_{ \substack {  1 \leq  j \leq 4 ,  \\  j \neq i    } } d ( X_i, X_j ) \Big)^2  \e \Big( \sum_{ \substack {   1 \leq  j \leq 4 ,  \\  j \neq i    } }  d ( Y_i, Y_j )    \Big)^2  \bigg] \nn \\
& \quad + \frac { 1 } { 16 }   \sum_{ \substack { 1 \leq   i , k \leq 4, \\  i \neq k } }  \bigg\{  \e \bigg[ \bigg( \sum_{ \substack { 1 \leq   j  \leq 4,  \\  j \neq i  } }  d ( X_i, X_j )  \bigg)   \bigg( \sum_{ \substack { 1 \leq   l  \leq 4,  \\  l \neq k  } }  d ( X_k, X_l )  \bigg)  \bigg] \\
& \hspace{5cm} \times   \e \bigg[ \bigg( \sum_{ \substack { 1 \leq   j  \leq 4,  \\  j \neq i  } }  d ( Y_i, Y_j )  \bigg)   \bigg( \sum_{ \substack { 1 \leq   l  \leq 4,  \\  l \neq k  } }  d ( Y_k, Y_l )  \bigg)  \bigg]  \bigg\} \nn \\
& =  \frac {1} { 16 } \times 4 \times 3 \e [ d^2 ( X_1, X_2 ) ] \times 3 \e [ d^2 ( Y_1, Y_2 ) ] + \frac {1} { 16 }    \sum_{ \substack { 1 \leq   i , k \leq 4, \\  i \neq k } } \e [ d^2 ( X_i, X_k )  ] \e [ d^2 ( Y_i, Y_k )  ]  \nn \\
& = 3 \e [ d^2 ( X_1, X_2 ) ] \e [ d^2 ( Y_1, Y_2 ) ].
\end{align*}

Similarly, we can obtain
\begin{align*}
I_3  & = \frac { 16 } { 576 }  \e  \Big(   \sum_{ 1 \leq  i < j \leq 4  }  d ( X_i, X_j )   \Big)^2  \e  \Big(   \sum_{ 1 \leq  i < j \leq 4  }  d ( Y_i, Y_j )   \Big)^2  \nn \\
& =    \e [ d^2 ( X_1, X_2 ) ] \e [ d^2 ( Y_1, Y_2 ) ], 
\end{align*}
\begin{align*}
I_4 & =  - \frac 1 4 \sum_{ i = 1 }^4  \e \bigg\{  \Big( \sum_{ 1 \leq k < l \leq 4 }  d ( X_k, X_l ) d ( Y_k, Y_l )    \Big)   \Big(  \sum_{ \substack { 1 \leq  j \leq 4, \\  j \neq i } }    d ( X_i, X_j )  \Big)  \Big(  \sum_{ \substack { 1 \leq  j \leq 4, \\  j \neq i } }    d ( Y_i, Y_j )  \Big)   \bigg\} \nn \\
& =  - \frac 1 4 \times 4 \times 3 \e [ d^2 ( X_1, X_2 ) ] \e [ d^2 ( Y_1, Y_2 ) ] = - 3 \e [ d^2 ( X_1, X_2 ) ] \e [ d^2 ( Y_1, Y_2 ) ],
\end{align*}
\begin{align*}
I_5  & =  \frac { 8 } { 48 } \sum_{ 1 \leq i < j \leq 4 }  \e \bigg[ d ( X_i, X_j ) d ( Y_i, Y_j )  \Big(  \sum_{ 1 \leq k < l \leq 4 }  d ( X_k, X_l )  \Big) \Big(   \sum_{ 1 \leq k < l \leq 4 }  d ( Y_k, Y_l )  \Big)  \bigg] \nn \\
& =  \e [ d^2 ( X_1, X_2 ) ] \e [ d^2 ( Y_1, Y_2 ) ] ,
\end{align*}
and
\begin{align*}
I_6 & =  - \frac { 4 } { 48 } \sum_{ i = 1 }^4  \bigg\{  \e \bigg[  \Big(   \sum_{ \substack {  1 \leq  j \leq 4, \\ j \neq i } }    d ( X_i, X_j )    \Big)   \Big(  \sum_{ 1 \leq k < l \leq 4 } d ( X_k,  X_l )  \Big)   \bigg]  \\
& \hspace{5cm} \times  \e \bigg[  \Big(   \sum_{ \substack {  1 \leq  j \leq 4, \\ j \neq i } }    d ( Y_i, Y_j )    \Big)   \Big(  \sum_{ 1 \leq k < l \leq 4 } d ( Y_k, Y_l )  \Big)   \bigg] \bigg\} \nn \\
& =  - \frac { 4 \times 4 \times 9 } { 48  }  \e [ d^2 ( X_1, X_2 ) ] \e [ d^2 ( Y_1, Y_2 ) ]  = - 3  \e [ d^2 ( X_1, X_2 ) ] \e [ d^2 ( Y_1, Y_2 ) ] .
\end{align*}
Consequently, it follows that 
\begin{align*}
\e [ h ( (X_1, Y_1), (X_2, Y_2), (X_3, Y_3), (X_4, Y_4) ) ]^2 & =  \sum_{ k = 1}^6 I_k = \frac 12  \e [ d^2 ( X_1, X_2 ) ] \e [ d^2 ( Y_1, Y_2 ) ]  \\
& = \frac { 1 }  { 2 } \V^2 ( X )  \V^2 ( Y ),
\end{align*}
which completes the proof of Lemma \ref{le-h2}.

\subsection{Lemma \ref{le-cru} and its proof} \label{SecC.7}
The following lemma provides some basic inequalities that are based on the Taylor expansion and serve as the fundamental ingredients for the proofs of Propositions \ref{prop1}--\ref{prop3}.

\begin{lemma}    \label{le-cru}
	For $ x \geq  -1 $, it holds that 
	\begin{align}
	\big| ( 1 + x )^{ 1/ 2 } - 1   \big|  & \leq  | x |,  \label{cru1} \\
	\big| ( 1 + x )^{ 1/ 2 } - ( 1 +   x / 2   )    \big|  & \leq     x^2 /2  ,  \label{cru2}   \\
	\big| ( 1 + x )^{ 1 /2 } -  ( 1 +    x/ 2  -   x^2 / 8   )  \big| & \leq  3 |  x^3 | / 8 ,  \label{cru3} \\
	\big|  (1 + x)^{1/2}  - (1 + x/ 2 - x^2 / 8  +  x^3/16)  \big| & \leq x^4. \label{cru4}
	\end{align}
\end{lemma}

\noindent \textit{Proof}. (i) We first prove \eqref{cru1}. 
It is evident that $  1 + x \leq ( 1 + x )^{1/2}  \leq 1  $ for $ x \in [ - 1, 0 ] $ and $ 1 < ( 1 + x )^{1/2}   < 1 + x   $ for $ x \in ( 0, \infty) $. Thus we can obtain \eqref{cru1} directly.

(ii) We next show \eqref{cru2}. 
Define $ u_1 (x) = ( 1 + x )^{1/2} - ( 1 + x /2  ) $. Then we have the derivative 
\begin{align*}
u_1' ( x ) =    [ ( 1 + x )^{ - 1 / 2 } - 1 ] / 2, 
\end{align*}
$ u_1' ( x ) >  0 $ for $ x \in [ -1, 0 )  $, and $ u_1' ( x ) < 0  $ for $ x \in ( 0 , \infty ) $. Since $ u_1 (0) = 0 $, it holds that $ u_1 (x) \leq 0 $ for $ x \geq - 1  $. It remains to show that for $ x  \geq -1 $,
\begin{align*}
u_1 ( x )  \geq - x^2 /2. 
\end{align*}
Denote by $ u_2 (x)  =  ( 1 + x )^{ 1/2 } - ( 1 +  x /2  ) + x^2 /2  $. Then we have 
\begin{align*}
u_2' (x) & =  \frac 1 2 ( 1 + x )^{ 1/2 } - \frac 1 2 + x , \\
u_2'' (x) &  = - \frac 1 4 ( 1 + x  )^{ - 3/2 } + 1,
\end{align*} 
$ u_2 '' (x) \leq 0 $ for $ -1 \leq x \leq  4^{ - 2/3 } - 1 $, and $ u_2 '' (x) > 0  $ for $ x > 4^{  - 2/3 } - 1 $. In addition, it holds that $ u_2 ' (- 1 ) = + \infty  $, $ u_2 ' (0 ) = 0  $, and $ u_2' ( + \infty) = + \infty $, which lead to $ u_2 (x) \geq \min \{  u_2 ( - 1 ) , u_2 (0)  \} = 0  $. Hence the proof of \eqref{cru2} is completed.

(iii) We now prove \eqref{cru3}.  
First, the result is trivial when $ x = 0 $. Define $ u_3 (x) = ( 1 + x )^{ 1 /2 } -  ( 1 +    x/ 2  -   x^2 / 8   ) $. Then we have
\begin{align*}
u_3' (x) & =  \frac 12 [  ( 1 + x )^{ - 1 / 2 }  - 1 +   x/ 2] \nn \\ 
& = \frac 1 2 ( 1 + x )^{ - 1 } \big[ ( 1 + x )^{1/2} - ( 1 + x /2 - x^2 /2  )  \big].
\end{align*}
It has been shown in the proof above that $ ( 1 + x )^{ 1/2 } - ( 1 + x/2 ) + x^2 /2 \geq 0  $ for $ x \geq - 1 $. Thus $ u_3 ' (x) \geq 0  $ for $ x \geq - 1 $. It follows that $ u_3 (x) \leq 0  $ for $ -1 \leq x \leq 0 $ and $ u_3 (x) > 0  $ for $ x >  0   $.
Now it remains to show that for $ x \in [ -1, 0 ) \cup (0, \infty) $, 
\begin{align*}
u_3 (x) / x^3 \leq  3/ 8.
\end{align*}
It is easy to show that  
\begin{align*}
\Big( \frac { u_3 (x) } { x^3 } \Big)' = \frac { u_4 (x)  } { 2 x^4  },
\end{align*}
where $ u_4 (x) = -  5 ( 1 + x )^{1/2} -  ( 1 + x )^{ - 1/2} + 2 x  -  x^2/4 + 6 $. 

Observe that 
\begin{align*}
u_4' (x) & =  - \frac 5 2  ( 1 + x )^{ - 1 /2  } + \frac 1 2 ( 1 + x  )^{ - 3/ 2 } + 2 - \frac x 2  , \\
u_4'' (x) & =  \frac 5 4 ( 1 + x  )^{ - 3 / 2 } - \frac 3 4 ( 1 + x  )^{ - 5 /2  } - \frac 1 2  , \\
u_4''' (x) &  = - \frac { 15 } { 8 }  x ( 1 + x  )^{ - 7 / 2 },
\end{align*}
$ u_4 ''' ( x )  > 0 $ for $ x \in [-1, 0) $, and $  u_4''' ( x ) < 0 $ for $ x > 0$. Furthermore, $ u_4 '' ( 0 )  = 0 $ and thus $ u_4'' (x) \leq  0  $ for any $ x \in [-1, 0) \cup (0, \infty)  $. In addition, $ u_4' ( 0 ) = 0 $ and thus $ u_4' (x) > 0  $ for $ x \in [ - 1, 0 ) $ and $ u_4' (x) <  0  $ for $ x \in ( 0 , \infty ) $. Since $ u_4 ( 0 )  = 0 $, it follows that $ u_4 ( x ) <  0 $ for any $ x \in [-1, 0) \cup (0, \infty)  $, which entails that 
\begin{equation*}
\frac { u_3 (x) } { x^3 }  \leq   \frac  { u_3 (x) }  { x^3 } \Big\vert_{x = -1 } =   \frac 3 8 .
\end{equation*}
Similay by taking derivatives, \eqref{cru4} can be proved. We omit its proof to avoid redundancy. This concludes the proof of Lemma \ref{le-cru}.

\subsection{Lemma \ref{le-v2} and its proof} \label{SecC.8}

\begin{lemma} \label{le-v2}
	If $ \e [  \Vert X \Vert^2  ] < \infty $, then we have 
	\begin{align}
	\e [ W_{12}^2 ] & = B_X^{ - 2 } \big( 2 [ \e  \Vert  X  \Vert^4   - ( \e   \Vert  X  \Vert^2 )^2 ] + 4 \e [ ( X_1^T X_2 )^2 ] \big),  \label{W12} \\
	\e [ W_{12} W_{13} ] & =   B_X^{ - 2 }  [ \e  \Vert  X  \Vert^4   - ( \e   \Vert  X  \Vert^2 )^2 ] .  \label{W123}
	\end{align}
\end{lemma}

\noindent \textit{Proof}. Define $ \alpha_1 (X) =  \Vert X  \Vert^2  - \e [ \Vert X \Vert^2  ] $ and $ \alpha_2 ( X_1, X_2 )  = X_1^T X_2 $. By the definition of $ W_{12} $ and $ W_{13} $, we have
\begin{align*}
\e [ W_{12}^2 ]  & =  B_X^{ - 2 } \e \big \{ \big[ \alpha_1 ( X_1 ) + \alpha_1 ( X_2 ) - 2 \alpha_2 ( X_1, X_2 )  \big]^2  \big\}  , \\
\e [ W_{12} W_{13} ] & = B_X^{ - 2 } \e \big \{ \big[ \alpha_1 ( X_1 ) + \alpha_1 ( X_2 ) - 2 \alpha_2 ( X_1, X_2 )  \big]  \big[ \alpha_1 ( X_1 ) + \alpha_1 ( X_3 ) - 2 \alpha_2 ( X_1, X_3 )  \big] \big\}. 
\end{align*}
Since $ \e [ \alpha_1 (X) ] = 0 $ and $ \e (X)  = 0  $, by expanding the products above we can deduce 
\begin{align}
\e [ W_{12}^2  ] &  = B_X^{ - 2 } \big( 2 \e [ \alpha_1^2 ( X_1 ) ]  + 4 \e [ \alpha_2^2 ( X_1, X_2 ) ] \big) \nn \\
&  = B_X^{ - 2 }  \big( 2 [ \e  \Vert  X  \Vert^4   - ( \e   \Vert  X  \Vert^2 )^2 ] + 4 \e [ ( X_1^T X_2 )^2 ] \big)  \label{W12_sq} \\
\intertext{and}
\e [ W_{12} W_{13} ] & = B_X^{ - 2 } \e [ \alpha_1^2 ( X_1 ) ] =   B_X^{ - 2 }  [ \e  \Vert  X  \Vert^4   - ( \e   \Vert  X  \Vert^2 )^2 ].  \label{W12W13}
\end{align}
The desired result then follows immediately. This completes the proof of Lemma \ref{le-v2}.

\subsection{Lemma \ref{le-g-mo1} and its proof} \label{SecC.9}

\begin{lemma}  \label{le-g-mo1}
	If $ \e [ \Vert X \Vert ^4  ] < \infty $, then we have 
	\begin{align}
	& \e [ W_{12} W_{13} W_{24} W_{34} ] - 4 \e [ W_{12} W_{13} W_{24} W_{45}  ] + 2 ( \e [ W_{12} W_{13} ] )^2  \nn \\
	& =  16 B_X^{ - 4 } \e [ ( X_1^T \Sigma_x X_2 )^2 ].
	\end{align}
\end{lemma}

\noindent \textit{Proof}.  By the definition of $ W_{i j } $, we have 
\begin{align*}
&  \e [ W_{12} W_{13} W_{24} W_{34} ]  \nn \\
& =   B_X^{ - 4}  \e \big \{ \big[ \alpha_1 ( X_1 ) + \alpha_1 ( X_2 ) - 2 \alpha_2 ( X_1, X_2 )  \big]  \big[ \alpha_1 ( X_1 ) + \alpha_1 ( X_3 ) - 2 \alpha_2 ( X_1, X_3 )  \big]    \nn \\
& \hspace{1.5cm}  \times  \big[ \alpha_1 ( X_2 ) + \alpha_1 ( X_4 ) - 2 \alpha_2 ( X_2, X_4 )  \big]  \big[ \alpha_1 ( X_3 ) + \alpha_1 ( X_4 ) - 2 \alpha_2 ( X_3, X_4 )  \big]  \big\}.
\end{align*} 
Noting that $ \e [ \alpha_1 (X_1) ]  = 0 $ and $ \e [ X ] = 0  $, it follows from expanding the above product and the symmetry of $ X_1, \cdots, X_4 $ that
\begin{align*}
\e [  W_{12} W_{13} W_{24} W_{34}   ] & = B_X^{ - 4 } \big\{ 2  ( \e [ \alpha_1^2 ( X ) ]  )^2 + 16 \e [ \alpha_2 (  X_1, X_2 ) \alpha_2 ( X_1, X_3 ) \alpha_1 ( X_2 ) \alpha_1 ( X_3 ) ] \nn \\
& \hspace{1cm} +  16 \e [ \alpha_2 ( X_1, X_2 )  \alpha_2 ( X_1, X_3 )  \alpha_2 ( X_2, X_4 )  \alpha_2 ( X_3, X_4 )  ] \big\}.
\end{align*}

By the same token, we can deduce 
\begin{align*}
&  \e [ W_{12} W_{13} W_{24} W_{45} ]  \nn \\
& =   B_X^{ - 4}  \e \big \{ \big[ \alpha_1 ( X_1 ) + \alpha_1 ( X_2 ) - 2 \alpha_2 ( X_1, X_2 )  \big]  \big[ \alpha_1 ( X_1 ) + \alpha_1 ( X_3 ) - 2 \alpha_2 ( X_1, X_3 )  \big]    \nn \\
& \hspace{1.5cm}  \times  \big[ \alpha_1 ( X_2 ) + \alpha_1 ( X_4 ) - 2 \alpha_2 ( X_2, X_4 )  \big]  \big[ \alpha_1 ( X_4 ) + \alpha_1 ( X_5 ) - 2 \alpha_2 ( X_4, X_5 )  \big]  \big\}, \nn \\
& =  B_X^{ - 4 } \big\{ ( \e [ \alpha_1^2 ( X ) ]  )^2  + 4 \e [ \alpha_2 (  X_1, X_2 ) \alpha_2 ( X_1, X_3 ) \alpha_1 ( X_2 ) \alpha_1 ( X_3 ) ] \big\}.  
\end{align*}
Therefore, combining the above expressions with \eqref{W12W13} results in 
\begin{align*}
&  \e [ W_{12} W_{13} W_{24} W_{34} ]  - 4 \e [ W_{12} W_{13} W_{24} W_{45} ] + 2 (   \e [ W_{12} W_{13} ]  )^2 \nn \\
& =  16 B_X^{ - 4 } \e [ \alpha_2 ( X_1, X_2 )  \alpha_2 ( X_1, X_3 )  \alpha_2 ( X_2, X_4 )  \alpha_2 ( X_3, X_4 )  ]  =  16 B_X^{ - 4 } \e [ ( X_1^T \Sigma_x X_2 )^2 ],
\end{align*}
which concludes the proof of Lemma \ref{le-g-mo1}.

\subsection{Lemma \ref{le-dcov} and its proof} \label{SecC.10}

\begin{lemma} \label{le-dcov}
	For any random vectors $ X \in \R^ p$ and  $ Y \in \R^q $ satisfying $ \e [ \Vert X \Vert^{12} ] + \e [ \Vert Y \Vert^{12} ] < \infty$, we have
	\begin{align*}
	\V^2 (X, Y) = I_1 + I_2 + I_3 + I_4 + I_5,
	\end{align*}
	where 
	\begin{align*}
	I_1 &=  \frac 1 4 B_X^{1/2} B_Y^{1/2} \big( \e [ W_{12} V_{12} ] -  2 \e [ W_{12} V_{13} ]  \big), \\
	I_2 &= - \frac {1} { 16  } B_X^{1/2} B_Y^{1/2} \big( \e [ W_{12} V_{12}^2  ] - 2 \e [ W_{12}  V_{13}^2 ]  +  \e [ W_{12}^2 V_{12} ] - 2 \e [ W_{12}^2 V_{13} ] \big), \\
	I_3 &=  \frac {1} {32} B_X^{1/2} B_Y^{1/2} \big(  \e [ W_{12} V_{12}^3 ] -  2\e [ W_{12} V_{13}^3 ]  +  \e [ V_{12 } W_{12}^3 ] - 2 \e [ V_{12} W_{13}^3  ]  \big),  \\
	I_4 & =  \frac {1} { 64 } B_X^{1/2} B_Y^{1/2} \big( \e [W_{12}^2 V_{12}^2] - 2 \e [ W_{12}^2 V_{13}^2 ]   +  \e [ W_{12}^2 ] \e [ V_{12}^2 ]  \big), \\
	I_5 & = O \Big\{ B_X^{1/2} B_Y^{1/2} \Big[  ( \e | W_{12} |^5 )^{2/5} ( \e | V_{12} |^5 )^{3/5}  + ( \e | W_{12} |^5 )^{3/5} ( \e | V_{12} |^5 )^{2/5} \nn \\
	& \quad +  ( \e | W_{12} |^5 )^{1/5} ( \e | V_{12} |^5 )^{4/5} + ( \e | W_{12} |^5 )^{4/5} ( \e | V_{12} |^5 )^{1/5}  \nn \\
	& \quad  + ( \e |W_{12}|^6 )^{1/2} ( \e | V_{12} |^6 )^{1/2} \Big] \Big\}.
	\end{align*}	    	      	       
\end{lemma}

\noindent \textit{Proof}. We will conduct the Taylor expansion to $ \V^2 (X, Y) = \e [ d(X_1, X_2) d (Y_1, Y_2) ] $. In light of \eqref{d}, some straightforward calculations lead to 
\begin{align*}
\V^2 (X, Y) = \e [ b ( X_1, X_2  ) b (Y_1, Y_2)  ] - 2 \e [ b  ( X_1, X_2  ) b (Y_1, Y_3) ]  +  \e [ b (X_1, X_2) ] \e [ b (Y_1, Y_2) ],
\end{align*}
where $ b( X_1, X_2 ) = \Vert X_1 - X_2 \Vert - B_X^{1/2} $ and $ b ( Y_1, Y_2  ) = \Vert Y_1 - Y_2 \Vert  - B_Y^{1/2}$. Define 
$$ 
W_{ij}=  B_X^{-1} ( \Vert X_i - X_j \Vert^2  - B_X ) \ \text{ and } \  V_{ij} = B_Y^{-1}  ( \Vert Y_i - Y_j \Vert^2 - B_Y ). $$
Observe that $  b (X_1, X_2)  =  B_X^{1/2} [ ( 1 + W_{12} )^{1/2} - 1 ]$. An application of similar arguments as those in the proof of Proposition \ref{prop2} by resorting to \eqref{cru4} in Lemma \ref{le-cru} yields 
\begin{align*}
\e & [ b(X_1, X_2) b(Y_1, Y_2) ] = B_X^{1/2} B_Y^{1/2} \Big\{  \frac  1 4  \e [ W_{12} V_{12} ]  - \frac {1} { 16 } \big( \e [ W_{12} V_{12}^2  ] + \e [ W_{12}^2 V_{12} ] \big) \\
& \quad +  \frac {1} { 64 } \e [ W_{12}^2 V_{12}^2 ] + \frac {1} {32} \big( \e [W_{12}V_{12}^3] + \e [ W_{12}^3 V_{12} ] \big) + O \Big( \e [ | W_{12} |  V_{12}^4 ]   +  \e [ W_{12}^4 | V_{12} | ]  \\
& \quad +  \e [ W_{12}^2 | V_{12}^3 | ]  + \e [ | W_{12}^3 | V_{12}^2 ]  +  \e [ | W_{12} |^3 | V_{12} |^3 ] \Big)
\Big\}.
\end{align*}

By the same token, we can deduce that 
\begin{align*}
\e & [ b (X_1 , X_2) b (Y_1, Y_3) ] = B_X^{1/2} B_Y^{1/2} \Big\{   \frac  1 4 \e [ W_{12} V_{13} ]  - \frac {1} {16} \big( \e [ W_{12} V_{13}^2 ] + \e [ W_{12}^2 V_{13} ] \big)   \\
& \quad + \frac {1} {64} \e [ W_{12}^2 V_{13}^2 ] + \frac {1} {32} \big(  \e [ W_{12} V_{13}^3 ]  + \e [ W_{12}^3 V_{13} ] \big)  + O \Big( \e [ | W_{12} | V_{13}^4 ] +  \e [ W_{12}^4 V_{13} ] \\
& \quad + \e [ W_{12}^2 | V_{13} |^3 ] + \e[ | W_{12} |^3 V_{13}^2 ] +  \e [ | W_{12} |^3 | V_{13} |^3 ] \Big)
\Big\} 
\end{align*}	
and 
\begin{align*}
\e [ b (X_1, X_2) ] \e[ b (Y_1, Y_2) ] & = B_X^{1/2} B_Y^{1/2} \Big\{   \frac {1} {64}  \e[W_{12}^2] \e [ V_{12}^2 ] + O \Big( \e [ W_{12}^2 ] \e [  | V_{12} | ^3]  \\
& \quad + \e   [ | W_{12}^3 | ]  \e [ V_{12}^2 ] \Big)
\Big\}.
\end{align*}
Therefore, the desired decomposition follows from a combination of the above three representations and the Cauchy--Schwarz inequality, which completes the proof of Lemma \ref{le-dcov}.

\renewcommand{\thesubsection}{D.\arabic{subsection}}
\section{Theoretical results for the case of $ 1/2 < \tau \leq 1 $} \label{tau} 

\subsection{Theory} \label{SecD.1}
In this section, we introduce our parallel results of Theorems \ref{thm2} and \ref{thm4} for the case of $ 1/2 < \tau  \leq 1 $. When $ \e [ \Vert X \Vert^{ 2 + 2 \tau  }  ] + \e [ \Vert Y \Vert ^{ 2 + 2 \tau }  ]  < \infty $ for a larger value of $\tau$ with $ 1/ 2 < \tau \leq 1 $, the key ingredient is that higher-order Taylor expansions can be applied while bounding $ \e [ g( X_1, X_2, X_3, X_4 ) ] $. We start with presenting the expansion of $ \e [ g( X_1, X_2, X_3, X_4 ) ] $ for $ 1/ 2 < \tau \leq 1 $. Let us define 
\begin{gather*}
\mathscr{G}_1 (X)    =     \big| \e [ ( X_1^T X_2 )^2  X_1^T \Sigma_x^2  X_2  ] \big|  ,  \quad
\mathscr{G}_2 ( X )  =   \e [  \Vert X_1 \Vert^2 ( X_1^T \Sigma_x X_2 )^2   ],   \\   
\mathscr{G}_3 ( X )   = \e [ X^T X X^T  ] \Sigma_x^2 \e [ X X^T X ] ,             \\ 
N_{\tau} (X)  =     \frac { \e [ ( X_1^T \Sigma_x X_2 )^2 ]  + B_X^{ - 2 \tau  } L_{x, \tau}^{ (2 + \tau) / ( 1 + \tau  ) } + B_X^{ - 1 }  \sum_{ i = 1 }^3 \mathscr{G}_i ( X ) } {  (  \e [ ( X_1^T X_2 )^2 ]  )^2 }.                  
\end{gather*}
We also have $ \mathscr{G}_1 (Y),  \mathscr{G}_2 (Y), \mathscr{G}_3 (Y)$, and $ N_{ \tau } ( Y ) $ that are defined in a similar way.

\begin{proposition}  \label{prop4}
	If $ \e [ \Vert X \Vert ^{ 4 + 4 \tau }   ]  < \infty $ for some $ 1/2 < \tau \leq  1 $, then there exists some absolute positive constant $ C   $ such that 
	\begin{align}
	& \big|  \e [ g ( X_1, X_2, X_3, X_4 ) ]  \big| \nn \\
	& \leq   C  \big\{   B_X^{ - 2 } \e [ ( X_1^T \Sigma_x X_2 )^2  ]  +   B_X^{ - 3 } \textstyle \sum_{i = 1}^3 \mathscr{G}_i (X)  + B_X^{ - ( 2 + 2 \tau ) }  L_{ x, \tau }^{  ( 2 + \tau )/ ( 1 + \tau )  }  \big\} .  \label{g2}
	\end{align}         
\end{proposition}

The proof of Proposition \ref{prop4} is given in Section \ref{SecD.2}. We can obtain the following central limit theorem and the associated rate of convergence for the case of $ 1/2 < \tau \leq 1 $ by substituting the bounds in Propositions \ref{prop1}--\ref{prop2} and \ref{prop4} into Theorem \ref{thm3}. 

\begin{theorem}   \label{thm8}
	Assume that $ \e  [ \Vert X \Vert ^{ 4 + 4 \tau }   ] + \e [  \Vert Y \Vert ^{ 4 + 4 \tau }   ] < \infty  $ for some $ 1/ 2 <  \tau \leq 1  $ and 
	\begin{align}
	B_X^{ - 1  } L_{x, 1/ 2} /  \e [ ( X_1 ^T X_2 )^2 ]  & \leq \frac {1} { 18 },    \label{cc1}\\
	B_Y^{  - 1 } L_{ y, 1/ 2 }  / \e [ ( X_1^T X_2 )^2 ] & \leq \frac {1} { 18 }.  \label{cc2}    
	\end{align}
	Then under the independence of $ X $ and $ Y $, we have 
	\begin{align}
	& \sup \limits_{ x \in \mathbb{R} }  | \mathbb{P}  ( T_n \leq x  ) - \Phi (x) |  \nn \\ 
	& \quad  \leq   C \bigg\{  [ N_{\tau} (X)  N_{\tau} (Y)  ]^{  ( 1 + \tau ) / 2 } +  \frac  { n^{ - \tau }  L_{ x, \tau }  L_{ y, \tau }  } { \big(  \e [ ( X_1^T X_2 )^2 ] \e [ ( Y_1^T Y_2 )^2 ]   \big)^{ 1 + \tau }  } \bigg\}^{  1 / ( 3 + 2 \tau )  }.  \label{rate-l}
	\end{align} 
\end{theorem}

The proof of Theorem \ref{thm8} is provided in Section \ref{SecD.2.new}. Theorem \ref{thm9} below is a direct corollary of Theorem \ref{thm8}.

\begin{theorem}   \label{thm9}
	Assume that $ \e [ \Vert X \Vert^{  4 + 4 \tau } ] + \e [ \Vert Y \Vert ^{ 4 + 4 \tau } ] < \infty $ for some $ 1/2 < \tau \leq 1  $ and \eqref{co1} holds as $ n \rightarrow  \infty$ and $ p + q \rightarrow \infty $. In addition, assume that \eqref{cc1} and $ N_{\tau} (X) \rightarrow 0 $ are satisfied as $ p \rightarrow \infty $, and that\eqref{cc2} and $ N_{\tau} (Y) \rightarrow 0 $ are satisfied as $ q \rightarrow \infty $. Then under the independence of $ X $ and $ Y $, we have 
	\begin{align*}
	T_n \stackrel{\mathscr{D}}{ \rightarrow } N(0, 1).
	\end{align*}
\end{theorem}

\subsection{Proof of Theorem \ref{thm8}} \label{SecD.2.new}
Note that \eqref{cc1}, \eqref{cc2}, and Proposition \ref{prop2} entail that 
\begin{align}
\V^2 (X)  \geq   B_X^{ - 1 } \e [ ( X_1^T X_2 )^2  ] /2 \ \text{ and } \    \V^2 (Y)  \geq   B_Y^{ - 1 } \e [ ( Y_1^T Y_2 )^2  ] /2, \label{V-lowb}
\end{align}
which together with \eqref{g2} leads to 
\begin{align*}
\frac  {    |  \e [ g (X_1,X_2, X_3, X_4 ) ]  \e [ g( Y_1, Y_2, Y_3, Y_4 ) ]  |  }  {  [  \V^2 ( X ) \V^2 ( Y )  ] ^2 }    \leq   N_{\tau} (X)  N_{\tau} (Y).
\end{align*} 
It follows from Proposition \ref{prop1} and \eqref{V-lowb} that
\begin{align*}
\frac {  \e [ | d (X_1, X_2 ) |^{2 + 2 \tau } ]  \e [  | d (Y_1, Y_2) |^{2 + 2 \tau } ] } { n^{\tau}  [ \V^2 ( X ) \V^2 ( Y ) ]^{ 1 + \tau  }    } \leq \frac  {  n^{ - \tau } L_{ x, \tau } L_{y, \tau} } {   ( \e [ ( X_1^T X_2 )^2 ]  \e [  (Y_1^T Y_2 )^2 ] )^{ 1 + \tau } }.
\end{align*}
Therefore, we can obtain the desired result \eqref{rate-l} by Theorem \ref{thm3}, which concludes the proof of Theorem \ref{thm8}.

\subsection{Proof of Proposition \ref{prop4} } \label{SecD.2}
It suffices to analyze the terms on the right hand side of \eqref{g-expan}. Compared to Proposition \ref{prop3}, we assume higher moments and thus we can conduct higher-order Taylor expansions for term $  ( 1 + W_{12} )^{1/2} $.

Let us first deal with term $ G_1 $. Denote by $ D_1 = \{  \max ( W_{12} , W_{ 13} , W_{24}, W_{34} ) \leq 1 \}  $ and $ D_1^{c} $ the complement of $ D_1 $. Following the notation in the proof of Proposition \ref{prop3}, by \eqref{cru1} and \eqref{cru3} we can deduce 
\begin{align*}
G_1 
& =  B_X^2 \e \Big(  \big[ \frac 1 2 W_{12} - \frac 1 8 W_{12}^2 + O(1) ( | W_{12} |^3 ) \big]  \big[ \frac 1 2 W_{13} - \frac 1 8 W_{13}^2 + O (1) ( | W_{13} |^3 ) \big]  \nn \\
& \hspace{2cm} \times \big[ \frac 1 2 W_{24} - \frac 1 8 W_{24}^2 + O(1) ( | W_{24} |^3 ) \big] \big[ \frac 1 2 W_{34} - \frac 1 8 W_{34}^2 + O(1) ( | W_{34} |^3 )  \big] \I \{ D_1 \} \Big)  \nn \\
& \qquad  +  O(1) B_X^2  \e [ | W_{12} W_{13} W_{24} W_{34} | \I \{ D_1^{c} \} ].  
\end{align*}
By expanding the products and reorganizing the terms, it holds that 
\begin{align*}
G_1 & = \frac { B_X^2 } { 16 } \Big(  \e [ W_{12} W_{13} W_{24} W_{34}  ] - \e [ W_{12}^2 W_{13} W_{24} W_{34} ]   + O(1) \e [ W_{12}^2 W_{13}^2 | W_{24}  W_{34} |  \I \{ D_1 \} ] \nn \\
& \hspace{2cm} +   O(1) \e [ W_{12}^2 W_{34}^2 | W_{13} W_{24} | \I \{ D_1 \} ] +  O(1) \e [ | W_{12} |^3 | W_{ 13 } W_{ 24  } W_{ 34 } |  \I \{ D_1 \} ]  \nn \\
& \hspace{2cm} +  O(1)  \e [ | W_{12} W_{13} W_{24} W_{34} | \I \{ D_1^{c} \} ]  \Big).
\end{align*}
Furthermore, if $ \e [ \Vert X \Vert^{ 4 + 4 \tau } ] < \infty $ for some $  1/2 < \tau \leq 1 $, then an application of Chebyshev's inequality and the Cauchy--Schwarz inequality results in 
\begin{align*}
& | \e [  W_{12}^2 W_{13}^2 | W_{24}  W_{34} |  \I ( D_1 ) ] | \nn \\
& \leq \e [ | W_{12}  |^{ 1 +  \tau }   | W_{13} |^{1 +  \tau } | W_{24} W_{34} | ]  \nn \\
& =  \e \big\{ \e [ | W_{12} |^{ 1 + \tau } | W_{13} |^{1 + \tau } \vert X_2, X_3 ] \e [ | W_{24} W_{34} | \vert X_2, X_3 ] \big\} \nn \\
&  = \e \big\{ ( \e [ | W_{12} |^{ 2  + 2 \tau  } \vert X_2 ]  ) ^{ \frac { 2 + \tau  } { 2 + 2 \tau  } } \big\} \e \big\{   ( \e [ | W_{13} |^{ 2  + 2 \tau  } \vert X_3 ]  ) ^{ \frac { 2 + \tau  } { 2 + 2 \tau  } }    \big\} \nn \\
& \leq  ( \e [ | W_{12} |^{ 2  + 2 \tau } ] )^{ \frac { 2 + \tau  } { 2 + 2 \tau  } }  ( \e [ | W_{13} |^{ 2  + 2 \tau } ] )^{ \frac { 2 + \tau  } { 2 + 2 \tau  } }   = ( \e [ | W_{12} |^{ 2  + 2 \tau } ] )^{ \frac { 2 + \tau  } { 1 +  \tau  } }. 
\end{align*}

By the same token, we can obtain
\begin{align*}
\e [ W_{12}^2 W_{34}^2 | W_{13} W_{24} | \I ( D_1 ) ]  & \leq \e [ | W_{12}  |^{ 1 +  \tau }   | W_{34} |^{1 +  \tau } | W_{13} W_{24} | ]  \nn \\
& \leq ( \e [ | W_{12} |^{ 2  + 2 \tau } ] )^{ \frac { 2 + \tau  } { 1 +  \tau  } }, \nn  \\
\e [ | W_{12} |^3 | W_{ 13 } W_{ 24  } W_{ 34 } |  \I ( D_1 ) ] &  \leq  \e [ | W_{12} |^{ 1 + 2 \tau  } | W_{13} W_{24} W_{34} | ] \nn \\
& \leq ( \e [ | W_{12} |^{ 2  + 2 \tau } ] )^{ \frac { 2 + \tau  } { 1 +  \tau  } },
\end{align*}
and 
\begin{align*}
\e [ | W_{12} W_{13} W_{24} W_{34} | \I ( D_1^c ) ] & \leq 4 \e [ | W_{12} |^{1 + 2 \tau} | W_{13} W_{24} W_{34}  | ]  \nn \\
& \leq 4  ( \e [ | W_{12} |^{ 2  + 2 \tau } ] )^{ \frac { 2 + \tau  } { 1 +  \tau  } }.
\end{align*}
In consequence, it follows that 
\begin{align}
G_1 
& = \frac { B_X^2 } { 16 } \Big(  \e [ W_{12} W_{13} W_{24} W_{34}  ] - \e [ W_{12}^2 W_{13} W_{24} W_{34} ]   +  O(1)  ( \e [ | W_{12} |^{ 2  + 2 \tau } ] )^{ \frac { 2 + \tau  } { 1 +  \tau  } } \Big).  \label{g-np1}
\end{align}

As for term $ G_2 $, let $ D_2 =  \{  \max ( W_{12} , W_{ 13} , W_{24}, W_{45} ) \leq 1 \}   $ and $ D_2^{c} $ be its complement. Similarly, by \eqref{cru1} and \eqref{cru3} we can obtain 
\begin{align}
G_2  
& =  \frac { B_X^2 } { 64 } \Big(  4 \e [ W_{12} W_{13} W_{24} W_{45}  ]  -  \e [ W_{12}^2 W_{13} W_{24} W_{45} ] - \e [ W_{12} W_{13}^2 W_{24} W_{45}  ] \nn \\
& \hspace{1cm}  - \e [ W_{12} W_{13} W_{24}^2 W_{45} ] - \e [ W_{12} W_{13} W_{24} W_{45}^2 ] + O(1) ( \e [ | W_{12} |^{ 2  + 2 \tau } ] )^{ \frac { 2 + \tau  } { 1 +  \tau  } } \Big).   \label{g-np2}
\end{align}
We now consider term $ G_3^2 $. Define $ D_3 = \{ \max ( W_{12} , W_{13} )  \leq 1 \} $ and $ D_3^c  $ its complement. Similarly, we can show that 
\begin{align*}
G_3 &  =  B_X \Big( \frac 14 \e [ W_{12} W_{13} ] - \frac { 1 } { 8 } \e [ W_{12}^2 W_{13} \I \{ D_3 \} ]  + O(1)  \delta_2 \Big),
\end{align*}
where $ \delta_2 =  \e [ W_{12}^2  W_{13}^2  \I \{ D_3  \} ] + \e [ | W_{12}  W_{13}^3 | \I \{ D_3 \} ] +  \e [ | W_{12}  W_{13} |  \I \{ D_3^c \} ] $.
Note that when $ \e \Vert X \Vert ^{ 4 + 4 \tau  } < \infty $ for some $ 1/2 < \tau \leq 1 $, it follows from Chebyshev's inequality that 
\begin{align*}
\delta_2 \cdot  | \e [ W_{12} W_{13} ] | & \leq \e [ | W_{12} |^{ 1 + \tau  }  | W_{13} |^{1 + \tau} ]  \e [ | W_{12}  W_{13} |   ] + 3 \e [ | W_{12} | | W_{13} |^{ 1 + 2 \tau  } ]  \e [ | W_{12}  W_{13} |   ]  \nn \\
& \leq  4   ( \e [ | W_{12} |^{ 2 + 2 \tau  } ] )^{ \frac { 2 + \tau } { 1 + \tau  } },
\end{align*}
\begin{align*}
\delta_2  \e [ |  W_{12}^2  W_{13} | ]   \leq  4  \e [ | W_{12} | | W_{13} |^{2 \tau} ] \e [ | W_{12}^2 W_{13} | ]  \leq  4 ( \e [ | W_{12} |^{ 2 + 2 \tau  } ] )^{ \frac { 2 + \tau } { 1 + \tau  } }, 
\end{align*}
\begin{align*}
\delta_2^2 & \leq  3  \big( \e [ W_{12}^2  W_{13}^2  \I ( D_3  ) ]   \big)^2 + 3 \big(  \e [ | W_{12}  W_{13}^3 | \I ( D_3 ) ]    \big)^2 + 3 \big( \e [ | W_{12}  W_{13} |  \I ( D_3^c )]  \big)^2 \nn \\
& \leq  18 \big(  \e [ | W_{12} | | W_{13} |^{1 + \tau} ]  \big)^2   \leq  18  ( \e [ | W_{12} |^{ 2 + 2 \tau  } ] )^{ \frac { 2 + \tau } { 1 + \tau  } },
\end{align*}
\begin{align*}
\e [ | W_{12}  W_{13} | ] \e [ | W_{12}^2 W_{13} | \I ( D_3^c ) ]  & \leq    \e [ | W_{12}  W_{13} | ]  \big( \e [ | W_{12} |^2 | W_{13} |^{ 2 \tau } ]  +   \e [ | W_{12} |^{ 1 + 2 \tau } | W_{13} | ] \big) \nn \\
& \leq 2 ( \e [ | W_{12} |^{ 2 + 2 \tau  } ] )^{ \frac { 2 + \tau } { 1 + \tau  } },
\end{align*}
and
\begin{align*}
\big( \e [ W_{12}^2 W_{13} \I ( D_3 ) ] \big)^2 \leq  \big(  \e [ W_{12}^2 | W_{13} |^{\tau}  ]  \big)^2 \leq  ( \e [ | W_{12} |^{ 2 + 2 \tau  } ] )^{ \frac { 2 + \tau } { 1 + \tau  } }.
\end{align*}
Thus we can deduce 
\begin{align}
G_3^2 & =  \frac {  B_X^2 } { 16 }   \big( ( \e [ W_{12} W_{13} ] )^2 -   \e [ W_{12} W_{13}] \e [ W_{12}^2 W_{13} ]  + O( 1 )   ( \e [ | W_{12} |^{ 2 + 2 \tau  } ] )^{ \frac { 2 + \tau } { 1 + \tau  } }   \big).   \label{g-np3}
\end{align}

Then we deal with term $  \Delta G_4$. Denote by $ D_4 = \{  \max ( W_{12}, W_{13}, W_{24} ) \leq 1  \} $ and $ D_4^c  $ its complement. By \eqref{cru1} and \eqref{cru2}, we have for $ 1/2 < \tau \leq 1 $, 
\begin{align*}
G_4
& = B_X^{3/2} \e \big\{ [ ( 1 + W_{12} )^{1/2} - 1 ]  [ ( 1 + W_{13} )^{1/2} - 1 ] [ ( 1 + W_{24} )^{1/2} - 1 ] \big\} \nn \\
& = B_X^{3/2} \e \big\{ [   W_{12} / 2 + O( W_{12}^2 ) ]   [  W_{13} / 2 + O( W_{13}^2 ) ]   [  W_{24} / 2 + O( W_{24}^2 ) ] \I ( D_4 ) \big\}   \nn \\
& \quad + O(1) \e [ | W_{12} W_{13} W_{24} | \I ( D_4^c ) ] \nn \\
& = B_X^{3/2}  \Big( \frac 18 \e [ W_{12} W_{13} W_{24} ] + O(1) \big(  \e [ | W_{12} |^{ 2 \tau } | W_{13} W_{24} | ] + \e [ | W_{13} |^{ 2 \tau } | W_{12} W_{24} | ] \big) \Big).
\end{align*}
Moreover, it holds that 
\begin{align*}
\Delta &  =  B_X^{1/2}  \e [  ( 1 + W_{12} )^{1/2}  - 1 ] \nn \\
& =  B_X^{1/2}  \e \big[  \big(    W_{12} /2  -  W_{12}^2/ 8 + O(1) ( W_{12}^3 )  \big)  \I \{ W_{12} \leq 1 \} \big]  +  O(1) B_X^{1/2} \e [ | W_{12} |  \I \{ W_{12} > 1 \} ]  \nn \\
& = B_X^{1/2} \Big( - \frac 1 8 \e [ W_{12}^2 ] + O(1) \big( \e [ | W_{12} |^3 \I \{  W_{12} \leq 1 \} ] +  \e [ | W_{12} |^2 \I \{ W_{12} > 1 \} ] \big) \Big).  
\end{align*}
Observe from \eqref{g-p4-1} that for $ 1/2 <  \tau \leq 1 $, we have 
\begin{align*}
\e [ | W_{12} W_{13} W_{24} | ]   \big( \e [ | W_{12} |^3 \I \{  W_{12} \leq 1 \} ] +  \e [ | W_{12} |^2 \I \{ W_{12} > 1 \} ] \big)  &  \leq  C ( \e  [ W_{12}^2 ] )^{3/2}  \e [ | W_{12} |^{1 + 2 \tau} ]  \nn \\
&  \leq  C ( \e [ | W_{12} |^{ 2 + 2 \tau  } ] )^{ \frac { 2 + \tau  } { 1 + \tau  } }
\end{align*}
and 
\begin{align*}
& \big(  \e [ | W_{12} |^{ 2 \tau } | W_{13} W_{24} | ] + \e [ | W_{13} |^{ 2 \tau } | W_{12} W_{24} | ] \big)  \big( \e [ | W_{12} |^3 \I \{  W_{12} \leq 1 \} ] +  \e [ | W_{12} |^2 \I \{ W_{12} > 1 \} ] \big) \nn \\
& \leq  \big(  \e [ | W_{12} |^{ 2 \tau } | W_{13} W_{24} | ] + \e [ | W_{13} |^{ 2 \tau } | W_{12} W_{24} | ] \big)  \e [ W_{12}^2 ]  \nn\\
& \leq  2 ( \e [ | W_{12} |^{ 2 + 2 \tau  } ] )^{ \frac { 2 + \tau  } { 1 + \tau  } }.
\end{align*}
Hence it follows that 
\begin{align}
\Delta G_4
&  =  \frac { B_X^2  } { 64 } \big(  - \e [ W_{12}^2 ] \e [ W_{12} W_{13} W_{24} ]  + O( 1 ) ( \e [ | W_{12} |^{ 2 + 2 \tau  } ] )^{ \frac { 2 + \tau  } { 1 + \tau  } } \big).   \label{g-np4}
\end{align}

As for term $ \Delta^2 G_3 $, by \eqref{part3} and \eqref{cru1} we have for $  1/2 < \tau \leq 1 $, 
\begin{align}
| \Delta^2 G_3 |  & \leq  C B_X^2 \e [ | W_{12} |^{ 2 + 2 \tau } ] \e [ | W_{12} W_{13} | ] \nn \\
& \leq C B_X^2  ( \e [ | W_{12} |^{ 2 + 2 \tau } ] )^{ \frac { 2 + \tau  } { 1 + \tau } }. \label{g-np5}  
\end{align}
Note that \eqref{delta-square} entails that 
\begin{align}
| \Delta^4 |   & \leq  C B_X^2 \big( \e [ | W_{12} |^3 \I \{  W_{12}  \leq 1 \} ] + \e [ W_{12}^2 \I \{  W_{12} > 1  \} ]  \big)^2   \nn \\
& \leq C B_X^2 ( \e [ | W_{12} |^{ 2 + \tau } ] )^2 \leq C  B_X^2 ( \e [ | W_{12} |^{ 2 + 2 \tau } ] )^{ \frac { 2 + \tau } { 1 + \tau  } }. \label{g-np6}
\end{align}
Finally, substituting \eqref{g-np1}--\eqref{g-np6} into \eqref{g-expan} yields 
\begin{align}
\e [ g ( X_1, X_2, X_3, X_4 )  ]   = \frac { B_X^2 } {16}  \Big( E_1 + E_2  + O(1)  ( \e [ | W_{12} |^{ 2 + 2 \tau } ] )^{ \frac { 2 + \tau } { 1 + \tau  } } \Big),  \label{g-nre}
\end{align}
where
\begin{align}
E_1 & =  \e [ W_{12} W_{13} W_{24} W_{34}  ]-  4 \e [ W_{12} W_{13} W_{24} W_{45}  ]  + 2 ( \e [ W_{12} W_{13} ] )^2 \nn \\
\intertext{and}
E_2 & = - \e [ W_{12}^2 W_{13} W_{24} W_{34} ] +  \e [ W_{12}^2 W_{13} W_{24} W_{45} ] + \e [ W_{12} W_{13}^2 W_{24} W_{45}  ]  \nn \\
& \hspace{1cm}   +  \e [ W_{12} W_{13} W_{24}^2 W_{45} ]  + \e [ W_{12} W_{13} W_{24} W_{45}^2  ]   - 2  \e [ W_{12} W_{13}] \e [ W_{12}^2 W_{13} ]  \nn \\
&  \hspace{1cm}   -  \e [ W_{12}^2 ] \e [ W_{12} W_{13} W_{24} ]  . \label{E2_expression}
\end{align}

By some algebra, we can obtain Lemma \ref{le-g-2} in Section \ref{SecD.3}. Recall that $ B_X = 2 \e [ \Vert X \Vert^2  ]  $. Then the  equalities obtained above along with Lemma \ref{le-g-mo1} lead to
\begin{align*}
E_1 + E_2 & =  16 B_X^{ - 5 } \Big(  6 \e [ \Vert X \Vert^2  ] \e [ ( X_1^T \Sigma_x X_2 )^2 ] +  \e [ X^T X X^T  ] \Sigma_x^2 \e [ X X^T X ]  \nn \\
&\hspace{3cm}   + 2 \e [ ( X_1^T X_2 )^2 X_1^T \Sigma_x^2  X_2 ]  - 4 \e [  \Vert X_1 \Vert^2 ( X_1^T \Sigma_x X_2 )^2   ]  \Big).
\end{align*}
\ignore{
	\noi In addition, observing that
	\begin{align*}
	| \e [ ( X_1^T X_2 )^2 X_1^T \Sigma_x^2 X_2 ] | & = | \e [ ( X_1^T X_2 ) ( X_2^T X_1 X_1^T X_3 X_3^T X_4 X_4^T X_2 )  ] | \nn \\
	& = | \e [ ( X_1^T X_2 ) ( X_2^T \Sigma_x X_3 )^2  ] | \leq  ( \e [ ( X_1^T X_2 )^2 ] )^{1/2} ( \e [ ( X_1^T \Sigma_x X_2 )^4 ] )^{1/2} \nn \\
	& \leq  ( \e [  \Vert X_1 \Vert^2    \Vert X_2 \Vert^2   ] )^{1/2} ( \e [ ( X_1^T \Sigma_x X_2 )^4 ] )^{1/2} \nn \\
	&   =  \e [ \Vert X \Vert^2  ]  ( \e [ ( X_1^T \Sigma_x X_2 )^4 ] )^{1/2}  \leq ( \e [ \Vert X \Vert^4  ] )^{1/2}   ( \e [ ( X_1^T \Sigma_x X_2 )^4 ] )^{1/2},
	\end{align*}  
	\begin{align*}
	\e [ \Vert X \Vert^2 ] \e [ ( X_1^T \Sigma_x X_2 )^2 ] & \leq  ( \e [ \Vert X \Vert^4  ] )^{1/2}   ( \e [ ( X_1^T \Sigma_x X_2 )^4 ] )^{1/2}, \nn \\
	\intertext{and}
	\e [  \Vert X_1 \Vert^2 ( X_1^T \Sigma_x X_2 )^2   ] &  \leq  ( \e [ \Vert X \Vert^4  ] )^{1/2}   ( \e [ ( X_1^T \Sigma_x X_2 )^4 ] )^{1/2},
	\end{align*}
}
Therefore, we can obtain the desired result \eqref{g2}, which completes the proof of Proposition \ref{prop4}.

\subsection{Lemma \ref{le-g-2} and its proof} \label{SecD.3}

\begin{lemma}   \label{le-g-2} 
	It holds that 
	\begin{align}
	E_2 & = 16 B_X^{ - 5 }  \Big(   \e [ X^T X X^T  ] \Sigma_x^2 \e [ X X^T X ]   + 2 \e [ ( X_1^T X_2 )^2 X_1^T \Sigma_x^2  X_2 ]  \nn \\
	& \hspace{2cm}  - 4 \e [  \Vert X_1 \Vert^2 ( X_1^T \Sigma_x X_2 )^2   ]  + 4 \e [ \Vert X \Vert^2 ] \e [ ( X_1^T \Sigma_x X_2 )^2 ]  \Big). \label{E2}
	\end{align}
\end{lemma} 

\noindent \textit{Proof}. In view of the notation in the proofs of Lemmas \ref{le-v2} and \ref{le-g-mo1}, it holds that $ \alpha_1 ( X )  = \Vert X\Vert^2 - \e [ \Vert X \Vert^2  ] $ and $ \alpha_2 (X_1, X_2)  = X_1^T X_2 $. Thus we have 
\begin{align*}
W_{12}  = B_X^{ - 1 } [ \alpha_1 (X_1) + \alpha_1 (X_2) - 2 \alpha_2 (X_1, X_2) ] .
\end{align*}
Then it follows that 
\begin{align*}
& \e [ W_{12}^2 W_{13} W_{24} W_{34} ] \nn \\
& = B_X^{ - 5 }   \e \Big\{  \big[ \alpha_1 ( X_1 ) + \alpha_1 ( X_2 ) - 2 \alpha_2 ( X_1, X_2 ) \big]^2 \big[  \alpha_1 ( X_1 ) + \alpha_1 ( X_3 ) - 2 \alpha_2 ( X_1, X_3 ) \big] \nn \\
& \hspace{2cm} \times  \big[ \alpha_1 ( X_2 ) + \alpha_1 ( X_4 ) - 2 \alpha_2 ( X_2, X_4 ) \big] \big[  \alpha_1 ( X_3 ) + \alpha_1 ( X_4 ) - 2 \alpha_2 ( X_3, X_4 ) \big] \Big\}.
\end{align*}
The idea of the proof is to expand the products. Since $ X_1, X_2, X_3 $, and $ X_4 $ are i.i.d., we can deduce 
\begin{align*}
&  \e [ W_{12}^2 W_{13} W_{24} W_{34} ]  \nn \\
& =   B_X^{ - 5 } (2 D_1  + 8  D_2  
- 20 D_3    - 16  D_4   - 8 D_5   + 24 D_6   + 32 D_7 + 16 D_8   - 48 D_9  - 32 D_{10}    + 64 D_{11} ) ,
\end{align*}
where 
\begin{align*}
D_1 & = \e [ \alpha_1^2 ( X ) ] \e [ \alpha_1^3 ( X ) ], \\
D_2 & = \e [ \alpha_1^2  ( X ) ] \e [ \alpha_2^2 ( X_1, X_2 ) \alpha_1 ( X_2 ) ],  \nn \\  
D_3 &= \e [ \alpha_1^2 (X) ] \e [ \alpha_1 ( X_1 ) \alpha_1 ( X_2 ) \alpha_2 ( X_1, X_2 ) ], \nn \\
D_4 & =  \e [ \alpha_1 ( X_3 ) \alpha_2 ( X_1, X_3 ) \alpha_2^2 ( X_1, X_2 ) \alpha_1 ( X_2 ) ],  \nn \\
D_5 &=  \e [ \alpha_2^2 (X_1, X_2 ) ]  \e [ \alpha_1 ( X_1 ) \alpha_1 ( X_2 ) \alpha_2 ( X_1, X_2 ) ],  \nn \\
D_6 &=   \e  [ \alpha_1 (X_3 ) \alpha_2 ( X_1, X_3 ) \alpha_2 ( X_1, X_2 ) \alpha_1^2 ( X_2 ) ], \nn \\
D_7 & = \e [ \alpha_2 ( X_3, X_4 ) \alpha_2 ( X_1, X_3 )  \alpha_1 ( X_4 ) \alpha_2^2 ( X_1, X_2 ) ],   \nn \\
D_8  & = \e [\alpha_1 ( X_1 ) \alpha_1 ( X_2 )  \alpha_1 ( X_3 ) \alpha_2 ( X_1, X_3 ) \alpha_2 ( X_1, X_2 ) ],  \nn \\
D_9 & =  \e [  \alpha_2 ( X_3, X_4 )  \alpha_1 ( X_4 )  \alpha_2 ( X_1, X_3 )  \alpha_2 ( X_1, X_2 ) \alpha_1 ( X_2 ) ], \nn \\
D_{10}  & =\e [ \alpha_2 ( X_3, X_4 )  \alpha_2 ( X_2, X_4 ) \alpha_2 ( X_1, X_3 ) \alpha_2^2 ( X_1, X_2 ) ], \nn  \\
D_{11}  &  = \e [ \alpha_2 ( X_3, X_4 ) \alpha_2 ( X_1, X_3 ) \alpha_2 ( X_2, X_4 ) \alpha_1 ( X_1 )  \alpha_2 ( X_1, X_2 ) ] \Big] .
\end{align*}

Similarly, we can show that 
\begin{align*}
\e [ W_{12}^2 W_{13} W_{24} W_{45}  ]  
& =   B_X^{ - 5 } ( D_1  + 4   D_2  - 8 D_3 - 8 D_4 + 8 D_6    +  8 D_8 ) ,   \\
\e [ W_{12} W_{13}^2 W_{24} W_{45}  ] 
& =   B_X^{ - 5 }  ( D_1 +  4 D_2    - 8  D_3 - 8 D_5 + 4     D_6 + 16 D_7 - 16   D_9  ) ,  \\
\e [ W_{12} W_{13} W_{24}^2  W_{45}  ]  
& = B_X^{  - 5 } (  D_1 +   4  D_2 - 8  D_3  - 8  D_4  + 8 D_6  + 8 D_8   ) ,  \\
\e [ W_{12} W_{13} W_{24} W_{45}^2  ]  
& =   B_X^{ - 5 } (   D_1 + 4   D_2 - 8 D_3 - 8 D_5  + 4  D_6 + 16  D_7    - 16  D_9 ), \\
\e [ W_{12} W_{13} ] \e [ W_{12}^2 W_{13} ]  & =  B_X^{ - 5 } (  D_1 + 4  D_2 - 4  D_3  ), \\
\e [ W_{12}^2 ] \e [ W_{12} W_{13} W_{24} ]  & =  - B_X^{-5} (  4 D_3 + 8 D_5 ). 
\end{align*}
Thus by plugging the above equalities into \eqref{E2_expression}, it holds that 
\begin{align*}
E_2 & =  16 B_X^{ - 5 } (  D_9 +  2 D_{10} - 4 D_{11} ).
\end{align*}
It is easy to see that 
\begin{align*}
D_{9} & =  \e [ X^T X X^T  ] \Sigma_x^2 \e [ X X^T X ],  \\
D_{10} &  =  \e [ ( X_1^T X_2 )^2 X_1^T \Sigma_x^2  X_2 ],   \\
D_{11} & = \e [  \Vert X_1 \Vert^2 ( X_1^T \Sigma_x X_2 )^2   ]  -  \e [ \Vert X \Vert^2 ] \e [ ( X_1^T \Sigma_x X_2 )^2 ].
\end{align*}
Therefore, we can obtain the desired result \eqref{E2}, which concludes the proof of Lemma \ref{le-g-2}.

\renewcommand{\thesubsection}{E.\arabic{subsection}}
\section{Connections between normal approximation and gamma approximation} \label{normal-gamma} 

For the test of independence based on the sample distance covariance, empirically one can use the gamma approximation to calculate the limiting p-values. \cite{HS2016} showed that under some moment conditions and the independence of $X$ and $Y$, it holds that 
\begin{align*}
n \V_n^* (X, Y)  \xrightarrow  [ n\rightarrow \infty ]  {\mathscr{D}} \sum_{i = 1} ^{\infty} \lambda_i (Z_i^2 - 1),
\end{align*}
where $ \{\lambda_i \}_{i \geq 1 }$ are some values depending on the underlying distribution and $ \{Z_i\}_{i \geq 1} $ are i.i.d. standard normal random variables. 
In practice, it is infeasible to apply this limiting distribution directly and thus the gamma approximation can serve as a surrogate. By \cite{HH2017}, it follows that 
\[ \sum_{i=1}^{\infty} \lambda_i  = \e [ \Vert X - X' \Vert ]  \e [ \Vert  Y - Y'\Vert ] \ \text{ and } \ \sum_{i = 1}^{\infty} \lambda_i^2 = \V^2 (X) \V^2(Y), \] 
and hence $ \sum_{i = 1}^{\infty} \lambda_i (Z_i^2 - 1) $ can be approximated by a centered gamma distribution $ \Gamma(\beta_1, \beta_2) - \beta_1 \beta_2^{ - 1 }  $, where the shape and rate parameters $ \beta_1 $ and $ \beta_2 $ are determined by matching the first two moments. To this end, we define 
\begin{align*}
\beta_1 = \frac { \big( \sum_{i = 1}^{\infty} \lambda_i \big) ^2 } { 2 \sum_{i = 1}^{\infty}  \lambda_i^2}  = \frac { \big( \e [ \Vert  X - X' \Vert ]  \e [ \Vert Y - Y' \Vert ]   \big)^2  } { 2 \V^2 (X) \V^2 (Y) }
\end{align*}
and 
\begin{align*}
\beta_2 = \frac { \sum_{i = 1}^{\infty}  \lambda_i } {  2 \sum_{i = 1}^{\infty} \lambda_i^2  } = \frac { \e [ \Vert  X - X' \Vert ]  \e [ \Vert Y - Y' \Vert ]  } { 2 \V^2 (X) \V^2 (Y) }.
\end{align*}

For a simple illustration, let us consider a specific case when both $X $ and $Y$ consist of i.i.d. components.  Then it holds that $ \e [ \Vert X - X' \Vert  ] = O (\sqrt p) $ and $ \e [ \Vert Y - Y' \Vert  ] = O (\sqrt {p}) $. Moreover, it follows from Proposition \ref{prop2} that $ \V^2 (X) $ and $ \V^2 (Y) $ are bounded from above and below by some positive constants, which entails that $ \beta_1 = O (p q) $ and $ \beta_1 \rightarrow \infty $ as $ \max \{p, q\} \to \infty $. Recall the fact that the gamma random variable can be represented as a sum of certain i.i.d. exponential random variables. Thus by the central limit theorem, we have 
\begin{align*}
\frac { \Gamma(\beta_1, \beta_2) - \beta_1 \beta_2^{ - 1}    } { \sqrt{  \beta_1 \beta_2^{- 2} } } \stackrel{\mathscr{D}}{ \rightarrow } N (0, 1)
\end{align*}
as $ \max \{ p, q \}  \to \infty$. Since $ \beta_1 \beta_2^{-2} = 2 V^2 (X) V^2 (Y) $ and Lemma \ref{le-consist} has provided the consistency of $ \V^*_n (X) $ and $ \V^*_n (Y) $, it holds that 
\begin{align*}
\mathbbm{P} ( T_n \leq  x  ) & = \mathbbm{P} \Big( \sqrt{\frac { n (n - 1) } {2}}  \frac { \V^*_n (X, Y) } { \sqrt{ \V^* (X)  \V^* (Y) } } \leq x  \Big)  \\
& \approx \mathbbm{P}  \Big( \frac {n \V_n^* (X, Y) }  {\sqrt{ 2 \V^2 (X) \V^2 (Y) } }  \leq x \Big) \approx \mathbbm{P} \Big(  \frac { \Gamma(\beta_1, \beta_2) - \beta_1 \beta_2^{ - 1}    } { \sqrt{  \beta_1 \beta_2^{- 2} } }  \leq x \Big) \\
& \to \Phi (x),
\end{align*}
where $ \Phi (x) $ stands for the standard normal distribution function. Therefore, the gamma approximation for $ n \V^* (X, Y) $ may be asymptotically equivalent to the normal approximation to $ T_n $ under certain scenarios. It is worth mentioning that the above analysis intends to build some connections between the normal approximation and the gamma approximation, but is not a rigorous proof. A rigorous theoretical foundation for the gamma approximation still remains undeveloped.

\section{Asymptotic normality of $T_R$} \label{SecF}
An anonymous referee asked a great question on whether similar asymptotic normality as in Theorem \ref{thm1} and associated rates of convergence as in Theorem \ref{thm3} hold for the studentized sample distance correlation $T_R$. The answer is affirmative as shown in the following proposition.

\begin{proposition} \label{prop-tr}
	Under the same conditions of Theorem \ref{thm1}, we  have $T_R \stackrel{\mathscr{D}}{\rightarrow} N(0, 1)$. Moreover, under the conditions of Theorem \ref{thm3}, the same rate of convergence as in \eqref{rate-general} holds for $T_R$.
\end{proposition}

\noindent \textit{Proof}. By Lemma \ref{le-consist}, we have $ \V_n^* (X) / \V^2 (X) \stackrel{p}{\rightarrow} 1 $ and $ \V_n^* (Y) / \V^2 (Y) \stackrel{p}{\rightarrow} 1  $ under condition \eqref{cond1}. In addition, it follows from \eqref{h-b1} and Lemma \ref{h-d} that for $ 0 < \tau \leq 1 $,
\begin{align*}
\e [ | \V_n^* (X, Y)  - \V^2 (X, Y)|^{1 + \tau } ] \leq C n^{ - \tau } \big[ \e (| d (X_1, X_2) |^{ 2 + 2 \tau } )  \e ( | d (Y_1, Y_2) |^{2 + 2 \tau } )  \big] ^{ 1/2 }.
\end{align*}
Hence under condition (18), it holds that 
\begin{align*}
& \e \Big[ \Big| \frac { \V^*_n (X, Y) } { \sqrt{\V^2 (X) \V^2 (Y)} } - \mathcal{R}^2 (X, Y) \Big|^{1 + \tau} \Big] \\
&\quad  \leq  \frac { C } { n^{\tau /2} } \Big( \frac {  \big[ \e (| d (X_1, X_2) |^{ 2 + 2 \tau } )  \e ( | d (Y_1, Y_2) |^{2 + 2 \tau } )  \big]   } { n^{\tau } [ \V^2 (X) \V^2 (Y)]^{1 + \tau} } \Big)^{1/2}  \rightarrow 0.
\end{align*}
This entails that $ \frac { \V^*_n (X, Y) } { \sqrt{V^2 (X) \V^2 (Y)} } \stackrel{p}{\rightarrow} \mathcal{R}^2 (X, Y) $ and thus $ \mathcal{R}_n^* (X, Y)  \stackrel{p}{\rightarrow } \mathcal{R}^2 (X, Y)$ as well. Under the null hypothesis, it holds that $ \mathcal{R}^2 (X, Y)  = 0 $ and hence $ \mathcal{R}_n ^* (X, Y) \stackrel{p}{\rightarrow } 0$. In light of the definition of $T_n^*$, it holds that 
$$
T_R = T_n \cdot \sqrt{\frac { n (n - 3) - 2 } { n (n - 1)   } } \cdot \frac { 1 } { \sqrt{ 1 - \mathcal{R}_n^* (X, Y)} }.
$$
By Theorem \ref{thm1}, we have $ T_n \stackrel{p}{\rightarrow} N (0, 1) $. As a consequence, under the conditions of Theorem \ref{thm1}, it holds that $ T_R \stackrel{\mathscr{D}}{\rightarrow} N (0, 1)$ as well. 

Next we proceed to show that the rates of convergence in Theorem \ref{thm3} also apply to $T_R$. It follows from the definitions of $T_n $ and $T_R$ that for $x > 0$ (similar analysis applies for $x \leq 0$), 
\begin{align*}
\mathbbm{P} (T_R >   x ) & =  \mathbbm{P} \bigg( T_n > x \cdot \sqrt { \frac { n (n - 1) } { 2 x^2 + n (n - 3) - 2} } \bigg) .
\end{align*}
Thus it holds that for $x > 0 $, 
\begin{align*}
& | \mathbbm{P} ( T_R > x)  - [ 1 - \Phi(x) ] | \\
& \leq  \Bigg|  \mathbbm{P} \bigg( T_n > x \cdot \sqrt{ \frac { n (n - 1) } { 2 x^2 + n (n - 3) - 2 } } \bigg) - \bigg[ 1 - \Phi \bigg( x \cdot \sqrt{ \frac { n (n - 1) } { 2 x^2 + n (n - 3) - 2 } }  \bigg) \bigg] \Bigg| \\
& \quad + \Bigg| \bigg[ 1 - \Phi \bigg( x \cdot \sqrt{ \frac { n (n - 1) } { 2 x^2 + n (n - 3) - 2 } }  \bigg) \bigg] - [ 1 - \Phi (x) ] \Bigg| .
\end{align*}
Note that the first term on the right hand side of the above inequality is bounded by the convergence rate in Theorem 3. As for the second term, observe that when $ 0 < x \leq c n $ for some small constant $ c > 0 $, we have 
\begin{align*}
& \bigg|  x \cdot \sqrt{ \frac { n (n - 1) } { 2x^2 + n (n - 3) - 2 } } - x \bigg|\\
& = \frac {x} { 1 + \sqrt{ \frac { n (n - 1) } { 2 x^2 + n (n - 3) - 2   } }  } \cdot \bigg| \frac { 2 x^2 - 2 n - 2  } { 2 x^2 + n (n - 3) - 2 } \bigg|\\
& \leq x \cdot \bigg| \frac { 2 x^2 - 2 n - 2  } { 2 x^2 + n (n - 3) - 2 } \bigg| = O \Big\{ x \big( \frac { x^2 } {n^2 } + \frac { 1 } {n} \big) \Big\}.
\end{align*}
By the properties of normal distribution function, we can obtain that for $ 0 < x \leq c n  $, 
\begin{align}
\Bigg| \bigg[ 1 - \Phi \bigg( x \cdot \sqrt{ \frac { n (n - 1) } { 2 x^2 + n (n - 3) - 2 } }  \bigg) \bigg] - [ 1 - \Phi (x) ] \Bigg|  = O\big(\frac 1 n \big). \label{bound-norm}
\end{align}

When $ x > c n $,  it is easy to see that $ 1 - \Phi (x) \leq e^{ - x^2 / 2 } \leq C_1 n^{-1} $ for some constant $C > 0$ depending on $c$. In addition, it holds that 
\begin{align*}
x \cdot \sqrt{ \frac { n (n - 1) } { 2 x^2 + n (n - 3) - 2 } }  =   \sqrt{ \frac {    n (n - 1) } { 2   + n (n - 3) / x^2 - 2/x^2  } }  \geq C_2 n ,
\end{align*}
where $C_2 > 0$ is some constant depending on $c$. Then it follows that for some positive constant $C_3$ depending on $c$, 
\begin{align*}
1 - \Phi \bigg( x \cdot \sqrt{ \frac { n (n - 1) } { 2 x^2 + n (n - 3) - 2 } }  \bigg)  \leq C_3 n^{-1}.
\end{align*}
Thus \eqref{bound-norm} still holds for the case of $x > c n $. 
In view of the convergence rate in Theorem 3, it is easy to see that $ O ( \frac {1} {n} ) $ is of a smaller order. Finally, we obtain that the same convergence rate as stated in Theorem 3 also applies to $T_R$, which completes the proof of Proposition \ref{prop-tr}.

\end{document}